\documentclass[reqno,12pt]{amsart}
\theoremstyle{plain}    
\newtheorem{thm}{Theorem}[section]
\theoremstyle{plain}    
\newtheorem{cor}[thm]{Corollary} 
\newtheorem{lemma}[thm]{Lemma} 
\newtheorem{prop}[thm]{Proposition}
\newtheorem{defi}[thm]{Definition}
\newtheorem{algorithm}[thm]{Algorithm}
\newtheorem{example}[thm]{Example}

\newtheorem{remark}[thm]{Remark}
\newtheorem{conj}[thm]{Conjecture}

\usepackage{amscd,amssymb,comment,epic,eufrak,euscript,graphics}
\usepackage[arrow,matrix]{xy}

\newcommand\Ac{{\mathcal{A}}}

\newcommand\Afr{{\mathfrak A}}

\newcommand\ah{{\hat a}}

\newcommand\at{{\tilde a}}

\newcommand\Bc{{\mathcal{B}}}

\newcommand\Bfr{{\mathfrak B}}

\newcommand\bh{{\hat b}}

\newcommand\bt{{\tilde b}}

\newcommand\btimes{\displaystyle\operatornamewithlimits\times}

\newcommand\Cpx{{\mathbf C}}

\newcommand\Dc{{\mathcal{D}}}

\newcommand\diag{\text{\rm diag}}

\newcommand\dif{\mbox{\it d}}

\newcommand\DT{\operatorname{DT}}

\newcommand\Eb{{\mathbf E}}

\newcommand\Ec{{\mathcal{E}}}

\newcommand\eps{\epsilon}

\newcommand\Fb{{\mathbf F}}

\newcommand\ft{{\tilde f}}

\newcommand\GRM{{\text{\rm GRM}}}

\newcommand\HEu{{\EuScript H}}                   

\newcommand\HURM{{\operatorname{HURM}}}

\newcommand\ImagPart{{\mathrm{Im}\;}}

\newcommand\intr{{\operatorname{int}}}

\newcommand\Ints{{\mathbf Z}}

\newcommand\iotah{\hat{\iota}}

\newcommand\lambdab{\overline\lambda}

\newcommand\Lambdao{{\Lambda\oup}}

\newcommand\lambdat{{\tilde\lambda}}

\newcommand\LEu{{\EuScript L}}                   

\newcommand\Mcal{{\mathcal{M}}} 

\newcommand\MEu{{\EuScript M}}                   

\newcommand\mubar{{\overline\mu}}

\newcommand\Nats{{\mathbf N}}

\newcommand\Nc{{\mathcal{N}}}

\newcommand\NC{\operatorname{NC}}

\newcommand\NCP{\operatorname{NCP}}

\newcommand\NTO{\operatorname{NTO}}

\newcommand\nut{{\tilde{\nu}}}

\newcommand\oneb{{\mathbf1}}

\newcommand\oneh{{\hat 1}}

\newcommand\oup{^{\mathrm o}}

\newcommand\Prob{{\mathrm{Prob}}}

\newcommand\pt{{\tilde p}}

\newcommand\Qc{{\mathcal{Q}}}

\newcommand\RealPart{{\mathrm{Re}\;}}

\newcommand\Reals{{\mathbf R}}

\newcommand\Res{{\operatorname{Res}}}

\newcommand\restrict{{\upharpoonright}}

\newcommand\SGRM{{\text{\rm SGRM}}}

\newcommand\sign{{\text{\rm sign}}}

\newcommand\supp{\operatorname{supp}}

\newcommand\Tc{{\mathcal{T}}}

\newcommand\Tcirc{{\mathbf T}}

\newcommand\Th{{\widehat T}}

\newcommand\tr{{\mathrm{tr}}}

\newcommand\Tr{{\mathrm{Tr}}}

\newcommand\Tt{{\widetilde T}}

\newcommand\Uc{{\mathcal{U}}}

\newcommand\UT{{\EuScript U\EuScript T}}

\newcommand\UTGRM{{\text{\rm UTGRM}}}

\newcommand\XEu{{\EuScript X}}

\newcommand\Xh{{\widehat X}}

\newcommand\Xt{{\widetilde X}}

\newcommand\YEu{{\EuScript Y}}

\newcommand\Yt{{\widetilde Y}}

\newcommand\zbar{{\overline z}}

\newcommand\zerob{{\mathbf0}}

\begin{document}

\pagestyle{myheadings}

\title{DT--operators and decomposability of Voiculescu's circular operator}
 
\author{Ken Dykema}
\address{\hskip-\parindent
Ken Dykema \\
Department of Mathamtics \\
Texas A\&M University \\
College Station TX 77843--3368, USA}
\email{Ken.Dykema@math.tamu.edu}

\author{Uffe Haagerup}
\address{\hskip-\parindent
Uffe Haagerup \\
Department of Mathematics and Computer Science \\
University of Southern Denmark --- Odense University \\
Campusvej 55 \\
5230 Odense M, Denmark}
\email{haagerup@imada.sdu.dk}

\thanks{K.D.\ supported in part by NSF grant DMS--0070558.
U.H.\ affiliated with MaPhySto,
Centre for Mathematical Physics and
Stochastics, which is funded by a grant from The Danish National Research Foundation.
Both authors thank also the Mathematical Sciences Research Institute in Berkeley,
where part of this work was carried out in the spring of 2001.
U.H. was employed by the Clay Mathematical Institute as a Clay Research Scholar while
he was at MSRI.
Research at MSRI is supported in part by NSF grant DMS--9701755.}

\date{8 May, 2002}

\begin{abstract}
The DT--operators are introduced, one for every pair $(\mu,c)$ consisting of
a compactly supported Borel
probability measure $\mu$ on the complex plane and a constant $c>0$.
These are operators on Hilbert space that
are defined as limits in $*$--moments of certain upper triangular random matrices.
The DT--operators include Voiculescu's circular operator
and elliptic deformations of it, as well as the circular free Poisson operators.
We show that every DT--operator is strongly decomposable.
We also show that a DT--operator generates a II$_1$--factor,
whose isomorphism class depends only on the number and sizes of atoms of $\mu$.
Those DT--operators that are also R--diagonal are identified.
For a quasi--nilpotent DT--operator $T$, we find the distribution of $T^*T$
and a recursion formula for general $*$--moments of $T$.
\end{abstract}

\maketitle


\markboth{\tiny DT--operators}{\tiny DT--operators}

\tableofcontents

\newpage

\section{Introduction}

\subsection{Local spectral theory and decomposability}

Let us begin by describing some ideas and results in local spectral theory related to
decomposability of operators.
Our exposition is drawn from the book~\cite{LN} by Laursen and Neumann,
and depends on work of Foia\c s \cite{F63}, \cite{F68}, Apostol \cite{A68a}, \cite{A68b},
Albrecht~\cite{Al}, Jafarian and Vasilescu~\cite{JV} and others
(see~\cite{LN} for some history, more detailed citations and proofs).
Theorem~\ref{thm:dec} below and the considerations surrounding it apply
to bounded operators on Banach spaces, but in keeping with the subject
of this paper we will restrict to discussion of a bounded operator $T$ on
a Hilbert space $\HEu$.

Such an operator $T$ is said to be {\em decomposable}
if, for every cover $\Cpx=U\cup V$ of the complex plane by two open subsets $U$ and $V$,
there are $T$--invariant closed subspaces $\HEu'$ and $\HEu''$ of $\HEu$ such that
the spectra of the restrictions of $T$ satisfy
$\sigma(T\restrict_{\HEu'})\subseteq U$ and $\sigma(T\restrict_{\HEu''})\subseteq V$,
and such that $\HEu=\HEu'+\HEu''$.

Given a bounded operator $T$ on $\HEu$, a {\em spectral capacity for} $T$
is a mapping $E$ from the set of all closed subsets of $\Cpx$ into the set of all
closed $T$--invariant subspaces of $\HEu$, such that
\renewcommand{\labelenumi}{(\roman{enumi})}
\begin{enumerate}

\item $E(\emptyset)=\{0\}$ and $E(\HEu)=\Cpx$,
\vskip1ex

\item $\HEu=E(\overline{U_1})+\dots+E(\overline{U_n})=\HEu$
for every open cover $\{U_1,\dots,U_n\}$ of $\Cpx$,
\vskip1ex

\item $E(\bigcap_{n=1}^\infty F_n)=\bigcap_{n=1}^\infty E(F_n)$
for any closed subsets $F_1,F_2,\dots$ of $\Cpx$,
\vskip1ex

\item $\sigma(T\restrict_{E(F)})\subseteq F$ for every closed subset $F$ of $\Cpx$,
(with the convention that the operator on the Hilbert space $\{0\}$ has empty spectrum).

\end{enumerate}

If $\xi\in\HEu$, the {\em local resolvent set}, $\rho_T(\xi)$, of $T$ at $\xi$ is the union
of all open subsets $U$ of $\Cpx$ for which there exist holomorphic vector--valued
functions $f_U:U\to\HEu$ such that $(T-\lambda)f_U(\lambda)=\xi$ for all $\lambda\in U$.
The {\em local spectrum} of $T$ at $\xi$ is then defined to be $\sigma_T(\xi)=\Cpx\backslash\rho_T(\xi)$.
For any subset $A$ of $\Cpx$, the corresponding {\em local spectral subspace} of $T$ is
\[
\HEu_T(A)=\{\xi\in\HEu\mid\sigma_T(\xi)\subseteq A\}.
\]
It is clear that $\HEu_T(A)$ is $T$--invariant, and even {\em $T$--hyperinvariant},
namely invariant under every operator that commutes with $T$.
Hence letting $p_T(A)$ deonte the projection of $\HEu$ onto the closure
of $\HEu_T(A)$, it follows that $p_T(A)$ lies in the von Neumann algebra $\operatorname{vN}(T)$
that is generated by $T$.
In fact, by using standard techniques like those in the proofs of Lemmas 2.3 and 2.4 of~\cite{DH},
one can show that the projection $p_T(A)$, as an element of the abstract
W$^*$--algebra $\Mcal$ that is isomorphic to $\operatorname{vN}(T)$,
is independent of the particular representation of $\Mcal$ on a Hilbert space.
(Our convention is that a von Neumann algebra is the image of a normal representation of
an abstract W$^*$--algebra on a Hilbert space.)

\begin{thm}{\rm (cf.~\cite[\S1.2]{LN}).}
\label{thm:dec}
Let $T$ be a bounded operator on a Hilbert space $\HEu$.
Then the following are equivalent:
\renewcommand{\labelenumi}{(\roman{enumi})}
\begin{enumerate}
\item $T$ is decomposable,
\item $T$ has a spectral capacity,
\item
for every closed subset $F$ of $\Cpx$, $\HEu_T(F)$ is closed and
\[
\sigma((1-p_T(F))T)\subseteq\overline{\sigma(T)\backslash F}\;,
\]
where the spectrum is of $(1-p_T(F))T$ considered as an operator on $\HEu_T(F)^\perp$.
\end{enumerate}
Moreover, if $T$ is decomposable then the map $F\mapsto\HEu_T(F)$ is the
unique spectral capacity for $T$.
\end{thm}

An operator $T$
is said to be {\em strongly decomposable} if it is decomposable and, morover, if
the restriction
$T\restrict_{\HEu_T(F)}$ is decomposible for every closed subset $F$ of $\Cpx$.
In light of the above comments about $p_T(A)$ in abstract W$^*$--algebras, it is clear that
decomposability and strong decomposability
of operators can be thought of as algebraic properties
of elements of W$^*$--algebras.

\subsection{DT--operators}

In this paper, we introduce the DT--operators.
These are operators on Hilbert
space that are defined as limits of certain large random matrices and are intimately related
to free probability theory, (see the book~\cite{VDN}).
We will introduce them first in the context of $*$--noncommutative probability spaces.
A {\em $*$--noncommutative probability space} is a pair
$(\Afr,\tau)$, where $\Afr$ is a unital algebra over $\Cpx$ with an involution usually denoted
$a\mapsto a^*$, and where $\tau:A\to\Cpx$ is a linear functional that is unital ($\tau(1)=1$)
and positive ($\tau(a^*a)\ge0$).
A $*$--noncommutative probability space is said to be {\em tracial} if $\tau$ satisfies the trace property.
An important example is the tracial $*$--noncommutative probability space of $n\times n$ random matrices, $(\MEu_n,\tau_n)$,
where $\MEu_n$ is an algebra of $n\times n$ matrices whose entries are random variables
over a classical probability space having moments of all orders and where,
viewing an element $a\in\MEu_n$ as an $M_n(\Cpx)$--valued random variable,
$\tau_n(a)$ is the expectation of $\tr_n(a)$, where $\tr_n$ is the normalized
trace on $M_n(\Cpx)$.
Let $T_n\in\MEu_n$ be a strictly upper triangular random matrix, the real and imaginary parts
of whose entries above the diagonal form a family
of $n(n-1)$ i.i.d.\ gaussian random variables, each having mean $0$
and variance $1/2n$.
Fix a Borel measure $\mu$ on $\Cpx$ having compact support, and let $D_n\in\MEu_n$ be a diagonal
random matrix having mutually independent $\mu$--distributed diagonal entries and
such that, as matrix--valued random variables, $D_n$ and $T_n$ are independent.
In~\S\ref{sec:DTdist}, we show that the pair $D_n,\,T_n$ converges in joint $*$--moments 
as $n\to\infty$.
This means that $\tau_n$ of any fixed word in $D_n$, $D_n^*$, $T_n$ and $T_n^*$ converges as $n\to\infty$;
in fact, we have a combinatorial formula involving non--crossing pairings giving the limiting value.
In particular, taking $Z_n=D_n+cT_n$ for any fixed $c>0$, it follows that $Z_n$ converges in $*$--moments
as $n\to\infty$.
This means that
for every $k\in\Nats$ and all $\eps(1),\ldots,\eps(k)\in\{*,1\}$, the limit
\begin{equation}
\label{eq:Znlim}
\lim_{n\to\infty}\tau_n(Z_n^{\eps(1)}Z_n^{\eps(2)}\cdots Z_n^{\eps(k)})
\end{equation}
exists.
We define a $\DT(\mu,c)${\em --element} to be an element $Z\in\Afr$ of a $*$--noncommutative probability space $(\Afr,\tau)$
whose $*$--moments $\tau(Z^{\eps(1)}\cdots Z^{\eps(k)})$
are given by the limits~\eqref{eq:Znlim}.
An element is a {\em DT--element} if it is a $\DT(\mu,c)$--element for some $\mu$ and $c$.

Drawing on results from~\cite{DH}, we have (Theorem~\ref{thm:DTflexdiag}) that the diagonal part
$D_n$ of $Z_n=D_n+cT_n$ can be modified quite extensively from the situation described above without changing
that $Z_n$ converges in $*$--moments to the $\DT(\mu,c)$--element $Z$.
This flexibility in the random matrix model for $Z$ is important in our investigation of DT--elements.

In~\S\ref{sec:DTex}, we will show that
the class of DT--elements includes some elements arising naturally in free probability theory,
namely Voiculsecu's circular element, the circular free Poisson elements, and
the so called elliptic elements, which are deformations of the circular element.

In~\S\ref{sec:DTinLF2}, we construct DT--elements as operators on Hilbert space;
more specifically, for every $\mu$ and $c$, we give a canonical construction of a $\DT(\mu,c)$--element
in the $*$--noncommutative probability space $(L(\Fb_2),\tau)$, where $L(\Fb_2)\subseteq B(\ell^2(\Fb_2))$ is the von Neumann algebra
generated by the left regular representation of the nonabelian free group on two generators and $\tau$
is the canonical tracial state.
A {\em W$^*$--noncommutative probability space} is a $*$--noncommutative probability space $(\Mcal,\tau)$,
where $\Mcal$ is a von Neumann algebra and $\tau$ is a normal state.
We say that a $\DT(\mu,c)${\em --operator} is a $\DT(\mu,c)$--element in W$^*$--noncommutative probability space $(\Mcal,\tau)$,
where the state $\tau$ is assumed to be faithful.
Of course, an element of a von Neumann algebra with a specified faithful normal state is called a {\em DT--operator}
if it is a $\DT(\mu,c)$--operator for some $\mu$ and $c$.

Using the flexibility in the random matrix model for a DT--element and our construction of DT--operators in $L(\Fb_2)$,
we prove (Theorem~\ref{thm:DTfreematrix}) that if
\begin{equation}
\label{eq:DTmat}
A=\left(\begin{matrix}
   a_1 & b_{12} & \cdots & b_{1N}    \\
     0 &  a_2   & \ddots & \vdots    \\
\vdots & \ddots & \ddots & b_{N-1,N} \\
     0 & \cdots &    0   & a_N
\end{matrix}\right)
\end{equation}
is an upper triangular $N\times N$
matrix of operators,
whose entries are mutually $*$--free,
where each $b_{ij}$ is a circular element
and where $a_j$ is a $\DT(\mu_j,\frac c{\sqrt N})$--operator,
then $A$ is itself a $\DT(\mu,c)$--operator, where $\mu$ is the average
of $\mu_1,\ldots,\mu_N$.
This result appears as a basic tool in our proofs of a number of interesting facts about
DT--operators.

In~\S\ref{sec:decomp}, we use the upper triangular picture~\eqref{eq:DTmat}
to investigate the local spectral subspaces of DT--operators.
We show that if $Z$ is a $\DT(\mu,c)$--operator in a W$^*$--noncommutative probability space $(\Mcal,\tau)$,
then
the spectrum of $Z$ is equal to the support of $\mu$.
Moreover,
if $B$ is any Borel set in $\Cpx$,
then the projection $p_Z(B)$ onto the closure of the spectral subspace of $B$ has trace $\tau(p_Z(B))=\mu(B)$.
Also, if $0<\mu(B)<1$, then $Zp_Z(B)$ and $(1-p_Z(B))Z$ are DT--operators
corresponding to the restrictions of
$\mu$ to $B$ and to its complement $B^c$, respectively.
From these facts, we show that the Brown measure of $Z$ is $\mu$
and $Z$ is strongly decomposable.

Since the projections $p_Z(B)$ belong to the von Neumann algebra generated by $Z$,
we have proved,
so long as the support of $\mu$ has more than one point,
that $Z$ has nontrivial invariant subspaces affiliated to the von Neumann algebra it generates.
On the other hand, let $T$ be a $\DT(\delta_0,1)$--operator.
Then $T$ is a natural candidate for an operator without any invariant subspaces relative to the von Neumann algebra
it generates.
This issue is related to questions about
the von Neumann algebra generated by the canonical copy of $T$ that we construct in $L(\Fb_2)$.
In~\S\ref{sec:vN}, we investigate this von Neumann algebra, and prove it is an irreducible subfactor of $L(\Fb_2)$.

The circular free Poisson operators were introduced in~\cite{DH}.
In that paper, we showed that circular free Poisson operators have nontrivial invariant subspaces
relative to the von Neumann algebras they generate.
As remarked above, every circular free Poisson operator is a DT--operator,
and it is clear that the invariant subspaces constructed in~\cite{DH} are the local spectral subspaces
for annular subsets of the spectrum.
The circular free Poisson operators are also examples of R--diagonal operators, in the sense
of Nica and Speicher~\cite{NS97}.
\'Sniady and Speicher~\cite{SnSp}
have extended the results of~\cite{DH} in another direction,
showing that every R--diagonal operator has an increasing, one--parameter family of
nontrivial invariant subspaces relative to its von Neumann algebra.
In~\S\ref{sec:DTRdiag}, we show that the only DT--operators that are also R-diagonal operators
are the circular free Poisson operators.
Hence, the overlap between the results of~\cite{SnSp} and of this paper lies in
their common antecedent~\cite{DH}.

In~\S\ref{sec:TstT}, we find the moments of $T^*T$:
\[
\tau((T^*T)^n)=\frac{n^n}{(n+1)!}.
\]
We show that these are the moments of a probability distribution
supported on $[0,e]$ which is Lebesgue absolutely
continuous, and we explicitly describe its density function.
As corollaries, we have $\|T\|=\sqrt e$ and $\ker T=\{0\}$.

In~\S\ref{sec:Tmom}, we prove a recursion formula for general $*$--moments
\[
\tau\big((T^*)^{k_1}T^{\ell_1}\dots(T^*)^{k_n}T^{\ell_n}\big)
\]
of $T$.
We conclude by conjecturing a pretty formula for the special $*$--moments
\[
\tau(((T^*)^kT^k)^n),
\]
which has been proved in some cases and checked by computation in others.

\section{The DT $*$--moment distributions}
\label{sec:DTdist}

In this section, we show convergence in $*$--moments of certain sequences of upper triangular
random matrices and we use them to define DT--elements.
We also prove some preliminary results about DT--elements.
The random matrix results seem to be related to those in~\cite{Sh96},
but go in a rather different direction.

As in the introduction, we will let $(\MEu_n,\tau_n)$ be the $*$--noncommutative probability space of
$n\times n$ random matrices, over a (large enough) classical probability space.
We will use the notations $\GRM$, $\UTGRM$, $\HURM$ and $\SGRM$ for special sorts
of random matrices in $\MEu_n$, as described in Notations~3.1 and~4.1 of~\cite{DH}.
In particular, 
we say that a random variable $a\in\LEu$ is {\em complex} $N(0,\sigma^2)$
if $\RealPart a$ and $\ImagPart a$ are independent real gaussian random variables each having first moment $0$ and second moment $\sigma^2/2$,
and for $T\in\MEu_n$, we have $T\in\UTGRM(n,\sigma^2)$ (the acronym is for upper triangular gaussian random matrix)
if the entries $t_{ij}$ of $T$ ($1\le i,j\le n$) satisfy that $t_{ij}=0$ whenever $i\ge j$
and that $(t_{ij})_{1\le i<j\le n}$ is an independent family
of random variables, each of which is complex $N(0,\sigma^2)$.

\begin{thm}
\label{thm:DTdef}
Let $\mu$ be a compactly supported Borel probability measure on $\Cpx$.
For $n\in\Nats$, let $D_n$ be a diagonal random matrix whose diagonal entries
are i.i.d.\ $\mu$--distributed random variables
and let $T_n\in\UTGRM(n,\frac1n)$ be such that $D_n$ and $T_n$ are independent
(as matrix--valued random variables).
Then the pair $D_n,T_n$ converges jointly in $*$--moments as $n\to\infty$.
\end{thm}

The proof will come later in this section, but we now use the result to define DT--elements.

\begin{defi}\rm
\label{def:DT}
Let $\mu$ be a compactly supported Borel probability measure on $\Cpx$ and let $c>0$.
Let $D_n$ and $T_n$ be as in Theorem~\ref{thm:DTdef} above, for this choice of $\mu$.
An element $z$ of a $*$--noncommutative probability space $(A,\phi)$
is called a $\DT(\mu,c)$--{\em element} if its $*$--moment distribution
is the limit $*$--moment distribution of $D_n+cT_n$ as $n\to\infty$,
i.e.\ if
\[
\phi(z^{\eps(1)}z^{\eps(2)}\cdots z^{\eps(k)})=
\lim_{n\to\infty}\tau_n(Z_n^{\eps(1)}Z_n^{\eps(2)}\cdots Z_n^{\eps(k)})
\]
for every $k\in\Nats$ and $\eps(1),\ldots,\eps(k)\in\{1,*\}$,
where $Z_n=D_n+cT_n$.

An element of a $*$--noncommutative probability space is called a {\em DT--element} if it is a $\DT(\mu,c)$--element
for some $\mu$ and $c>0$.
\end{defi}

We note that the letters ``DT'' signify diagonal + (gaussian upper) triangular.

\begin{remark}\rm
\label{rem:DTbdd}
We could use the Gelfand--Naimark--Segal representation to show,
at the same time we prove Theorem~\ref{thm:DTdef}, that every DT--element can be realized
as a bounded operator;
more precisely, for every $\mu$ and $c$ there is a von Neumann algebra $\Mcal$ equipped with a normal tracial state $\tau$
such that there is a $\DT(\mu,c)$--element in the $*$--noncommutative probability space $(\Mcal,\tau)$.
However, we will postpone a proof of this to~\S\ref{sec:DTinLF2}, where we will give a construction
of an arbitrary DT--element in the free group factor $(L(\Fb_2),\tau)$;
see Theorem~\ref{thm:DTinLF2} and Remark~\ref{rem:normT}.
\end{remark}

Our proof of Theorem~\ref{thm:DTdef} begins with a combinatorial analysis of the limiting $*$--moments
of $T_n$.

Recall that a \emph{pairing} of $\{1,2,\ldots,k\}$, (for $k$ even)
is
\[\sigma=\{\{i_1,j_i\},\{i_2,j_2\},\ldots,\{i_{k/2},j_{k/2}\}\}
\]
where $\{i_1,\ldots,i_{k/2},j_1,\ldots,j_{k/2}\}=\{1,\ldots,k\}$; a pairing $\sigma$ is said to be
\emph{crossing} if $i_1<i_2<j_1<j_2$ for some $\{i_1,j_1\}\in\sigma$ and $\{i_2,j_2\}\in\sigma$
and is said to be \emph{non--crossing} otherwise.

\begin{lemma}
\label{lem:utmom}
For every $n\in\Nats$ let $T_n\in\UTGRM(n,\linebreak[1]1/n)$,
Let $k\in\Nats$ and $\eps(1),\ldots,\eps(k)\in\{*,1\}$.
Then the limit
\begin{equation}
\label{eq:Tbnlim}
\lim_{n\to\infty}\tau_n(T_n^{\eps(1)}T_n^{\eps(2)}\cdots T_n^{\eps(k)})
\end{equation}
exists.
The value of the limit~\eqref{eq:Tbnlim} is zero unless there is a non--crossing pairing $\sigma$
of $\{1,2,\ldots,k\}$ such that whenever $\{i,j\}\in\sigma$ we have
$\eps(i)\ne\eps(j)$;
we will call such a non--crossing pairing $\sigma$
\emph{compatible with} $\eps(1),\ldots,\eps(k)$.
The value of the limit~\eqref{eq:Tbnlim} is
\begin{equation}
\label{eq:Tbnlimval}
\frac1{(\frac k2+1)!}\sum_{\substack{\sigma\in\NCP(k)\\ \text{compatible}}}\NTO(\sigma;\eps(1),\ldots,\eps(k)),
\end{equation}
where the sum is over all non--crossing pairings $\sigma$ of $\{1,2,\ldots,k\}$ that are compatible with
$\eps(1),\ldots,\eps(k)$, and where $\NTO(\sigma;\eps(1),\ldots,\eps(k))$ is the positive
integer obtained via the algorithm below.
\end{lemma}

\begin{algorithm}\rm
\label{alg:NTO}
Suppose $k\in\Nats$ and $\sigma=\{\{i_1,j_i\},\{i_2,j_2\},\ldots,\{i_{k/2},j_{k/2}\}\}$ is a pairing of $\{1,2,\ldots,k\}$
and let $\eps(1),\ldots,\eps(k)\in\{*,1\}$ be such that $\eps(i)\ne\eps(j)$ for every $\{i,j\}\in\sigma$.
Then the integer $\NTO(\sigma;\eps(1),\ldots,\eps(k))$ is determined by the following algorithm.

\medskip
\noindent
(A).  Let $G$ be the $k$--gon graph with consecutive vertices $v_1,\ldots,v_k$ and consecutive edges $e_1,\ldots,e_k$,
with edge $e_j$ having vertices $v_j$ and $v_{j+1}$;
(here and throughout the algorithm we take indices of vertices modulo $k$).
Orient each edge $e_j$ negatively (i.e.\ with arrow pointing from $v_{j+1}$ to $v_j$) if $\eps(j)=1$
and positively (i.e.\ with arrow pointing from $v_j$ to $v_{j+1}$) if $\eps(j)=*$.

\medskip
\noindent
(B).  Let $Q=Q(\sigma;\eps(1),\ldots,\eps(k))$ be the quotient graph of $G$ obtained by identifying edges $e_i$ and $e_j$
whenever $\{i,j\}\in\sigma$,
in such a way that the orientations given the edges in part~(A) are matched.
Thus vertex $v_i$ is identified with $v_{j+1}$ and vertex $v_{i+1}$ is identified with $v_j$.
The edges of $Q$ are given the orientations (i.e.\ arrows) inherited from the edges of $G$.
Note that taken without orientation of its edges, $Q$ is the quotient graph considered by Voiculescu in~\cite{V:RM}.

\medskip
\noindent
(C). Set $\NTO(\sigma;\eps(1),\ldots,\eps(k))=0$ if $\sigma$ is a crossing pairing.
If $\sigma$ is a non--crossing pairing then by Lemma~\ref{lem:nctree},
$Q$ is a tree and has $\frac k2+1$ vertices;
let $w_1,w_2,\ldots,w_{\frac k2+1}$ be the vertices of $Q$.
Consider the relation $\leftarrow$ on $\{w_1,w_2,\ldots,w_{\frac k2+1}\}$ defined by $w_i\leftarrow w_j$
if there is an edge in $Q$ with vertices $w_i$ and $w_j$ and whose orientation is an arrow pointing to $w_i$
from $w_j$.
Let $\NTO(\sigma;\eps(1),\ldots,\eps(k))$ be the number of different total orderings $<$ of
$\{w_1,w_2,\ldots,w_{\frac k2+1}\}$ that extend the relation $\leftarrow$.
Because $Q$ is a tree it is clear that the transitive relation generated by $\leftarrow$ is a partial ordering
and hence that $\NTO(\sigma;\eps(1),\ldots,\eps(k))\ge1$.
\end{algorithm}

\begin{example}\rm
\label{ex:NTO}
We will use Algorithm~\ref{alg:NTO} to find $\NTO(\sigma;*,1,*,1,*,1,1,*,1,*)$
when $\sigma=\{\{1,6\},\{2,3\},\{4,5\},\{7,10\},\{8,9\}\}$.
Performing parts~(A) and~(B) yields
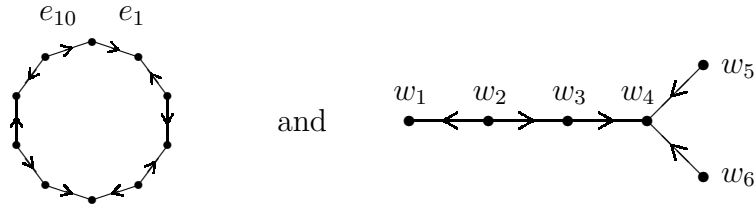
\begin{center}
\begin{picture}(270,90)(0,-10)


  \put(58.5317,39.2705){\circle*{3}} 
  \put(47.6336,54.2705){\circle*{3}}
  \put(30.,60.){\circle*{3}}
  \put(12.3664,54.2705){\circle*{3}}
  \put(1.4683,39.2705){\circle*{3}}
  \put(1.4683,20.7295){\circle*{3}}
  \put(12.3664,5.72949){\circle*{3}}
  \put(30.,0){\circle*{3}}
  \put(47.6336,5.72949){\circle*{3}}
  \put(58.5317,20.7295){\circle*{3}}

  \drawline(58.5317,39.2705)(47.6336,54.2705)(30.,60.)(12.3664,54.2705) 
           (1.4683,39.2705)(1.4683,20.7295)(12.3664,5.72949)(30.,0)
           (47.6336,5.72949)(58.5317,20.7295)(58.5317,39.2705)

  \drawline(36.0798,55.9216)(39.9923,56.7533)(37.3158,59.7259)  
  \drawline(56.0103,46.1436)(52.3561,47.7705)(52.7742,43.7924)
  \drawline(56.5317,32.228)(58.5317,28.7639)(60.5317,32.228)
  \drawline(53.3911,10.2514)(53.8092,14.2295)(50.155,12.6025)
  \drawline(40.3177,5.45536)(37.6412,2.48278)(41.5538,1.65113)
  \drawline(18.4462,1.65113)(22.3588,2.48278)(19.6823,5.45536)
  \drawline(3.98973,13.8564)(7.64392,12.2295)(7.2258,16.2076)
  \drawline(3.4683,27.772)(1.4683,31.2361)(-0.531695,27.772)
  \drawline(6.60894,49.7486)(6.19083,45.7705)(9.84501,47.3975)
  \drawline(19.6823,54.5446)(22.3588,57.5172)(18.4462,58.3489)

  \put(40,68.0423){$e_1$}  
  \put(10,68.0423){$e_{10}$}

  \put(100,26){and}


  \put(150,30){\circle*{4}} 
  \put(180,30){\circle*{4}}
  \put(210,30){\circle*{4}}
  \put(240,30){\circle*{4}}
  \put(261.213,51.213){\circle*{4}}
  \put(261.213,8.787){\circle*{4}}

  \drawline(150,30)(180,30)(210,30)(240,30)(261.213,51.213) 
  \drawline(240,30)(261.213,8.787)

   \drawline(168.196,33)(163,30)(168.196,27) 
   \drawline(191.804,33)(197,30)(191.804,27)
   \drawline(221.804,33)(227,30)(221.804,27)
   \drawline(250.745,44.988)(249.192,39.192)(254.988,40.745)
   \drawline(250.745,15.012)(249.192,20.808)(254.988,19.255)

   \put(144,38){$w_1$} 
   \put(174,38){$w_2$}
   \put(204,38){$w_3$}
   \put(230,38){$w_4$}
   \put(268,48){$w_5$}
   \put(268,7){$w_6$}

\end{picture}
\end{center}
for $G$ and $Q$, respectively.
The relation $\leftarrow$ can be drawn schematically as
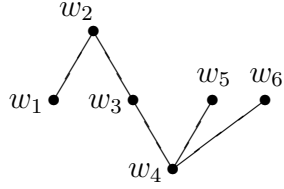
\begin{center}
\begin{picture}(130,80)(0,0)

  \put(0,30){\circle*{4}}  
  \put(15,55.981){\circle*{4}}
  \put(30,30){\circle*{4}}
  \put(45,4.019){\circle*{4}}
  \put(60,30){\circle*{4}}
  \put(80,30){\circle*{4}}

  \drawline(0,30)(15,55.981)(30,30)(45,4.019)(60,30) 
  \drawline(45,4.019)(80,30)

  \put(-17,27){$w_1$}  
  \put(2,62){$w_2$}
  \put(13,27){$w_3$}
  \put(28,1){$w_4$}
  \put(54,37){$w_5$}
  \put(74,37){$w_6$}

\end{picture}
\end{center}
and we count the number of total orderings extending $\leftarrow$ to be
\[
\NTO(\sigma;*,1,*,1,*,1,1,*,1,*)=2\cdot4\cdot5+3\cdot4=52.
\]
The term $2\cdot4\cdot5$ is for when $w_1$ is placed above $w_4$
and $3\cdot4$ is for when $w_1$ is placed below $w_4$.
\end{example}

Given a pairing $\sigma$ of $\{1,\ldots,k\}$, we say that $\sigma$ has \emph{paired neighbors}
if $\{i,i+1\}\in\sigma$ for some $i\in\{1,2,\ldots,k-1\}$.
In that case, $\sigma-\{i,i+1\}$ will denote the pairing  of $\{1,\ldots,k-2\}$ obtained from $\sigma\backslash\{\{i,i+1\}\}$
by applying the map
\[
j\mapsto\begin{cases}j&\text{if }j<i\\j-2&\text{if }j>i+1.\end{cases}
\]
We say that $\sigma-\{i,i+1\}$ is the pairing obtained by
\emph{removing paired neighbors} $\{i,i+1\}$ from $\sigma$.
For example, if $\sigma$ is as in Example~\ref{ex:NTO} then
\[
\sigma-\{4,5\}=\{\{1,4\},\{2,3\},\{5,8\},\{6,7\}\}.
\]

The following lemma is well known;
however we include a short proof for convenience.

\begin{lemma}
\label{lem:prdnbs}
A pairing $\sigma$ of $\{1,\ldots,k\}$ is non--crossing if and only if it can be reduced to the empty set
by successively removing paired neighbors.
\end{lemma}
\begin{proof}
It is easy to see by induction that every non--crossing pairing
has paired neighbors.
Indeed, if $\sigma$ is non--crossing and if $\{i,j\}\in\sigma$ with $i<j$ then the restriction of
$\sigma$ to $\{i+1,i+2,\ldots,j-1\}$ is non-crossing.
Moreover, if a pairing $\sigma$ has paired neighbors $\{i,i+1\}$
then $\sigma-\{i,i+1\}$ is crossing if and only if $\sigma$ is crossing.

Now the statement of the lemma is easily proved by induction on $k/2$.
The case $k=2$ is trivial.
Let $k\ge4$.
If $\sigma$ has no paired neighbors then it cannot be reduced at all by removing paired neighbors, and $\sigma$ is also crossing.
If $\sigma$ has paired neighbors $\{i,i+1\}$ then $\sigma$ is non--crossing if and only if $\sigma-\{i,i+1\}$ is non--crossing,
and by the induction hypothesis $\sigma-\{i,i+1\}$ can be reduced to the empty set by removing paired neighbors if and only if it is
non--crossing.
\end{proof}

\begin{remark}\rm
The above proof actually shows that every non--crossing pairing can be reduced to the empty set by removing paired
neighbors and at every stage choosing arbitrary paired neighbors to be removed.
\end{remark}

The following lemma is also well known.
Again, for convenience we include a proof.

\begin{lemma}
\label{lem:nctree}
Let $\sigma$ be a pairing of $\{1,\ldots,k\}$ and let $\eps(1),\ldots,\eps(k)\in\{*,1\}$ be such that
$\eps(i)\ne\eps(j)$ whenver $\{i,j\}\in\sigma$, and let $Q=Q(\sigma;\eps(1),\ldots,\eps(k))$ be the quotient
graph obtained in step~(B) of Algorithm~\ref{alg:NTO}.
If $\sigma$ is non--crossing then $Q$ is a tree having $\frac k2+1$ vertices, while if $\sigma$ is crossing then
$Q$ is not a tree and has $\le\frac k2$ vertices.
\end{lemma}
\begin{proof}
If $\sigma$ has no paired neighbors then every edge of the $k$--gon graph $G$ is identified with another edge not its neighbor.
Hence every vertex of $G$ gets identified with at least one other vertex and thus $Q$
has $\frac k2$ edges and $\le\frac k2$ vertices, and is therefore not a tree.
If $\sigma$ has paired neighbors $\{i,i+1\}$ then identifying edges $e_i$ and $e_{i+1}$ of $G$ we get a $(k-2)$--gon graph
with a tail consisting of one edge:
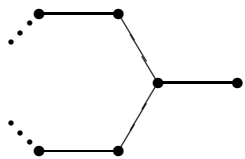
\begin{center}
\begin{picture}(100,60)


  \put(15,4.019){\circle*{4}}  
  \put(45,4.019){\circle*{4}}
  \put(60,30){\circle*{4}}
  \put(45,55.981){\circle*{4}}
  \put(15,55.981){\circle*{4}}
  \put(90,30){\circle*{4}}

  \drawline(15,4.019)(45,4.019)(60,30)(45,55.981)(15,55.981)  
  \drawline(60,30)(90,30)

  \put(11.465,52.445){\circle*{1.5}}  
  \put(7.929,48.910){\circle*{1.5}}
  \put(4.393,45.374){\circle*{1.5}}

  \put(11.465,7.555){\circle*{1.5}}  
  \put(7.929,11.090){\circle*{1.5}}
  \put(4.393,14.626){\circle*{1.5}}

\end{picture}
\end{center}
(where we have omitted to indicate the edges' orientations).
Then $Q$ is obtained from the above graph by identifying edges of the inner $(k-2)$--gon graph according to the pairing
$\sigma-\{i,i+1\}$.
Continuing in this way to identify neighboring edges as many times as possible, after $p$ steps, $1\le p\le\frac k2-2$, we
will have an inner $(k-2p)$--gon with adjoined tails, each consisting of one or more trees,
the tails together having $p$ additional vertices.
If $\sigma$ is crossing then by Lemma~\ref{lem:prdnbs} for some $p$ the inner $(k-2p)$--gon will have no neighboring edges identified,
and thus each of its $k-2p$ vertices will be paired with at least one other vertex.
The graph $Q$ will then have at most $\frac k2$ vertices, and $Q$ will not be a tree.
If $\sigma$ is non--crossing then continuing until $p=\frac k2-2$ the graph will be
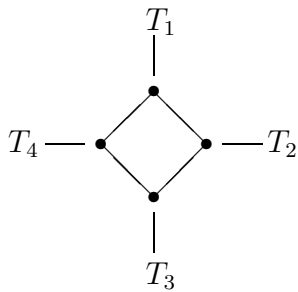
\begin{center}
\begin{picture}(100,100)(-30,-30)


  \put(0,20){\circle*{4}}  
  \put(20,0){\circle*{4}}
  \put(40,20){\circle*{4}}
  \put(20,40){\circle*{4}}

  \drawline(0,20)(20,0)(40,20)(20,40)(0,20)  

  \drawline(-6,20)(-21,20)  
  \put(-35,17){$T_4$}
  \drawline(20,-6)(20,-21)
  \put(17,-33){$T_3$}
  \drawline(46,20)(61,20)
  \put(63,17){$T_2$}
  \drawline(20,46)(20,61)
  \put(17,63){$T_1$}

\end{picture}
\end{center}
where each of $T_1,T_2,T_3$ and $T_4$ represents one or more trees adjoined for a total of $\frac k2-2$ additional vertices.
The pairing of the remaining four edges of the inner $4$--gon will be non--crossing, i.e.
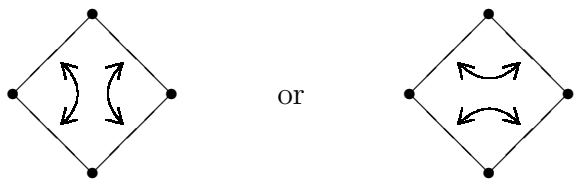
\begin{center}
\begin{picture}(220,70)(0,0)


  \put(0,30){\circle*{4}}  
  \put(30,0){\circle*{4}}
  \put(60,30){\circle*{4}}
  \put(30,60){\circle*{4}}

  \drawline(0,30)(30,0)(60,30)(30,60)(0,30)  

  \qbezier(18.536,41.465)(30,30)(18.536,18.536)          
  \drawline(20.088,35.667)(18.536,41.465)(24.331,39.912)
  \drawline(24.331,20.088)(18.536,18.536)(20.088,24.331)

  \qbezier(41.464,41.465)(30,30)(41.464,18.536)          
  \drawline(39.912,35.667)(41.454,41.465)(35.669,39.912)
  \drawline(35.669,20.088)(41.454,18.536)(39.912,24.331)

  \put(100,26){or}

  \put(150,30){\circle*{4}}  
  \put(180,0){\circle*{4}}
  \put(210,30){\circle*{4}}
  \put(180,60){\circle*{4}}

  \drawline(150,30)(180,0)(210,30)(180,60)(150,30)  

  \qbezier(168.536,41.465)(180,30)(191.464,41.465)          
  \drawline(170.088,35.667)(168.536,41.465)(174.331,39.912)
  \drawline(174.331,20.088)(168.536,18.536)(170.088,24.331)

  \qbezier(168.536,18.536)(180,30)(191.464,18.536)          
  \drawline(189.912,35.667)(191.454,41.465)(185.669,39.912)
  \drawline(185.669,20.088)(191.454,18.536)(189.912,24.331)

\end{picture}
\end{center}
giving rise to the trees
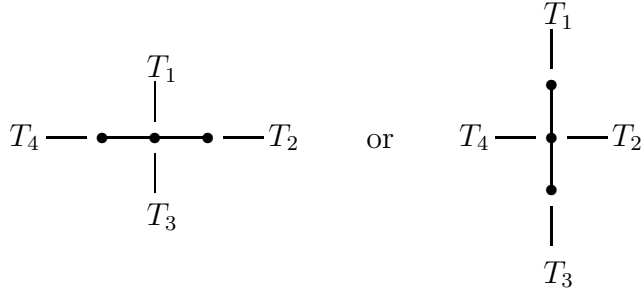
\begin{center}
\begin{picture}(250,120)(-20,-40)


  \put(0,20){\circle*{4}}  
  \put(20,20){\circle*{4}}
  \put(40,20){\circle*{4}}

  \drawline(0,20)(20,20)(40,20)  

  \drawline(-6,20)(-21,20) 
  \put(-35,17){$T_4$}
  \drawline(20,14)(20,-1)
  \put(17,-13){$T_3$}
  \drawline(46,20)(61,20)
  \put(63,17){$T_2$}
  \drawline(20,26)(20,41)
  \put(17,43){$T_1$}

  \put(100,16){or}


  \put(170,0){\circle*{4}} 
  \put(170,20){\circle*{4}}
  \put(170,40){\circle*{4}}

  \drawline(170,0)(170,20)(170,40) 

  \drawline(164,20)(149,20) 
  \put(135,17){$T_4$}
  \drawline(170,-6)(170,-21)
  \put(167,-35){$T_3$}
  \drawline(176,20)(191,20)
  \put(193,17){$T_2$}
  \drawline(170,46)(170,61)
  \put(167,63){$T_1$}

\end{picture}
\end{center}
respectively, each of which has $\frac k2+1$ vertices.
\end{proof}

\begin{proof}[Proof of Lemma~\ref{lem:utmom}]
Let $t(i,j,\eps;n)$ denote the $ij$th entry of $T_n^\eps$.
Thus
\begin{equation}
\label{eqn:tij}
\big(t(i,j,1;n)\big)_{1\le i<j\le n}
\end{equation}
is a family of i.i.d.\ complex $N(0,\frac1n)$--random variables, while $t(i,j,1;n)=0$ if $i\ge j$ and
\[
t(i,j,*\,;n)=\overline{t(j,i,1;n)}.
\]

We have
\begin{multline*}
\tau_n(T_n^{\eps(1)}T_n^{\eps(2)}\cdots T_n^{\eps(k)})= \\
=n^{-1}\sum_{i_1,\ldots,i_k\in\{1,\ldots,n\}}
\Eb\big(t(i_1,i_2,\eps(1);n)t(i_2,i_3,\eps(2);n)\cdots t(i_k,i_{k+1},\eps(k);n)\big),
\end{multline*}
where $i_{k+1}$ means $i_1$.
For a given choice of $i_1,\ldots,i_k\in\{1,\ldots,n\}$ the quantity
\[
\Eb\big(t(i_1,i_2,\eps(1);n)t(i_2,i_3,\eps(2);n)\cdots t(i_k,i_{k+1},\eps(k);n)\big)
\]
is zero unless the $k$ random variables $t(i_j,i_{j+1},\eps(j);n)$, $1\le j\le k$, can be paired off so
that
each pair consists of mutually conjugate random variables.
This happens if and only if there is a pairing $\sigma$ of $\{1,\ldots,k\}$ that is compatible with $\eps(1),\ldots,\eps(k)$
and such that
\begin{equation}
\label{eq:ipqcond}
i_p=i_{q+1}\text{ and }i_{p+1}=i_q\text{ whenever }\{p,q\}\in\sigma.
\end{equation}
In this case, let $Q=Q(\sigma;\eps(1),\ldots,\eps(k))$ be the quotient graph of the $k$--gon graph $G$
as described in Algorithm~\ref{alg:NTO}.
Label the $j$th vertex $v_j$ of $G$ with the value of $i_j$.
Then condition~\eqref{eq:ipqcond} is equivalent to the condition that vertices of $G$ that are sent by the quotient
mapping to the same vertex in $Q$ are labeled with the same value.
Thus choosing $i_1,\ldots,i_k\in\{1,\ldots,n\}$ so that~\eqref{eq:ipqcond} holds is equivalent to labeling the vertices of $Q$
with values from $\{1,\ldots,n\}$.

For a given choice of values of $i_1,\ldots,i_k$ there may be more than one pairing $\sigma$ of $\{1,\ldots,k\}$ that
is compatible with $\eps(1),\ldots,\eps(k)$ and so that~\eqref{eq:ipqcond} holds.
However, if for some {\it non--crossing} pairing $\sigma$, the values of $i_1,\ldots,i_k$ are such that
the corresponding labeling of the vertices of $Q$ is with {\it distinct} values from $\{1,\ldots,n\}$, then $\sigma$
is the unique pairing of $\{1,\ldots,k\}$ such that~\eqref{eq:ipqcond} holds.
Indeed, since $Q$, being a tree, has at most one edge connecting any two vertices, the different edges of $Q$
are distinguished by the labels of their vertices;
since $Q$ is formed by pairing the edges of $G$, in the vertex labeling $i_1,\ldots,i_k$ of $G$, each edge of $G$
must have exactly one other edge with the same set of vertex labels.
Thus any choice of $i_1,\ldots,i_n$ yielding a
distinct vertex labeling of the tree $Q$ determines a unique pairing $\sigma$ of the edges of $G$.

An upper bound for
\begin{equation}
\label{eq:Tmom}
\tau_n(T_n^{\eps(1)}T_n^{\eps(2)}\cdots T_n^{\eps(k)})
\end{equation}
is
\begin{equation}
\label{eq:Tmomub}
n^{-1}\sum_{\substack{\sigma\in\operatorname{P}(k)\\ \text{compatible}}}
\sum_{\substack{\text{vertex labelings of}\\ Q(\sigma;\eps(1),\ldots,\eps(k))\\ \text{from }\{1,\ldots,n\}}}
\Eb\big(t(i_1,i_2,\eps(1);n)\cdots t(i_k,i_{k+1},\eps(k);n)\big),
\end{equation}
where the first sum is over all pairings $\sigma$ of $\{1,\ldots,k\}$ that are compatible with
$\eps(1),\ldots,\eps(k)$ and the second sum is over all vertex labelings of the quotient
graph $Q=Q(\sigma;\eps(1),\ldots,\eps(k))$ with elements of $\{1,\ldots,n\}$;
then $i_1,\ldots,i_k$ appearing in~\eqref{eq:Tmomub} are assigned values based on the vertex labeling of $Q$
by giving $i_j$ the value of the label of the vertex of $Q$ to which the vertex $v_j$ in $G$ is mapped
under the quotient mapping.
A lower bound for~\eqref{eq:Tmom} is
\begin{equation}
\label{eq:Tmomlb}
n^{-1}\sum_{\substack{\sigma\in\NCP(k)\\ \text{compatible}}}
\sum_{\substack{\text{distinct vertex}\\ \text{labelings of }\\ Q(\sigma;\eps(1),\ldots,\eps(k))\\ \text{from }\{1,\ldots,n\}}}
\Eb\big(t(i_1,i_2,\eps(1);n)\cdots t(i_k,i_{k+1},\eps(k);n)\big),
\end{equation}
where the first sum is over all non--crossing pairings $\sigma$ of $\{1,\ldots,k\}$ that are compatible with
$\eps(1),\ldots,\eps(k)$ and the second sum is over all labelings of the $\frac k2+1$ vertices
of $Q(\sigma;\eps(1),\ldots,\eps(k))$ with distinct elements from $\{1,\ldots,n\}$, the corresponding values of $i_1,\ldots,i_k$
in~\eqref{eq:Tmomlb} being determined as described above.
The difference between the upper bound~\eqref{eq:Tmomub} and the lower bound~\eqref{eq:Tmomlb} is
\begin{equation}
\label{eq:Tmomdif}
\begin{aligned}
n^{-1}&\sum_{\substack{\sigma\in\NCP(k)\\ \text{compatible}}}
\;\sum_{\substack{\text{non--distinct}\\ \text{vertex labelings of }\\ Q(\sigma;\eps(1),\ldots,\eps(k))\\ \text{from }\{1,\ldots,n\}}}
\Eb\big(t(i_1,i_2,\eps(1);n)\cdots t(i_k,i_{k+1},\eps(k);n)\big) \\
+n^{-1}&\sum_{\substack{\sigma\in\operatorname{P}(k)\backslash\NCP(k)\\ \text{compatible}}}
\;\sum_{\substack{\text{vertex labelings of}\\ Q(\sigma;\eps(1),\ldots,\eps(k))\\ \text{from }\{1,\ldots,n\}}}
\Eb\big(t(i_1,i_2,\eps(1);n)\cdots t(i_k,i_{k+1},\eps(k);n)\big),
\end{aligned}
\end{equation}
where $\sigma\in\operatorname{P}(k)\backslash\NCP(k)$ means that $\sigma$ is a crossing pairing of $\{1,\ldots,k\}$.
There is a constant $C_k$ depending only on $k$ such that for all $n\in\Nats$ and all $i_1,\ldots,i_k\in\{1,\ldots,n\}$,
\[
0\le\Eb\big(t(i_1,i_2,\eps(1);n)\cdots t(i_k,i_{k+1},\eps(k);n)\big)\le C_kn^{-k/2}.
\]
For a given $\sigma\in\NCP(k)$, the number of different non--distinct vertex labelings of $Q(\sigma;\eps(1),\ldots,\eps(k)$
is bounded above by $\binom{\frac k2+1}2n^{k/2}$.
Hence
\[
0\le\sum_{\substack{\text{non--distinct}\\ \text{vertex labelings of }\\ Q(\sigma;\eps(1),\ldots,\eps(k))\\ \text{from }\{1,\ldots,n\}}}
\Eb\big(t(i_1,i_2,\eps(1);n)\cdots t(i_k,i_{k+1},\eps(k);n)\big)\le C_k\binom{\frac k2+1}2.
\]
For a given $\sigma\in\operatorname{P}(k)\backslash\NCP(k)$, since by Lemma~\ref{lem:nctree}
$Q=Q(\sigma;\eps(1),\ldots,\eps(k))$ has $\le\frac k2$ vertices, the number of different vertex labelings of $Q$
is bounded above by $n^{k/2}$ and therefore
\[
0\le\sum_{\substack{\text{vertex labelings of}\\ Q(\sigma;\eps(1),\ldots,\eps(k))\\ \text{from }\{1,\ldots,n\}}}
\Eb\big(t(i_1,i_2,\eps(1);n)\cdots t(i_k,i_{k+1},\eps(k);n)\big)\le C_k\;.
\]
Thus we see that the difference~\eqref{eq:Tmomdif} tends to zero as $n\to\infty$, and the limit
as $n\to\infty$ of~\eqref{eq:Tmom}
will equal the limit as $n\to\infty$ of the lower bound~\eqref{eq:Tmomlb}, if the latter exists.

We will now show that the lower bound~\eqref{eq:Tmomlb} converges as $n\to\infty$ to the desired value~\eqref{eq:Tbnlimval}.
Fix $\sigma\in\NCP(k)$ compatible with $\eps(1),\ldots,\eps(k)$, fix a distinct vertex labeling
of $Q=Q(\sigma;\eps(1),\ldots,\eps(k))$ from $\{1,\ldots,n\}$ and let $i_1,\ldots,i_k\in\{1,\ldots,n\}$ be the
values determined by this vertex labeling.
Then
\begin{equation}
\label{eq:Ebt}
\Eb\big(t(i_1,i_2,\eps(1);n)\cdots t(i_k,i_{k+1},\eps(k);n)\big)
\end{equation}
is equal to $n^{-k/2}$ if $i_j<i_{j+1}$ whenever $\eps(j)=1$ and $i_j>i_{j+1}$ whenever $\eps(j)=*$,
and is equal to zero otherwise.
In terms of the orientations of the edges of the graph $G$ as determined by $\eps(1),\ldots,\eps(k)$ in part~(A)
of Algorithm~\ref{alg:NTO}, the expectation~\eqref{eq:Ebt} is $n^{-k/2}$ if $i_j<i_{j+1}$ whenever
the arrow points to vertex $v_j$ from vertex $v_{j+1}$ and $i_j>i_{j+1}$ whenever the arrow points from
$v_j$ to $v_{j+1}$, and the quantity~\eqref{eq:Ebt} is zero otherwise.
In terms of the orientations of the edges of the graph $Q$, the quantity~\eqref{eq:Ebt} is $n^{-k/2}$ if the arrows
always point from the vertex labeled with the larger value to the vertex labeled with the smaller value.
Thus for a fixed $\sigma$, the quantity
\begin{equation}
\label{eq:Ebd}
n^{-1}\sum_{\substack{\text{distinct vertex}\\ \text{labelings of}\\ Q(\sigma;\eps(1),\ldots,\eps(k))\\ \text{from }\{1,\ldots,n\}}}
\Eb\big(t(i_1,i_2,\eps(1);n)\cdots t(i_k,i_{k+1},\eps(k);n)\big)
\end{equation}
is equal to $n^{-(\frac k2+1)}$ times the number of distinct vertex labelings of $Q$
such that the label of vertex $w_i$ of $Q$ is less than the label of vertex $w_j$ of $Q$
whenever $w_i\leftarrow w_j$,
where $\leftarrow$ is the relation defined in part~(C) of Algorithm~\ref{alg:NTO}.
Because $Q$ has $\frac k2+1$ vertices, as $n\to\infty$ the quantity~\eqref{eq:Ebd} converges to the volume with respect to
Lebesgue measure of the set $V$ of all $(\frac k2+1)$--tuples $(t_1,\ldots,t_{\frac k2+1})$ in the cube $[0,1]^{\frac k2+1}$
such that $t_i<t_j$ whenever $w_i\leftarrow w_j$.
Partitioning this cube into $(\frac k2+1)!$ wedges of equal measure corresponding to the different total orderings
of the $\frac k2+1$ co--ordinates, we see that $V$ is the union of $\NTO(\sigma;\eps(1),\ldots,\eps(k))$ different wedges;
hence~\eqref{eq:Ebd} converges as $n\to\infty$ to
\[
\frac1{(\frac k2+1)!}\NTO(\sigma;\eps(1),\ldots,\eps(k)).
\]
Summing over all non--crossing pairings yields the expression~\eqref{eq:Tbnlimval}
for the limit $*$--moment~\eqref{eq:Tbnlim}.
\end{proof}

\begin{proof}[Proof of Theorem~\ref{thm:DTdef}]
For any $a,b\in\Nats\cup\{0\}$, let
\[
D_n(a,b)=D_n^a(D_n^*)^b.
\]
Then we must show that the limit
\begin{equation}
\label{eq:DTlim}
\lim_{n\to\infty}\tau_n(D_n(a_1,b_1)T_n^{\eps(1)}\cdots D_n(a_k,b_k)T_n^{\eps(k)})
\end{equation}
exists, for any $k\in\Nats$, $a_1,\ldots,a_k,b_1,\ldots,b_k\in\Nats\cup\{0\}$
and $\eps(1),\ldots,\eps(k)\in\{*,1\}$.
Writing $t(i,j,\eps;n)$ for the $(i,j)$th element of $T_n^\eps$
as in the proof of Lemma~\ref{lem:utmom}
and letting $d(i;a,b,n)$ denote the $i$th diagonal entry of $D_n(a,b)$,
using the independence of $D_n$ and $T_n$ we have
\begin{align*}
\tau_n(D_n(a_1,b_1)&T_n^{\eps(1)}\cdots D_n(a_k,b_k)T_n^{\eps(k)})= \\[1ex]
=n^{-1}&\sum_{i_1,\ldots,i_k\in\{1,\ldots,n\}}\Eb\big(d(i_1;a_1,b_1,n)\cdots d(i_k;a_k,b_k,n)\big)\cdot \\[1ex]
&\qquad\cdot\Eb\big(t(i_1,i_2,\eps(1);n)t(i_2,i_3,\eps(2);n)\cdots t(i_k,i_{k+1},\eps(k);n)\big)\;.
\end{align*}
Since the hypotheses imply that
\[
\sup_{n\in\Nats}\quad\sup_{i_1,\ldots,i_k\in\{1,\ldots,n\}}|\Eb(d(i_1;a_1,b_1,n)\cdots d(i_k;a_k,b_k,n))|<\infty\;,
\]
the same arguments as in the proof of Lemma~\ref{lem:utmom} imply that the limit~\eqref{eq:DTlim}
is equal to the limit as $n\to\infty$ of
\begin{equation}
\label{eq:limTD2}
n^{-1}\sum_{\substack{\sigma\in\NCP(k)\\ \text{compatible}}}
\sum_{
\substack{\text{distinct vertex}\\
 \text{labelings of }\\
 Q(\sigma;\eps(1),\ldots,\eps(k))\\
 \text{from }\{1,\ldots,n\}}}
\begin{aligned}[t]
\bigg(&\Eb\big(d(i_1;a_1,b_1,n)\cdots d(i_k;a_k,b_k,n)\big)\cdot \\
&\cdot\Eb\big(t(i_1,i_2,\eps(1);n)\cdots t(i_k,i_{k+1},\eps(k);n)\big)\bigg)\;,\end{aligned}
\end{equation}
if the latter limit exists,
where the first sum in~\eqref{eq:limTD2} is over all non--crossing pairings $\sigma$ of $\{1,\ldots,k\}$
that are compatible with $\eps(1),\ldots,\eps(k)$ and the second sum is over all labelings
of the $\frac k2+1$ vertices of $Q=Q(\sigma;\eps(1),\ldots,\eps(k))$ with distinct elements of $\{1,\ldots,n\}$,
with each $i_j$ being assigned the value of the label of the vertex of $Q$ to which the vertex $v_j$ of the $k$--gon graph $G$ is mapped
via the quotient mapping $G\to Q$.

For fixed $\sigma\in\NCP(k)$ compatible with $\eps(1),\ldots,\eps(k)$, let
$w_1,\linebreak[2]\ldots,\linebreak[1]w_{\frac k2+1}$
be the vertices of $Q=Q(\sigma;\eps(1),\ldots,\eps(k))$
and for $1\le p\le\frac k2+1$ let $E(p)$ be the set of all $j\in\{1,\ldots,k\}$ such that the vertex $v_j$
of $G$ is mapped to $w_p$ by the quotient map.
Fix an arbitrary distinct vertex labeling of $Q$ from $\{1,\ldots,n\}$
and let $\iotah_p$ be the label of the vertex $w_p$.
Because the entries of $D_n$ are mutually independent, we have
\begin{equation}
\label{eq:Edlim}
\begin{aligned}
\Eb(d(i_1;a_1,b_1,n)\cdots d(i_k;a_k,b_k,n))
&=\prod_{p=1}^{\frac k2+1}\Eb({\textstyle \prod_{j\in E(p)}d(\iotah_p;a_j,b_j},n)) \\
&=\prod_{p=1}^{\frac k2+1}\int\lambda^{r(p)}\lambdab^{\,s(p)}\dif\mu(\lambda)\;,
\end{aligned}
\end{equation}
where $r(p)=\sum_{j\in E(p)}a_j$ and $s(p)=\sum_{j\in E(p)}b_j$.
The quantity~\eqref{eq:Edlim} is thus independent of the particular distinct vertex labeling,
and we will denote it
\[
\Ec(\sigma;a_1,b_1,\ldots,a_k,b_k;\eps(i),\ldots\eps(k)).
\]
Using the same analysis as in the proof of Lemma~\ref{lem:utmom}, we see that
the quantity~\eqref{eq:limTD2} converges as $n\to\infty$ to the complex number
\begin{equation}
\label{eq:ENTOsum}
\frac1{(\frac k2+1)!}\sum_{\substack{\sigma\in\NCP(k)\\ \text{compatible}}}
\Ec(\sigma;a_1,b_1,\ldots,a_k,b_k;\eps(i),\ldots\eps(k))
\NTO(\sigma;\eps(1),\ldots,\eps(k))\;.
\end{equation}
In particular, the limit~\eqref{eq:DTlim} exists.
\end{proof}

In the case of $*$--moments
of a DT--element, we easily obtain the following corollary.

\begin{cor}
\label{cor:poly}
For every $\ell\in\Nats$ and $\eps(1),\ldots,\eps(\ell)\in\{*,1\}$ there is a polynomial
$P=P_{\mu;\eps(1),\ldots,\eps(\ell)}$ in $1+(\ell+1)^2$
variables and having nonnegative real coefficients such that if $z$ is a $\DT(\mu,c)$ element
in a $*$--noncommutative probability space $(\Afr,\tau)$ then
\[
\tau(z^{\eps(1)}z^{\eps(2)}\cdots z^{\eps(\ell)})=P\big(c^2,(M_\mu(r,s))_{r,s=0}^\ell\big),
\]
where $M_\mu(r,s)=\int_\Cpx\lambda^r\lambdab^{\,s}\dif\mu(\lambda)$.
Furthermore, if we assign $\deg(M_\mu(r,s))=r+s$ and $\deg(c)=1$, then $P$ is homogeneous of degree $\ell$.
\end{cor}

\begin{cor}
\label{cor:DTconv}
Let $\mu$ and $\mu_n$ ($n\in\Nats$) be compactly supported Borel probability measures on $\Cpx$
such that $\mu_n$ converges in $*$--moments as $n\to\infty$ to $\mu$.
Let $c_n>0$ ($n\in\Nats$) be such that $\lim_{n\to\infty}c_n=c\in(0,\infty)$.
Let $z_n$ be a $\DT(\mu_n,c_n)$ element.
Then $z_n$ converges in $*$--moments as $n\to\infty$ to a $\DT(\mu,c)$ element.
\end{cor}

\begin{prop}
\label{prop:constadj}
Let $z$ be a $\DT(\mu,c)$ element.
\renewcommand{\labelenumi}{(\roman{enumi})}
\begin{enumerate}
\item If $\lambda\in\Cpx\backslash\{0\}$ then $\lambda z$ is $\DT(M_\lambda\mu,|\lambda|c)$,
where $M_\lambda\mu$ is the measure $M_\lambda\mu(B)=\mu(\lambda^{-1}B)$.
\item $z^*$ is $\DT(\overline\mu,c)$ where $\overline\mu$ is the measure $\overline\mu(B)=\mu(\overline B)$.
\end{enumerate}
\end{prop}
\begin{proof}
Let $Z_n=D_n+cT_n$ where $D_n\in\MEu_n$ is a diagonal random matrix whose $n$ diagonal entries
are independent, each having distribution $\mu$, where $T_n\in\UTGRM(n,\frac1n)$ and where $D_n$ and $T_n$
are independent.
Then $Z_n$ converges in $*$--moments to $z$.

For~(i), $\lambda Z_n$ converges to $\lambda z$ in $*$--moments, as $n\to\infty$.
But each diagonal entry of $\lambda D_n$ has distribution $M_\lambda\mu$ and $\frac\lambda{|\lambda|}T_n\in\UTGRM(n,\frac1n)$,
so $\lambda Z_n$ converges in $*$--moments to a $\DT(M_\lambda\mu,|\lambda|c)$ element.

For~(ii), $Z_n^*$ converges to $z^*$ in $*$--moments.
Note that if $U_n$ is the $n\times n$ permutation matrix effecting the permutation
\[
1\mapsto n,\quad2\mapsto n-1,\quad\ldots,\quad n\mapsto1
\]
then $U_nZ_n^*U_n^*$ has the same distribution as $\overline{D_n}+cT_n$,
hence $Z_n^*$ converges to a $\DT(\mubar,c)$ element as $n\to\infty$.
\end{proof}

Applying Theorem~3.6 of~\cite{DH},
we gain some important flexibility in choosing the diagonal
part in random matrix approximations of DT--elements.

\begin{thm}
\label{thm:DTflexdiag}
Let $\mu$ be a compactly supported Borel probability measure on $\Cpx$.
For every $n\in\Nats$ let $D_n\in\MEu_n$ be a diagonal random matrix and denote by
$\nu_n$ the joint distribution of the $n$ diagonal entries of $D_n$.
Let $\nut_n$ be the symmetrixation of $\nu_n$ and for $p\in\{1,\ldots,n\}$
let $\nut_n^{(p)}$ be the $p$th marginal distribution of $\nut_n$.
Suppose for every $p\in\Nats$ the measure $\nut_n^{(p)}$ converges in $*$--moments to $\btimes_1^p\mu$.

For every $n\in\Nats$ let $T_n\in\UTGRM(n,\frac1n)$ be such that $T_n$ and $D_n$ are independent matrix--valued random variables.
Let $c>0$.
Then $D_n+cT_n$ converges in $*$--moments as $n\to\infty$ to a $\DT(\mu,c)$--element.
\end{thm}

\begin{cor}
\label{cor:nonrandD}
Let $\mu$ be a compactly supported Borel measure on $\Cpx$ and let $D_n\in M_n(\Cpx)$ be a sequence
of non--random diagonal matrices with $\sup_n\|D_n\|<\infty$.
Assume that the $*$--moments of $D_n$ converge to the the $*$--moments of $\mu$ as $n\to\infty$.
Let $T_n\in\UTGRM(n,\frac1n)$ and let $c>0$.
Then $D_n+cT_n$ converges in $*$--moments to a $\DT(\mu,c)$--element.
\end{cor}
\begin{proof}
We have $D_n=\diag(d_1^{(n)},\ldots,d_n^{(n)})$ where $d_j^{(n)}\in\Cpx$.
Let $\mu_j^{(n)}$ denote the dirac measure at $d_j^{(n)}$.
Then with the notation of Theorem~\ref{thm:DTflexdiag},
\[
\nu_n=\mu_1^{(n)}\times\cdots\times\mu_n^{(n)}\qquad
\nut_n=\frac1{n!}\sum_{\pi\in S_n}\mu_{\pi(1)}^{(n)}\times\cdots\times\mu_{\pi(n)}^{(n)}\;.
\]
Hence the first marginal distribution is
\[
\nut_n^{(1)}=\frac1n(\mu_1^{(n)}+\cdots+\mu_n^{(n)})\;.
\]
By the assumptions, $\nut_n^{(1)}$ converges in $*$--moments to $\mu$ as $n\to\infty$.
Since $\supp(\nut_n^{(1)})\subseteq K$ and $\supp(\mu)\subseteq K$, where $K$ is the closed disk of radius
$R=\sup_n\|D_n\|$ centered at the origin, by the Stone--Weierstrass theorem, we also have $\nut_n^{(1)}\to\mu$
in w$^*$--topology as $n\to\infty$.

We next consider the $p$th marginal distribution $\nu_n^{(p)}$ for $p\ge2$.
If $p=2$ we get
\[
\nut_n^{(2)}=\frac1{n(n-1)}\sum_{i\ne j}\mu_i^{(n)}\times\mu_j^{(n)}
\le\frac1{n(n-1)}\sum_{i,j=1}^n\mu_i^{(n)}\times\mu_j^{(n)}=\frac n{n-1}\,\nut_n^{(1)}\times\nut_n^{(1)}\;.
\]
Hence $\nut_n^{(2)}-\nut_n^{(1)}\times\nut_n^{(1)}=\lambda-\rho$, where $\lambda$ and $\rho$ are the positive measures
\[
\lambda=\frac1{n-1}\nut_n^{(1)}\times\nut_n^{(1)}\qquad
\rho=\frac n{n-1}\nut_n^{(1)}\times\nut_n^{(1)}-\nut_n^{(2)}\;,
\]
both of which have total mass $1/(n-1)$.
Thus the total variation $\|\nut_n^{(2)}-\nut_n^{(1)}\times\nut_n^{(1)}\|$
of $\nut_n^{(2)}-\nut_n^{(1)}\times\nut_n^{(1)}$ is at most $2/(n-1)$ and
\[
\text{w$^*$-}\lim_{n\to\infty}\nut_n^{(2)}=\text{w$^*$-}\lim_{n\to\infty}(\nut_n^{(1)}\times\nut_n^{(1)})=\mu\times\mu\;.
\]
Since all the measures involved have support in $K\times K$, convergence in $*$--moments of $\nut_n^{(2)}$ to $\mu\times\mu$
follows immediately.

The above argument can easily be generalized to $p>2$.
For $n\ge p$ we have
\[
\nut_n^{(p)}=\frac{(n-p)!}{n!}\sum\mu_{i_1}^{(n)}\times\cdots\times\mu_{i_p}^{(n)}\;,
\]
where the sum is over all $p$--tuples $(i_1,\ldots,i_p)$ of distinct integers in $\{1,\ldots,n\}$.
Hence
\[
\nut_n^{(p)}\le\frac{(n-p)!\,n^p}{n!}\btimes_1^p\nut_n^{(1)}\;.
\]
Arguing as above, we get
\[
\|\nut_n^{(p)}-\btimes_1^p\nut_n^{(1)}\|\le2\big(\frac{(n-p)!\,n^p}{n!}-1\big)\to0\text{ as }n\to\infty,
\]
from which we obtain convergence of $\nut_n^{(p)}$ to $\btimes_1^p\mu$ both in w$^*$--topology and in $*$--moments.
\end{proof}

\section{Examples of DT--elements}
\label{sec:DTex}

Recall that a circular element, as defined by Voiculescu (cf. \cite{V90}),
is an element of a $*$--noncommutative probability space
whose $*$--moment distribution is the same as $x_1+ix_2$,
where $x_1$ and $x_2$ are free, centered semicircular elements
having the same second moment.
The \emph{circular free Poisson} elements, one of which is a circular element,
were defined in~\cite{DH}.
Specifically, the circular free Poisson element of parameter $c\ge1$ is an element
of a $*$--noncommutative probability space having the same $*$--moment distribution
as $uh_c$ where $u$ is a Haar unitary, $h_c\ge0$, $h_c^2$ has the free Poisson
distribution of parameter $c$ and the pair $u,h_c$ is $*$--free.

There it was shown, using a result of F.\ Dyson for the uppertriangularization
of a nonsymmetric gaussian random matrix,
that circular free Poisson elements are
limits in $*$--moments of certain upper triangular random matrices,
and that these random matrices satisfy properties
implying the following theorem.

\begin{thm}
\label{thm:cfPDT}
Circular elements and circular free Poisson elements are DT--ele\-ments.
More specifically,
letting $\tau$ be the functional in the $*$--noncommutative probability space,
\renewcommand{\labelenumi}{(\roman{enumi})}
\begin{enumerate}

\item
a circular element $z$ satisfying $\tau(z^*z)=r^2$, 
is $\DT(\rho_r,r)$ where $\rho_r$ is uniform measure supported
on the disk of radius $r$ centered at $0$;

\item
a circular free Poisson element of parameter $c$ is $\DT(\rho_{\sqrt{c-1},\sqrt c},1)$,
where $\rho_{\sqrt{c-1},\sqrt c}$ is uniform measure supported
on the annulus having radii $\sqrt{c-1}$ and $\sqrt c$.
\end{enumerate}
\end{thm}

The elliptic elements form another class
which includes the circular element.
Their Brown measures were computed by F.~Larsen in~\cite{L},
(see also~\cite{BL} and~\cite{HL})
and their matrix models have been considered in~\cite{HP}.
We will say that an \emph{elliptic element} is an element of
a $*$--noncommutative probability space whose $*$--moment distribution
is the same as $ax_1+ibx_2$ for some $a,b>0$, where $x_1$ and $x_2$ are
free centered semicircular elements having second moment equal to $1$.
The remainder of this section is dedicated to showing that
elliptic elements are DT--elements.

\begin{thm}
\label{thm:elliptic}
An elliptic element $ax_1+bx_2$ as above is a $\DT(\nu_{a,b},c_{a,b})$--element,
where $\nu_{a,b}$ is uniform measure supported on the solid ellipse
\[
\bigg\{z\in\Cpx\,\bigg|\,\frac{(\RealPart z)^2}{4a^4}+\frac{(\ImagPart z)^2}{4b^4}
\le\frac1{a^2+b^2}\bigg\}
\]
and where
\[
c_{a,b}=\frac{2ab}{\sqrt{a^2+b^2}}\;.
\]
\end{thm}

The proof of this theorem will use the following result concerning
random matrices, which is derived from Dyson's result mentioned above.

Before stating this result, let us introduce some notation
for subsets of $M_n(\Cpx)$:
\renewcommand{\labelenumi}{$\bullet$}
\begin{enumerate}

\item
$M_n(\Cpx)_{s.a.}$, the self--adjoint matrices;

\item
$\Uc_n$, the unitary matrices;

\item
$\Dc_n$, the diagonal matrices;

\item
$\Tc_n$, the strictly upper triangular matrices,
$\{(t_{ij})_{1\le i,j\le n}\in M_n(\Cpx)\mid t_{ij}=0\text{ if }i\ge j\}$.

\end{enumerate}
By Lebesgue measure on $M_n(\Cpx)$, $\Dc_n$ and $\Tc_n$, we will mean the product
of Lebesgue measure on the real and imaginary parts of all their elements
that are not {\em a priori} zero.
By Lebesgue measure on $M_n(\Cpx)_{s.a.}$, we will mean of course
the product of Lebesgue measure on the real diagonal entries and
Lebesgue measure of the real and imaginary parts of all entries strictly
above the diagonal.

An $n\times n$ random matrix $H$ is said to belong to the class $\SGRM(n,\sigma^2)$
if it is self--adjoint, if its $\frac{n(n-1)}2$ entries strictly above the diagonal
are complex $N(0,\sigma^2)$ random variables, if its $n$ diagonal entries are real $N(0,\sigma^2)$
random variables and if together they form a collection of $\frac{n(n+1)}2$ of mutually independent random variables.
On the other hand, an element of $\HURM(n)$ is an $n\times n$ random unitary matrix that is distributed according to Haar measure
on the group of $n\times n$ unitaries.

\begin{thm}
\label{thm:rmelliptic}
Let $0<\theta<\frac\pi2$, let $n\in\Nats$ and let
\[
Y_\theta(n)=(\cos\theta)H_1(n)+i(\sin\theta)H_2(n),
\]
where $H_1(n)$ and $H_2(n)$ are random matrices in the class $\SGRM(n,\frac1n)$
which are independent as matrix--valued random variables.
Let
\[
Z_\theta(n)=U^*(n)(D_\theta(n)+\sin(2\theta)T(n))U(n)\;,
\]
where
$U(n)\in\HURM(n)$, $T(n)\in\UTGRM(n,\frac1n)$ and $D_\theta(n)\in\MEu_n$
is a diagonal random matrix whose $n$ diagonal entries have the distribution whose
density with respect to Lebesgue measure on $\Dc_n$ is
\begin{equation}
\label{eq:Dthetadist}
\begin{aligned}
\rho_\theta(D)=K_\theta\bigg(\prod_{1\le i<j\le n}|d_i-d_j|^2\bigg)
\exp\bigg(-\frac n2\sum_{1\le i\le n}
\Big(\frac{(\RealPart\,d_i)^2}{\cos^2\theta}
+\frac{(\ImagPart\,d_i)^2}{\sin^2\theta}\Big)\bigg),& \\
(D=\diag(d_1,\ldots,d_n)\in\Dc_n&)
\end{aligned}
\end{equation}
for the appropriate constant $K_\theta$.
Then $Y_\theta$ and $Z_\theta$ have the same $M_n(\Cpx)$--valued distribution.
\end{thm}

We would like to point out that the
eigenvalue distribution~\eqref{eq:Dthetadist} of $Y_\theta(n)$
was previously
found by Hiai and Petz~\cite[Lemma 4.1.10]{HP} using different techniques.

Consider the action of $\Uc_n$ on $M_n(\Cpx)$ by conjugation and let
$M_n(\Cpx)/\Uc_n$ denote the measure space of equivalence classes.
Every element of $M_n(\Cpx)/\Uc_n$ has a representative belonging to $\Dc_n+\Tc_n$,
(and usually several of them).
With respect to the quotient maps
\[
\xymatrix{
M_n(\Cpx) \ar[dr] & & \Dc_n+\Tc_n  \ar[dl] \\
& M_n(\Cpx)/\Uc_n
}
\]
every measure $\mu$ on $M_n(\Cpx)$ induces a measure on $\Dc_n+\Tc_n$ having
the same push--forward to $M_n(\Cpx)/\Uc_n$ as does $\mu$.

\begin{proof}[Proof of Theorem~\ref{thm:rmelliptic}]
F.\ Dyson proved this (see~\cite[A.35]{Mehta}) when $\theta=\pi/4$.
For general $\theta$, the distribution of $Y_\theta(n)$ has density with respect to Lebesgue measure
on $M_n(\Cpx)$
\[
\psi_\theta(A+iB)=K_\theta'
\exp\bigg(-\frac n2\Tr\Big(\frac{A^2}{\cos^2\theta}+\frac{B^2}{\sin^2\theta}\Big)\bigg),
\qquad (A,B\in M_n(\Cpx)_{s.a.}),
\]
while the distribution of $D_\theta(n)+\sin(2\theta)T(n)$ has density with respect to Lebesgue
measure on $\Dc_n+\Tc_n$
\[
\phi_\theta(D+S)=K_\theta''\rho_\theta(D)\exp\Big(-\frac n{\sin^2(2\theta)}\Tr(S^*S)\Big)
\qquad(D\in\Dc_n,\,S\in\Tc_n)
\]
for appropriate constants $K_\theta'$ and $K_\theta''$.
Since $\psi_\theta$ is invariant under unitary conjugation, as is the distribution of $Z_\theta(n)$,
and since $\phi_\theta$ agrees on elements of $\Dc_n+\Tc_n$ that have the same image in
$M_n(\Cpx)/\Uc_n$, in order to prove the theorem it will suffice to show
\begin{equation*}
\forall X\in\Dc_n+\Tc_n\qquad
\frac{\psi_\theta(X)}{\psi_{\frac\pi4}(X)}=c_\theta\frac{\phi_\theta(X)}{\phi_{\frac\pi4}(X)}\;,
\end{equation*}
where $c_\theta$ is a constant depending only on $\theta$.
Let $X=D+S$ for $D\in\Dc_n$ and $S\in\Tc_n$.
Then $X=A+iB$ for $A,B\in M_n(\Cpx)_{s.a.}$ given by
\begin{align*}
A&=\RealPart(D)+\tfrac12(S+S^*) \\
B&=\ImagPart(D)+\tfrac1{2i}(S-S^*)\;.
\end{align*}
Since $\RealPart(D)$, $S$ and $S^*$ are orthogonal in $M_n(\Cpx)$ with respect to the inner product
determined by the trace, we have
\[
\Tr(A^2)=\Tr(\RealPart(D)^2)+\tfrac12\Tr(S^*S)
\]
and similarly
\[
\Tr(B^2)=\Tr(\ImagPart(D)^2)+\tfrac12\Tr(S^*S)\;.
\]
Hence
\[
\psi_\theta(D+S)=K_\theta'\exp\Big(-\frac n2\Tr
\Big(\frac{\RealPart(D)^2}{\cos^2\theta}+\frac{\ImagPart(D)^2}{\sin^2\theta}\Big)\Big)
\exp\Big(-\frac n{\sin^2(2\theta)}\Tr(S^*S)\Big)\;.
\]
From this one obtains that $\frac{\psi_\theta(D+S)}{\psi_{\frac{\pi}4}(D+S)}$ and
$\frac{\phi_\theta(D+S)}{\phi_{\frac{\pi}4}(D+S)}$ coincide up to multiplication by a constant
depending only on $\theta$.
\end{proof}

We will need Hadamard's Determinant Theorem:
\begin{thm}[\cite{Had}]
\label{thm:Hadamard}
Let $A=(a_{ij})_{i,j=1}^n$ be a positive, semidefinite $n\times n$ matrix, i.e.\
$A\in M_n(\Cpx)_+\,$.
Then
\[
\det(A)\le\prod_{i=1}^n a_{ii}\;.
\]
\end{thm}

\begin{lemma}
\label{lem:margdistinequal}
Let $n\in\Nats$, $n\ge2$.
Let $\dif m(z)=\dif(\RealPart z)\dif(\ImagPart z)$
denote the Lebesgue measure on $\Cpx$ and let $g:\Cpx\to[0,\infty)$ be a Borel function
such that
\[
0<\int_\Cpx g(z)|z|^k\dif m(z)<\infty
\]
for $k=0,1,2,\ldots$.
Put
\[
\sigma(z_1,\ldots,z_n)=c\prod_{1\le i<j\le n}|z_i-z_j|^2\prod_{j=1}^ng(z_j)\;,
\]
where the normalization constant $c$ is chosen so that $\sigma$ is a probability density on $\Cpx^n$, i.e.
\[
\int_{\Cpx^n}\sigma(z_1,\ldots,z_n)\dif m(z_1)\cdots\dif m(z_n)=1\;.
\]
For $p\in\{1,\ldots,n\}$, let $\sigma^{(p)}$ be the marginal density of the first $p$ coordinates, i.e.
\[
\sigma^{(p)}(u_1,\ldots,u_p)=\int_{\Cpx^{n-p}}\sigma(u_1,\ldots,u_p,z_{p+1},\ldots,z_n)\dif m(z_{p+1})\cdots\dif m(z_n)
\]
for $p<n$ and $\sigma^{(n)}=\sigma$.
Then
\[
\sigma^{(p)}(u_1,\ldots,u_p)\le\frac{(n-p)!}{n!}n^p\sigma^{(1)}(u_1)\cdots\sigma^{(1)}(u_p)
\]
for all $p\in\{1,\ldots,n\}$ and $u_1,\ldots,u_p\in\Cpx$.
\end{lemma}
\begin{proof}
The proof follows standard methods in random matrix theory, see~\cite[\S5.2]{Mehta}
for the real variable case and~\cite[\S15.1]{Mehta} for the complex variable case with $g(z)=e^{-|z|^2}$.
Let $(P_j(z))_{j=0}^\infty$ be the sequence of orthonormal polynomials obtained from the Gram--Schmidt orthonomalization
process applied to $1,z,z^2,\ldots$ in $L^2(\Cpx,g\,\dif m)$, and put
\[
\phi_j(z)=\sqrt{g(z)}P_j(z)\qquad j=0,1,\ldots\;.
\]
Then $(\phi_j(z))_{j=0}^\infty$ is an orthonormal sequence in $L^2(\Cpx,\dif m)$.
Using the Vandermonde determinant, one finds
\[
\sigma(z_1,\ldots,z_n)=\frac1{n!}\,\big|\det\big((\phi_{i-1}(z_j))_{i,j=1}^n\big)\big|^2
=\frac1{n!}\,\bigg|\sum_{\pi\in S_n}\sign(\pi)\prod_{j=1}^n\phi_{\pi(j)-1}(z_j)\bigg|^2\;,
\]
where $S_n$ is the permitation group of $\{1,\ldots,n\}$.
The normalization constant $\frac1{n!}$ can be determined from the fact that $\sigma$ is a probability density.
Now put
\begin{equation}
\label{eq:psiuv}
\psi(u,v)=\sum_{j=0}^{n-1}\phi_j(u)\overline{\phi_j(v)},\qquad(u,v\in\Cpx)\;.
\end{equation}
Then as in~\cite[pp.\ 80, 91-92]{Mehta}, one finds 
\[
\sigma^{(p)}(u_1,\ldots,u_p)=\frac{(n-p)!}{n!}\det\big((\psi(u_i,u_j))_{i,j=1}^p\big)\;.
\]
In particular,
\[
\sigma^{(1)}(u)=\frac1n\psi(u,u)\;.
\]
By~\eqref{eq:psiuv}, $(\psi(u_i,u_j))_{i,j=1}^p$ is a positive semidefinite matrix.
Hence by Theorem~\ref{thm:Hadamard},
\[
\sigma^{(p)}(u_1,\ldots,u_p)\le\frac{(n-p)!}{n!}\prod_{i=1}^p\psi(u_i,u_i)
=\frac{(n-p)!}{n!}n^p\prod_{i=1}^p\sigma^{(1)}(u_i)\;.
\]
\end{proof}

\begin{lemma}
\label{lem:ellipticmarg}
For each $n\in\Nats$, let $\nu_n$ denote the probability measure on $\Cpx^n$ with density given
by~\eqref{eq:Dthetadist} in Theorem~\ref{thm:rmelliptic}, and for $p\in\{1,\ldots,n\}$ let
$\nu_n^{(p)}$ be the joint marginal distribution of $\nu_n$ on the first $p$ coordinates.
Let $\mu$ be the uniform distribution on the solid ellipse
\[
\bigg\{z\in\Cpx\,\bigg|\,\frac{(\RealPart z)^2}{4\cos^4\theta}+\frac{(\ImagPart z)^2}{4\sin^4\theta}
\le1\bigg\}\;.
\]
Then for every $p\in\Nats$, $\nu_n^{(p)}$ converges to $\btimes_1^p\mu$ as $n\to\infty$ both in $*$--moments
and in the w$^*$--topology on $\Prob(\Cpx^p)$.
\end{lemma}
\begin{proof}
The case $p=1$ is proved in~\cite{HL}.
Consider now the case $p\ge2$.
By Lemma~\ref{lem:margdistinequal},
\[
\nu_n^{(p)}\le\frac{(n-p)!}{n!}n^p\nu_n^{(1)}\times\cdots\times\nu_n^{(1)},\qquad(n\ge p)\;.
\]
Hence, as in the proof of Corollary~\ref{cor:nonrandD}, we get
\[
\|\nu_n^{(p)}-\btimes_1^p\nu_n^{(1)}\|\to0\quad\text{as }n\to\infty\;,
\]
so in particular
\[
\lim_{n\to\infty}\nu_n^{(p)}=\lim_{n\to\infty}\btimes_1^p\nu_n^{(1)}=\btimes_1^p\mu
\]
in the weak$^*$--topology on $\Prob(\Cpx^p)$.
However, $\nu_n^{(p)}$ is not compactly supported, so we must argue further to prove convergence in
$*$--moments.
Clearly, it suffices to consider convergence on real valued polynomials in $z_1,\ldots,z_p,\overline{z_1},\ldots,\overline{z_p}$.
Let $h:\Cpx^p\to\Reals$ be such a polynomial.
Choose $C>0$ and $d\in\Nats$ such that
\[
|h(z_1,\ldots,z_p)|\le C(1+|z_1|^2+\cdots+|z_p|^2)^d\;.
\]
Then 
\[
|h(z_1,\ldots,z_p)|\le C\bigg(\sum_{i=1}^p(1+|z_i|^2)\bigg)^d
\le Cp^{d-1}\sum_{i=1}^p(1+|z_i|^2)^d\;,
\]
where the last inequality follows from the convexity of $x\mapsto x^d$ on $[0,\infty)$.
Put
\[
g(z_1,\ldots,z_p)=Cp^{d-1}\sum_{i=1}^p(1+|z_i|^2)^d\;.
\]
Then $g+h$ and $g-h$ are non--negative, continuous functions of $\Cpx^p$.
By appoximating $g\pm h$ from below with functions in $C_c(\Cpx^p)_+$, it follows from the weak$^*$--convergence
of $\nu_n^{(p)}$ to $\btimes_1^p\mu$, that
\begin{equation}
\label{eq:gpmh}
\int_{\Cpx^p}(g\pm h)\dif(\btimes_1^p\mu)\le\liminf_{n\to\infty}\int_{\Cpx^p}(g\pm h)\dif\nu_n^{(p)}\;.
\end{equation}
Note that $g$ is of the form $g(z_1,\ldots,z_p)=\sum_{i=1}^pf(z_i)$,
where $f:\Cpx\to\Reals$ is a polynomial in $z$ and $\overline z$.
Moreover, the one--dimensional marginal distributions of $\nu_n^{(p)}$ are all equal to $\nu_n^{(1)}$.
Using that $\nu_n^{(1)}$ converges to $\nu$ in $*$--moments, we therefore get
\[
\int_{\Cpx^p}g\dif(\btimes_1^p\mu)=p\int_\Cpx f\dif\mu=\lim_{n\to\infty}\bigg(p\int_\Cpx f\dif\nu_n^{(1)}\bigg)
=\lim_{n\to\infty}\bigg(\int_\Cpx g\dif\nu_n^{(p)}\bigg)\;.
\]
Hence by~\eqref{eq:gpmh}, we have
\[
\int_{\Cpx^p}h\,\dif(\btimes_1^p\mu)\le\liminf_{n\to\infty}\int_{\Cpx^p}h\,\dif\nu_n^{(p)}
\le\limsup_{n\to\infty}\int_{\Cpx^p}h\,\dif\nu_n^{(p)}\le\int_{\Cpx^p}h\,\dif(\btimes_1^p\mu)\;.
\]
This proves the convergence of $\nu_n^{(p)}$ to $\btimes_1^p\mu$ in $*$--moments.
\end{proof}

\begin{proof}[Proof of Theorem~\ref{thm:elliptic}]
In light of Proposition~\ref{prop:constadj}(i), we may without loss of generality
assume $a^2+b^2=1$, say $a=\cos\theta$ and $b=\sin\theta$, some $0<\theta<\frac\pi2$.
By foundational results of Voiculescu~\cite{V:RM},
$Y_\theta(n)$ from Theorem~\ref{thm:rmelliptic} converges in $*$--moments as $n\to\infty$ to the elliptic element
$(\cos\theta)x_1+i(\sin\theta)x_2$.
By Theorem~\ref{thm:rmelliptic}, also $D_\theta(n)+\sin(2\theta)T(n)$ converges in $*$--moments
as $n\to\infty$ to this elliptic element.
Let $\nu_n$ be the joint distribution of $n$ complex variables having
density~\eqref{eq:Dthetadist} with respect to Lebesgue measure.
In Lemma~\ref{lem:ellipticmarg}, it is shown that for every $p\in\Nats$, the $p$th marginal
distribution $\nu_n^{(p)}$ of $\nu_n$ converges in $*$--moments as $n\to\infty$
to the $p$--fold Cartesian product of the uniform distribution on the solid ellipse
\[
\bigg\{z\in\Cpx\,\bigg|\,\frac{(\RealPart z)^2}{4\cos^4\theta}+\frac{(\ImagPart z)^2}{4\sin^4\theta}\le1\bigg\}.
\]
Consequently, the elliptic element $ax_1+ibx_2$ is a $\DT(\nu_{a,b},\sin(2\theta))$--element.
\end{proof}

\section{DT--operators in finite von Neumann algebras}
\label{sec:DTinLF2}

Our first task in this section is to show (Theorem~\ref{thm:DTinLF2}) how every DT--element can be constructed
in the von Neumann algebra $L(\Fb_2)$, (or more precisely, in the $*$--noncommutative probability space
$(L(\Fb_2),\tau)$, where $\tau$ is the tracial state on $L(\Fb_2)$), as $D+cT$ where $T$ is ``the
upper triangular part'' of a semicircular element $X$ and where $D$ is from a ``diagonal'' algebra free from $X$.
Then we prove a series of results culminating in what may be called our main technical theorem
(Theorem~\ref{thm:DTfreematrix}), which shows that every DT--element can be realized as an upper triangular matrix
having mutually free entries that are themselves circular elements or DT--elements.
This theorem, as well as being pretty, is the main technical tool for proving decomposability of DT-operators
in~\S\ref{sec:decomp}.

Let $\Mcal$ be a von Neumann algebra with a normal faithful state $\tau$,
let $\lambda:L^\infty([0,1])\to\Mcal$ be a normal, unital, injective $*$--homomorphism and let
$X\in\Mcal$ be a centered semicircular element with $\tau(X^2)=\nu^2>0$ and such that $X$ and the image of
$\lambda$ are free with respect to $\tau$.
(Thus the subalgebra of $\Mcal$ generated by $X$ and the image of $\lambda$ is isomorphic to $L(\Fb_2)$, to which
the restriction of $\tau$ is a trace.)
For $0\le a<b\le1$ let $p[a,b]=\lambda(1_{[a,b]})\in\Mcal$.

\begin{lemma}
\label{lem:Tconstr}
For $\ell\in\Nats$ let
\[
T_\ell=\sum_{j=1}^{2^\ell-1}p[\tfrac{j-1}{2^\ell},\tfrac j{2^\ell}]Xp[\tfrac j{2^\ell},1]\;.
\]
Then $T_\ell$ converges in norm as $\ell\to\infty$ to an element $T\in\Mcal$, satisfying $T+T^*=X$.
\end{lemma}
\begin{proof}
For example, $T_2$ and $T_3$ are represented by the shaded regions below.
\begin{center}
\begin{picture}(200,100)(0,-20)


  \linethickness{0.7pt}
  \drawline(0,0)(0,80)(80,80)(80,0)(0,0)

  \linethickness{0.5pt}
  \drawline(20,80)(20,60)(40,60)(40,40)(60,40)(60,20)(80,20)

  \linethickness{0.3pt}

  \drawline(22,80)(22,60)
  \drawline(24,80)(24,60)
  \drawline(26,80)(26,60)
  \drawline(28,80)(28,60)
  \drawline(30,80)(30,60)
  \drawline(32,80)(32,60)
  \drawline(34,80)(34,60)
  \drawline(36,80)(36,60)
  \drawline(38,80)(38,60)
  \drawline(40,80)(40,60)

  \drawline(20,78)(80,78)
  \drawline(20,76)(80,76)
  \drawline(20,74)(80,74)
  \drawline(20,72)(80,72)
  \drawline(20,70)(80,70)
  \drawline(20,68)(80,68)
  \drawline(20,66)(80,66)
  \drawline(20,64)(80,64)
  \drawline(20,62)(80,62)

  \drawline(42,80)(42,40)
  \drawline(44,80)(44,40)
  \drawline(46,80)(46,40)
  \drawline(48,80)(48,40)
  \drawline(50,80)(50,40)
  \drawline(52,80)(52,40)
  \drawline(54,80)(54,40)
  \drawline(56,80)(56,40)
  \drawline(58,80)(58,40)
  \drawline(60,80)(60,40)

  \drawline(40,60)(80,60)
  \drawline(40,58)(80,58)
  \drawline(40,56)(80,56)
  \drawline(40,54)(80,54)
  \drawline(40,52)(80,52)
  \drawline(40,50)(80,50)
  \drawline(40,48)(80,48)
  \drawline(40,46)(80,46)
  \drawline(40,44)(80,44)
  \drawline(40,42)(80,42)

  \drawline(62,80)(62,20)
  \drawline(64,80)(64,20)
  \drawline(66,80)(66,20)
  \drawline(68,80)(68,20)
  \drawline(70,80)(70,20)
  \drawline(72,80)(72,20)
  \drawline(74,80)(74,20)
  \drawline(76,80)(76,20)
  \drawline(78,80)(78,20)

  \drawline(60,40)(80,40)
  \drawline(60,38)(80,38)
  \drawline(60,36)(80,36)
  \drawline(60,34)(80,34)
  \drawline(60,32)(80,32)
  \drawline(60,30)(80,30)
  \drawline(60,28)(80,28)
  \drawline(60,26)(80,26)
  \drawline(60,24)(80,24)
  \drawline(60,22)(80,22)

  \put(35,-20){$T_2$}


  \linethickness{0.7pt}
  \drawline(120,0)(120,80)(200,80)(200,0)(120,0)

  \linethickness{0.5pt}
  \drawline(130,80)(130,70)(140,70)(140,60)(150,60)(150,50)(160,50)
           (160,40)(170,40)(170,30)(180,30)(180,20)(190,20)(190,10)(200,10)

  \linethickness{0.3pt}

  \drawline(130,78)(200,78)
  \drawline(130,76)(200,76)
  \drawline(130,74)(200,74)
  \drawline(130,72)(200,72)

  \drawline(132,80)(132,70)
  \drawline(134,80)(134,70)
  \drawline(136,80)(136,70)
  \drawline(138,80)(138,70)
  \drawline(140,80)(140,70)

  \drawline(140,70)(200,70)
  \drawline(140,68)(200,68)
  \drawline(140,66)(200,66)
  \drawline(140,64)(200,64)
  \drawline(140,62)(200,62)

  \drawline(142,80)(142,60)
  \drawline(144,80)(144,60)
  \drawline(146,80)(146,60)
  \drawline(148,80)(148,60)
  \drawline(150,80)(150,60)

  \drawline(150,60)(200,60)
  \drawline(150,58)(200,58)
  \drawline(150,56)(200,56)
  \drawline(150,54)(200,54)
  \drawline(150,52)(200,52)

  \drawline(152,80)(152,50)
  \drawline(154,80)(154,50)
  \drawline(156,80)(156,50)
  \drawline(158,80)(158,50)
  \drawline(160,80)(160,50)

  \drawline(160,50)(200,50)
  \drawline(160,48)(200,48)
  \drawline(160,46)(200,46)
  \drawline(160,44)(200,44)
  \drawline(160,42)(200,42)

  \drawline(162,80)(162,40)
  \drawline(164,80)(164,40)
  \drawline(166,80)(166,40)
  \drawline(168,80)(168,40)
  \drawline(170,80)(170,40)

  \drawline(170,40)(200,40)
  \drawline(170,38)(200,38)
  \drawline(170,36)(200,36)
  \drawline(170,34)(200,34)
  \drawline(170,32)(200,32)

  \drawline(172,80)(172,30)
  \drawline(174,80)(174,30)
  \drawline(176,80)(176,30)
  \drawline(178,80)(178,30)
  \drawline(180,80)(180,30)

  \drawline(180,30)(200,30)
  \drawline(180,28)(200,28)
  \drawline(180,26)(200,26)
  \drawline(180,24)(200,24)
  \drawline(180,22)(200,22)

  \drawline(182,80)(182,20)
  \drawline(184,80)(184,20)
  \drawline(186,80)(186,20)
  \drawline(188,80)(188,20)
  \drawline(190,80)(190,20)

  \drawline(190,20)(200,20)
  \drawline(190,18)(200,18)
  \drawline(190,16)(200,16)
  \drawline(190,14)(200,14)
  \drawline(190,12)(200,12)

  \drawline(192,80)(192,10)
  \drawline(194,80)(194,10)
  \drawline(196,80)(196,10)
  \drawline(198,80)(198,10)

  \put(155,-20){$T_3$}

\end{picture}
\end{center}
We have
\[
|T_{\ell+1}-T_\ell|^2=\sum_{j=1}^{2^\ell}p[\tfrac{2j-1}{2^{\ell+1}},\tfrac{2j}{2^{\ell+1}}]
Xp[\tfrac{2j-2}{2^{\ell+1}},\tfrac{2j-1}{2^{\ell+1}}]
Xp[\tfrac{2j-1}{2^{\ell+1}},\tfrac{2j}{2^{\ell+1}}]\;.
\]
By results of Voiculescu~\cite{V:RM}, 
\[
p[\tfrac{2j-1}{2^{\ell+1}},\tfrac{2j}{2^{\ell+1}}]
Xp[\tfrac{2j-2}{2^{\ell+1}},\tfrac{2j-1}{2^{\ell+1}}]
Xp[\tfrac{2j-1}{2^{\ell+1}},\tfrac{2j}{2^{\ell+1}}]
\]
is the square of a circular element of norm $2^{(1-\ell)/2}\nu$
in $p[\tfrac{2j-1}{2^{\ell+1}},\tfrac{2j}{2^{\ell+1}}]\Mcal p[\tfrac{2j-1}{2^{\ell+1}},\tfrac{2j}{2^{\ell+1}}]$
with respect to the renormalization of $\tau$.
Thus $\|T_{\ell+1}-T_\ell\|=2^{(1-\ell)/2}\nu$
and $T_\ell$ converges in norm as $\ell\to\infty$.

In a similar manner, since 
\[
X-T_\ell-T_\ell^*=\sum_{j=1}^{2^\ell}p[\tfrac{j-1}{2^\ell},\tfrac j{2^\ell}]Xp[\tfrac{j-1}{2^\ell},\tfrac j{2^\ell}]
\]
and $p[\tfrac{j-1}{2^\ell},\tfrac j{2^\ell}]Xp[\tfrac{j-1}{2^\ell},\tfrac j{2^\ell}]$
is a semicircular element of norm $2^{(2-\ell)/2}\nu$, we find $T+T^*=X$.
\end{proof}

\begin{defi}\rm
\label{def:UT}
The element $T$ constructed in Lemma~\ref{lem:Tconstr} from $X$ and $\lambda$ will be denoted
$T=\UT(X,\lambda)$.
\end{defi}

Note that for any nonzero real number $t$, $\UT(tX,\lambda)=|t|\UT(X,\lambda)$.

For the remainder of this section, we will let $X$ and $\lambda$ be as described above,
with the added convention that the second moment of $X$ will be $\nu^2=1$.

Given $n\in\Nats$ and $0\le a<b\le1$, let
\[
P_n[a,b]=\diag(\overset{[nb]}{\overbrace{\underset{[na]}{\underbrace{0,\ldots,0}},1,\ldots,1}},0,\ldots,0)\in M_n(\Cpx)\;.
\]

\begin{lemma}
\label{lem:rmDT}
Let $T=\UT(X,\lambda)$.
Let $D_s\in\lambda(L^\infty[0,1])$ for all $s$ in some set $I$.
Suppose that $D_s^{(k)}\in M_{2^k}(\Cpx)$ is a diagonal matrix ($s\in I$, $k\in\Nats$), such that for all $s\in I$,
$\|D_s^{(k)}\|$ remains bounded as $k\to\infty$ and such that the family
\[
\big(D_s^{(k)}\big)_{s\in I}\;,\quad\big(P_{2^k}[\tfrac i{2^\ell},\tfrac j{2^\ell}]\big)_{\ell\in\Nats,\,i,j\in\{0,1,\ldots,2^\ell\}}
\]
converges in $*$--moments as $k\to\infty$ to the family 
\[
\big(D_s\big)_{s\in I}\;,\quad\big(p[\tfrac i{2^\ell},\tfrac j{2^\ell}]\big)_{\ell\in\Nats,\,i,j\in\{0,1,\ldots,2^\ell\}}\;.
\]
Let $T^{(k)}\in\UTGRM(2^k,2^{-k})$.
Then the family
\begin{equation}
\label{eq:TkDk}
T^{(k)},\quad\big(D_s^{(k)}\big)_{s\in I}
\end{equation}
converges in $*$--moments as $k\to\infty$ to the family
\begin{equation}
\label{eq:TD}
T,\quad\big(D_s\big)_{s\in I}\;.
\end{equation}
\end{lemma}
\begin{proof}
We may assume the family $(D_s)_{s\in I}$ contains the identity and is closed under taking adjoints and under multiplication.
Let $X^{(k)}\in\SGRM(2^k,2^{-k})$.
By results of Voiculescu~\cite{V:RM}, the family
\[
X^{(k)},\quad\big(D_s^{(k)}\big)_{s\in I}\;,\quad\big(P_{2^k}[\tfrac i{2^\ell},\tfrac j{2^\ell}]\big)_{\ell\in\Nats,\,i,j\in\{0,1,\ldots,2^\ell\}}
\]
converges in $*$--moments as $k\to\infty$ to the family
\[
X,\quad\big(D_s\big)_{s\in I}\;,\quad\big(p[\tfrac i{2^\ell},\tfrac j{2^\ell}]\big)_{\ell\in\Nats,\,i,j\in\{0,1,\ldots,2^\ell\}}\;.
\]
Therefore, letting
\[
T_\ell^{(k)}=\sum_{j=1}^{2^\ell-1}P_{2^k}[\tfrac{j-1}{2^\ell},\tfrac j{2^\ell}]X^{(k)}P_{2^k}[\tfrac j{2^\ell},1]\;,
\]
the family
\begin{equation}
\label{eq:TlkDk}
T_\ell^{(k)},\quad\big(D_s^{(k)}\big)_{s\in I}
\end{equation}
converges in $*$--moments as $k\to\infty$ to the family
\begin{equation}
\label{eq:TlD}
T_\ell,\quad\big(D_s\big)_{s\in I}\;.
\end{equation}
Let $\delta_\ell T^{(k)}=T^{(k)}-T_\ell^{(k)}$ and for $q\in\{1,2,3,4\}$ let
\[
A_\ell^{(k)}(q)=\begin{cases}
T_\ell^{(k)}&\text{if }q=1, \\
(T_\ell^{(k)})^*&\text{if }q=2, \\
\delta_\ell T^{(k)}&\text{if }q=3, \\
(\delta_\ell T^{(k)})^*&\text{if }q=4.
\end{cases}
\]
We will show that if $m\in\Nats$, if $q_1,\ldots,q_m\in\{1,2,3,4\}$ with $q_j\in\{3,4\}$ for at least one
$j\in\{1,\ldots,m\}$ and if $s_1,\ldots,s_m\in I$, then there is a constant $C>0$, independent of $k$ and $\ell$,
such that
\begin{equation}
\label{eq:AqD}
\big|\tau_{2^k}\big(A_\ell^{(k)}(q_1)D_{s_1}^{(k)}\cdots A_\ell^{(k)}(q_m)D_{s_m}^{(k)}\big)\big|\le C2^{-\ell}
\end{equation}
for all $k\ge\ell\ge1$.
This estimate will finish the proof, because it will imply that for every $\eps(1),\ldots,\eps(m)\in\{*,1\}$, there is
a constant $C'$, independent of $\ell$, such that 
\[
\big|\lim_{k\to\infty}\tau_{2^k}\big((T^{(k)})^{\eps(1)}D_{s_1}^{(k)}\cdots (T^{(k)})^{\eps(m)}D_{s_m}^{(k)}
-(T_\ell^{(k)})^{\eps(1)}D_{s_1}^{(k)}\cdots (T_\ell^{(k)})^{\eps(m)}D_{s_m}^{(k)}
\big)\big|\le C'2^{-\ell}\;.
\]
Combined with the convergence of~\eqref{eq:TlkDk} to~\eqref{eq:TlD}
and the norm convergence of $T_\ell$ to $T$ as $\ell\to\infty$,
this will imply convergence in $*$--moments of the family~\eqref{eq:TkDk}
to the family~\eqref{eq:TD}.

To prove the upper bound~\eqref{eq:AqD}, we will make an analysis using graphs and pairings as did Voiculescu in~\cite{V:RM}
and we will use the fact that the proportion of nonzero entries in $\delta_\ell T^{(k)}$ is on the order of $2^{-\ell}$.
Fixing $\ell$, we have
\begin{align*}
&\tau_{2^k}\big(A_\ell^{(k)}(q_1)D_{s_1}^{(k)}\cdots A_\ell^{(k)}(q_m)D_{s_m}^{(k)}\big) \\
&=2^{-k}\sum_{i_1,\ldots,i_m\in\{1,\ldots,2^k\}}
\begin{aligned}[t]
\Eb(b^{(k)}(q_1;i_0,i_1)d^{(k)}(s_1;i_1)&b^{(k)}(q_2;i_1,i_2)d^{(k)}(s_2;i_2)\cdots \\
\cdots&b^{(k)}(q_m;i_{m-1},i_m)d^{(k)}(s_m;i_m))\;,
\end{aligned}
\end{align*}
where we use the convention $i_0=i_m$, where $d^{(k)}(s;i)$ is the $i$th diagonal entry of $D_s^{(k)}$ and
where $b(q;i,j)$ is the $(i,j)$th entry of $A_\ell^{(k)}(q)$.
Let $K=\prod_{j=1}^m\sup_{k\ge1}\|D_{s_j}^{(k)}\|$.
Then
\begin{align*}
\big|\Eb(b^{(k)}(q_1;i_0,i_1)&d^{(k)}(s_1;i_1)b^{(k)}(q_2;i_1,i_2)d^{(k)}(s_2;i_2)\cdots
b^{(k)}(q_m;i_{m-1},i_m)d^{(k)}(s_m;i_m))\big|\le \\[0.5ex]
&\le K|\Eb(b^{(k)}(q_1;i_0,i_1)b^{(k)}(q_2;i_1,i_2)\cdots b^{(k)}(q_m;i_{m-1},i_m))|\;.
\end{align*}
By a generalized H\"older inequality,
\begin{equation}
\label{eq:Ebb}
|\Eb(b^{(k)}(q_1;i_0,i_1)b^{(k)}(q_2;i_1,i_2)\cdots b^{(k)}(q_m;i_{m-1},i_m))|\le2^{-km/2}\;.
\end{equation}
Moreover, the LHS of~\eqref{eq:Ebb} is nonzero only if there is a pairing $\sigma$ of $\{1,2,\ldots,m\}$
such that $\{r,s\}\in\sigma$ implies
\renewcommand{\labelenumi}{(\roman{enumi})}
\begin{enumerate}
\item
$i_{r-1}=i_s$ and $i_r=i_{s-1}$,
\item
either $\{q_r,q_s\}=\{1,2\}$ or $\{q_r,q_s\}=\{3,4\}$.
\end{enumerate}
Therefore, an upper bound for the LHS of~\eqref{eq:AqD} is
\begin{equation}
\label{eq:KNsig}
K2^{-k(1+m/2)}\sum_\sigma N(\sigma)\;,
\end{equation}
where the sum is over all pairings $\sigma$ satisfying~(ii) for every $\{r,s\}\in\sigma$ and where $N(\sigma)$ is the number
of choices of $i_1,i_2,\ldots,i_m\in\{1,2,\ldots,2^k\}$ so that~(i) holds for every $\{r,s\}\in\sigma$, (with $i_0=i_m$)
and so that $b^{(k)}(q_j;i_{j-1},i_j)\ne0$ for every $j\in\{1,\ldots,m\}$.
Choosing $i_1,\ldots,i_m$ such that~(i) holds is equivalent to taking the quotient graph $Q$ of the $m$--gon graph $G$ according
to $\sigma$ as in Algorithm~\ref{alg:NTO} part~(B), (where of course $k$ over there is $m$ here),
and assigning values from $\{1,\ldots,2^k\}$ to the vertices of $Q$
(thus assigning values to $i_1,\ldots,i_m$ according to the values assigned the images in $Q$ of the corresponding vertices
$v_1,\ldots,v_m$ in $G$).
If $\sigma$ is crossing, then from Lemma~\ref{lem:nctree}, $Q$ has $\le m/2$ vertices and hence $N(\sigma)\le2^{km/2}$.
If $\sigma$ is non--crossing, then $Q$ has exactly $1+m/2$ vertices;
however, examining Lemma~\ref{lem:prdnbs}, we find that in this case
there is $j$ with $q_j\in\{3,4\}$ and with $i_{j-1}$ and $i_j$ mapped to distinct vertices of $Q$.
Since for a given value of $i_j$ there are at most $2^{k-\ell}$ values of $i_{j-1}$ making $b^{(k)}(q_j;i_{j-1},i_j)$
nonzero, we find that $N(\sigma)\le K2^{k(1+m/2)-\ell}$.
Hence using~\eqref{eq:KNsig}, when $k\ge\ell$ we get an upper bound of the form~\eqref{eq:AqD}.
\end{proof}

\begin{thm}
\label{thm:DTinLF2}
Let $T=\UT(X,\lambda)$, (where $X$ has second moment $1$) and let $c>0$.
Take $f\in L^\infty[0,1]$ and let $D=\lambda(f)$.
Then $D+cT$ is a $\DT(\mu,c)$--element, where $\mu$ is the push--forward measure
of Lebesgue measure by $f$.
\end{thm}
\begin{proof}
One easily finds diagonal matrices $D^{(k)}\in M_{2^k}(\Cpx)$ such that the family
\[
D^{(k)},\quad\big(P_{2^k}[\tfrac i{2^\ell},\tfrac j{2^\ell}]\big)_{\ell\in\Nats,\,i,j\in\{0,1,\ldots,2^\ell\}}
\]
converges in $*$--moments as $k\to\infty$ to the family 
\[
D,\quad\big(p[\tfrac i{2^\ell},\tfrac j{2^\ell}]\big)_{\ell\in\Nats,\,i,j\in\{0,1,\ldots,2^\ell\}}\;.
\]
Let $T^{(k)}\in\UTGRM(2^k,2^{-k})$.
Then by Lemma~\ref{lem:rmDT}, $D^{(k)}+cT^{(k)}$ converges in $*$--moments to $D+cT$.
On the other hand, by Corollary~\ref{cor:nonrandD}, $D^{(k)}+cT^{(k)}$ converges in $*$--moments
to a $\DT(\mu,c)$--element.
\end{proof}

\begin{defi}\rm
A $\DT(\mu,c)$--{\em operator}\/ is a $\DT(\mu,c)$--element in a W$^*$--noncommutative probability space $(\Mcal,\phi)$
where $\phi$ is faithful.
A {\em DT--operator} is an element that is a $\DT(\mu,c)$--operator for some $\mu$ and $c$.
\end{defi}

Theorem~\ref{thm:DTinLF2} shows that for every pair $(\mu,c)$ there exist $\DT(\mu,c)$--operators.

\begin{lemma}
\label{lem:Tbasics}
If $T=\UT(X,\lambda)$ and $0<t<1$ then
\begin{align}
Tp[0,t]&=p[0,t]Tp[0,t] \label{eq:Tp} \\
p[0,t]Tp[t,1]&=p[0,t]Xp[t,1]\;. \label{eq:pTp}
\end{align}
\end{lemma}
\begin{proof}
The identity~\eqref{eq:Tp} clearly holds when $T$ is replaced by $T_\ell$.
Taking the limit as $\ell\to\infty$ proves~\eqref{eq:Tp}.
Then the identity~\eqref{eq:pTp} follows from~\eqref{eq:Tp} and $X=T+T^*$.
\end{proof}

\begin{lemma}
\label{lem:iT}
Let $T=\UT(X,\lambda)$ and let $X'=iT^*-iT$.
Then $X'$ is a semicircular element of second moment $1$, and $X'$ and the image of $\lambda$ form a free pair.
Furthermore,
\begin{equation}
\label{eq:iT}
-iT=\UT(X',\lambda)\;.
\end{equation}
\end{lemma}
\begin{proof}
Using Lemma~\ref{lem:rmDT} and taking, for example, the family $(D_s^{(k)})_{s\in I}$ to be the family
$(P_{2^k}[\tfrac i{2^\ell},\tfrac j{2^\ell}])_{\ell\in\Nats,\,i,j\in\{0,1,\ldots,2^\ell\}}$ itself,
converging to the family
\[
(D_s)_{s\in I}=(p[\tfrac i{2^\ell},\tfrac j{2^\ell}])_{\ell\in\Nats,\,i,j\in\{0,1,\ldots,2^\ell\}}\;,
\]
we find that the family
\[
i(T^{(k)})^*-iT^{(k)},\quad(P_{2^k}[\tfrac i{2^\ell},\tfrac j{2^\ell}])_{\ell\in\Nats,\,i,j\in\{0,1,\ldots,2^\ell\}}
\]
converges in $*$--moments to
\[
X',\quad(p[\tfrac i{2^\ell},\tfrac j{2^\ell}])_{\ell\in\Nats,\,i,j\in\{0,1,\ldots,2^\ell\}}\;,
\]
where $T^{(k)}\in\UTGRM(2^k,2^{-k})$.
By Voiculescu's matrix model~\cite{V:RM}, it follows that $X'$ is a semicircular element of second moment $1$ and $X'$
and the image of $\lambda$ are free.

If $0<t<1$, from Lemma~\ref{lem:Tbasics} we get
\[
p[0,t]X'p[t,1]=p[0,t](-iT)p[t,1]=-ip[0,t]Xp[t,1]\;.
\]
Hence
\[
\sum_{j=1}^{2^\ell-1}p[\tfrac{j-1}{2^\ell},\tfrac j{2^\ell}]X'p[\tfrac j{2^\ell},1]
=-iT_\ell\;.
\]
Letting $\ell\to\infty$ yields~\eqref{eq:iT}.
\end{proof}

\begin{remark}\rm
\label{rem:normT}
What we have shown implies that if $T=\UT(X,\lambda)$ then $\|T\|\le2$, because
$\RealPart\,T=X/2$ and $\ImagPart\,T=X'/2$ both have norm $1$.
In~\S\ref{sec:TstT}, we will show that actually $\|T\|=\sqrt e$.
\end{remark}

\begin{lemma}
\label{lem:Xpieces}
Let $m\in\Nats$, $0=s_0<s_1<\cdots<s_m=1$ and
\[
S=\sum_{j=1}^mp[s_{j-1}s_j]Xp[s_j,1]\;.
\]
Then with $T=\UT(X,\lambda)$, we have $\|T-S\|\le2\max_{1\le j\le m}(s_j-s_{j-1})^{1/2}$.
\end{lemma}
\begin{proof}
We have
\[
2\,\RealPart(T-S)=X-(S+S^*)=\sum_{j=1}^mp[s_{j-1},s_j]Xp[s_{j-1},s_j]\;.
\]
By results of Voiculescu~\cite{V:RM},
$\|p[s_{j-1},s_j]Xp[s_{j-1},s_j]\|=2(s_j-s_{j-1})^{1/2}$.
Hence
\[
\|\RealPart(T-S)\|=\max_{1\le j\le m}(s_j-s_{j-1})^{1/2}\;.
\]
On the other hand, we have
$2\,\ImagPart(T-S)=X'-(iS^*-iS)$, where $X'=iT^*-iT$.
Appealing to Lemmas~\ref{lem:Tbasics} and~\ref{lem:iT},
\[
iS=i\sum_{j=1}^mp[s_{j-1},s_j]Tp[s_j,1]
=\sum_{j=1}^mp[s_{j-1},s_j]X'p[s_j,1]
\]
and $\|X'-(iS^*-iS)\|=2\max_{1\le j\le m}(s_j-s_{j-1})^{1/2}$.
\end{proof}

\begin{lemma}
\label{lem:qTq}
Let $T=\UT(X,\lambda)$.
Suppose $B\subseteq[0,1]$ is a Borel set of nonzero Lebesgue measure.
Let $\alpha:[0,1]\to[0,1]$ be defined by
\[
\alpha(t)=\frac{\rho([0,t)\cap B)}{\rho(B)}\;,
\]
where $\rho$ denotes Lebesgue measure on $[0,1]$.
Let $q=\lambda(1_B)$ and let
\[
\lambdat:L^\infty[0,1]\to q\lambda(L^\infty[0,1])\cong L^\infty(B)
\]
be 
\[
\lambdat(f)=q\lambda(f\circ\alpha)\;.
\]
Then $\lambdat$ is an injective, normal $*$--homomorphism and
\[
qTq=\UT(qXq,\lambdat)\;.
\]
\end{lemma}
\begin{proof}
The map $\alpha$ is monotone and continuous, hence measureable;
therefore, $\lambdat$ is a normal $*$--homomorphism.
For any measureable set $A\subseteq[0,1]$, we have $\rho(\alpha^{-1}(A)\cap B)=\rho(A)\rho(B)$.
Hence if $\rho(A)>0$, then $\lambdat(1_A)=\lambda(1_{\alpha^{-1}(A)\cap B})\ne0$
and $\lambdat$ is injective.
If $0\le s<t\le 1$, then up to measure zero,
$[s,t]\cap B=\alpha^{-1}([\alpha(s),\alpha(t)])\cap B$, so
\[
\lambda(1_{[s,t]})q=\lambda(1_{[\alpha(s),\alpha(t)]}\circ\alpha)q=\lambdat(1_{[\alpha(s),\alpha(t)]})\;.
\]
Writing $p[s,t]=\lambda(1_{[s,t]})$ as usual and letting $\pt[a,b]=\lambdat(1_{[a,b]})$, we thus have
$qp[s,t]=\pt[\alpha(s),\alpha(t)]$ and
\[
qT_\ell q=\sum_{j=1}^{2^\ell-1}qp[\tfrac{j-1}{2^\ell},\tfrac j{2^\ell}]Xp[\tfrac j{2^\ell},1]q
=\sum_{j=1}^{2^\ell-1}\pt[\alpha(\tfrac{j-1}{2^\ell}),\alpha(\tfrac j{2^\ell})](qXq)\pt[\alpha(\tfrac j{2^\ell}),1]\;.
\]
By Voiculescu's results~\cite{V:RM}, we have that, with respect to the trace $\tau(q)^{-1}\tau\restrict_{q\Mcal q}$,
$qXq$ is a semicircular element with second moment $\tau(q)$ and is free from the image of $\lambdat$.
Let $\Tt=\UT(qXq,\lambdat)$.
Since
\[
\alpha(t)-\alpha(s)\le\frac{t-s}{\rho(B)}\;,\quad(0\le s<t\le1)\;,
\]
Lemma~\ref{lem:Xpieces} gives
$\|qT_\ell q-\Tt\|\le2^{1-(\ell/2)}$.
Letting $\ell\to\infty$ yields $\Tt=qTq$.
\end{proof}

We will use the following freeness result concerning R--diagonal elements and Haar unitaries.
\begin{prop}
\label{prop:freeelts}
Let $(\Bfr,\phi)$ be a tracial W$^*$--noncommutative probability space and let $1\in A\subseteq\Bfr$ be a unital subalgebra.
Let $N\in\Nats$, let $b_{ij}$ ($1\le i<j\le N$) be R--diagonal elements and let $u_i$ ($1\le i\le N$) be Haar
unitaries such that
\[
A,\quad(\{b_{ij}\})_{1\le i<j\le N}\;,\quad(\{u_i\})_{1\le i\le N}
\]
is a $*$--free family.
Then the family
\[
(u_i^*Au_i)_{1\le i\le N}\;,\quad(\{u_i^*b_{ij}u_j\})_{1\le i<j\le N}
\]
is $*$--free.
\end{prop}
\begin{proof}
We may without loss of generality assume $b_{ij}=v_{ij}|b_{ij}|$, where each $v_{ij}$ is a Haar unitary and the pair
$v_{ij}$, $|b_{ij}|$ is $*$--free.
Then
\[
u_i^*b_{ij}u_j=(u_i^*v_{ij}u_j)(u_j^*|b_{ij}|u_j)\;.
\]
Since for a given $j$, the family
\[
u_j^*Au_j\;,\quad(\{u_j^*|b_{ij}|u_j\})_{1\le i<j}
\]
is $*$--free, letting $\Afr=W^*(A\cup\{|b_{ij}|\mid1\le i<j\le N\})$,
it will suffice to show the family
\[
(u_i^*\Afr u_i)_{1\le i\le N}\;,\quad(\{u_i^*v_{ij}u_j\})_{1\le i<j\le N}
\]
is $*$--free.
We have that the family
\[
\Afr,\quad(\{u_i\})_{1\le i\le N}\;,\quad(\{v_{ij})_{1\le i<j\le N}
\]
is $*$--free.

We will use the notation for any subalgebra $C$ of $\Bfr$,
\[
C\oup=C\cap\ker\phi.
\]
Using freeness, we have
\[
(u_i^*\Afr u_i)\oup=u_i^*\Afr\oup u_i.
\]
Given subsets $X_1,\ldots,X_n$ of $\Bfr$, we will use the notation
\begin{equation}
\label{eq:Lambdao}
\Lambdao(X_1,\ldots,X_n)
\end{equation}
for the set of all words $x_1x_2\ldots x_k$ with $k\ge 1$, $x_j\in X_{i(j)}$ and $i(1)\ne i(2),\,i(2)\ne i(3),\ldots,i(k-1)\ne i(k)$.
When refering to elements of~\eqref{eq:Lambdao},
we will frequently and intentionally shift between formally defined words and the corresponding
elements of $\Bfr$ obtained by performing the product operation.

Consider first $\Theta=\Lambdao((u_i^*\Afr\oup u_i)_{1\le i\le N})$.
Clearly $\Theta\subseteq\Psi$, where $\Psi$ is the set of all words
belonging to $\Lambdao(\Afr\oup,(\{u_i,u_i^*\})_{1\le i\le N})$ of length three or greater
\renewcommand{\labelenumi}{$\bullet$}
\begin{enumerate}
\item whose first letter is $u_i^*$, some $1\le i\le N$,
\item whose second letter is from $\Afr\oup$,
\item whose penultimate letter is from $\Afr\oup$,
\item whose last letter is $u_j$, some $1\le j\le N$.
\end{enumerate}
Thus it will suffice to show
\begin{equation}
\label{eq:LamPsi}
\Lambdao\big(\Psi,\big(\{(u_i^*v_{ij}u_j)^n\mid n\in\Ints\backslash\{0\}\}\big)_{1\le i<j\le N}\big)
\subseteq\ker\phi\;.
\end{equation}
However, taking a word $x$ belonging to the LHS of~\eqref{eq:LamPsi}, we see that some cancellation of
neighboring letters of the form $u_ju_j^*$ may possibly be performed, however only in the following situations:
\begin{alignat*}{3}
&\bullet\qquad&\cdots\Afr\oup u_j)&(u_j^*v_{jk}\cdots&\qquad&(j<k) \\
&\bullet&\cdots\Afr\oup u_j)&(u_j^*v_{ij}^*\cdots&&(i<j) \\
&\bullet&\cdots v_{ij}u_j)&(u_j^*\Afr\oup\cdots&&(i<j) \\
&\bullet&\cdots v_{jk}^*u_j)&(u_j^*\Afr\oup\cdots&&(j<k) \\
&\bullet&\cdots v_{ij}u_j)&(u_j^*v_{jk}\cdots&&(i<j<k) \\
&\bullet&\cdots v_{ij}u_j)&(u_j^*v_{pj}^*\cdots&&(i,p<j,\,i\ne p) \\
&\bullet&\cdots v_{jk}^*u_j)&(u_j^*v_{jq}\cdots&&(j<k,q,\,k\ne q) \\
&\bullet&\cdots v_{jk}^*u_j)&(u_j^*v_{pj}^*\cdots&&(p<j<k).
\end{alignat*}
Clearly, after making all such possible cancellations of $u_j^*u_j$ for all $j$, no further cancellations
are possible, and we are left with an element of
\[
\Lambdao\big(\Afr\oup,\;(\{u_i,u_i^*\})_{1\le i\le N}\;,\;(\{v_{ij},v_{ij}^*\})_{1\le i<j\le N}\big)\;.
\]
By freeness, this implies $\phi(x)=0$.
\end{proof}

\begin{thm}
\label{thm:DTfreematrix}
Let $N\in\Nats$ and let $\mu_1,\ldots,\mu_N$ be compactly supported Borel probability measures on $\Cpx$.
Let $c>0$ and suppose that in a W$^*$--noncommutative probability space $(\Mcal,\tau)$,
\[
(a_k)_{k=1,}^N\quad(b_{ij})_{1\le i<j\le N}
\]
is a $*$--free family, where $a_k$ is $\DT(\mu_k,\frac c{\sqrt N})$
and each $b_{ij}$ is circular with $\tau(|b_{ij}|^2)=\frac{c^2}N$.
Let
\[
Z=\left(\begin{matrix}
   a_1 & b_{12} & \cdots & b_{1N}    \\
     0 &  a_2   & \ddots & \vdots    \\
\vdots & \ddots & \ddots & b_{N-1,N} \\
     0 & \cdots &    0   & a_N
\end{matrix}\right) \in M_N(\Mcal)\;.
\]
Then with respect to the state $\tau\circ\tr_N$ on $M_N(\Mcal)$, $Z$ is a $\DT(\mu,c)$--element,
where $\mu=\frac1N(\mu_1+\cdots+\mu_N)$.
\end{thm}
\begin{proof}
In light of Proposition~\ref{prop:constadj}, we may without loss of generality fix $c=1$.
Let $(\Afr,\tau)$ be a tracial W$^*$--noncommutative probability space having semicircular elements
$x_i$ with $\tau(x_i^2)=1$, ($1\le i\le N$), circular elements $y_{ij}$ with $\tau(|y_{ij}|^2)=1$,
($1\le i<j\le N$),
and with a normal, injective $*$--homomorphism $\lambda:L^\infty[0,1]\to\Afr$ such that the family
\[
(\{x_i\})_{1\le i\le N}\;,\quad(\{y_{ij},y_{ij}^*\})_{1\le i<j\le N}\;,\quad\lambda(L^\infty[0,1])
\]
is free.
By random matrix results of Voiculescu~\cite{V:RM}, it follows that, with respect to the obvious tracial state
on $M_N(\Afr)$, the element
\[
X=\frac1{\sqrt N}\left(\begin{matrix}
   x_1     & y_{12} &      \cdots & y_{1N}    \\
  y_{12}^* &  x_2   &      \ddots & \vdots    \\
    \vdots & \ddots &      \ddots & y_{N-1,N} \\
  y_{1N}^* & \cdots & y_{N-1,N}^* & a_N
\end{matrix}\right) \in M_N(\Afr)
\]
is semicircular with second moment $1$ and is free from the image of $\kappa$,
where $\kappa:L^\infty[0,1]\to M_N(\Afr)$ is given by
\[
\kappa(f)=\diag\big(\lambda(f\cdot1_{[0,1/N]}),\lambda(f\cdot1_{[1/N,2/N]}),\ldots,\lambda(f\cdot1_{[(N-1)/N,1]})\big)\;.
\]
Let $f_i\in L^\infty[0,1]$ be such that the push--forward of Lebesgue measure under $f_i$ is $\mu_i$,
Let $d_i=\lambda(f_i)$ and let $D=\diag(d_1,\ldots,d_N)\in M_N(\Afr)$.
Then $D=\kappa(f)$ for some $f\in L^\infty[0,1]$ where the push forward of Lebesgue measure under $f$ is $\mu$.
Let $T=\UT(X,\kappa)$.
By Theorem~\ref{thm:DTinLF2}, $D+T$ is a $\DT(\mu,1)$--element.
From Lemma~\ref{lem:Tbasics},
\[
D+T=\sum_{1\le i\le j\le N}e_{ii}(D+T)e_{jj}
=\left(\begin{matrix}
   a_1 & b_{12} & \cdots & b_{1N}    \\
     0 &  a_2   & \ddots & \vdots    \\
\vdots & \ddots & \ddots & b_{N-1,N} \\
     0 & \cdots &    0   & a_N
\end{matrix}\right)\;,
\]
where $b_{ij}=\frac1{\sqrt N}y_{ij}$ and $a_i=d_i+e_{ii}Te_{ii}$.
By Lemma~\ref{lem:qTq},
\[
e_{ii}Te_{ii}=\UT(e_{ii}Xe_{ii},\lambda)=\UT(\tfrac1{\sqrt N}x_i,\lambda)=\tfrac1{\sqrt N}\UT(x_i,\lambda)\;,
\]
so by Theorem~\ref{thm:DTinLF2} again, $a_i$ is a $\DT(\mu_i,\frac1{\sqrt N})$--element.

Clearly, the family
\[
\{a_i\mid1\le i\le N\},\quad(\{b_{ij}\})_{1\le i<j\le N}
\]
is $*$--free, however the various $a_i$ are not free among themselves.
We may without loss of generality assume there are Haar unitaries $u_1,\ldots,u_N\in\Afr$ so that 
\begin{equation}
\label{eq:uabfree}
(\{u_i\})_{1\le i\le N}\,\quad\{a_i\mid1\le i\le N\},\quad(\{b_{ij}\})_{1\le i<j\le N}
\end{equation}
is a $*$--free family.
Let $U=\diag(u_1,\ldots,u_N)\in M_N(\Afr)$.
Then $U^*(D+T)U$ is a $\DT(\mu,1)$--element and 
\[
U^*(D+T)U=\left(\begin{matrix}
 \at_1 & \bt_{12} & \cdots & \bt_{1N}    \\
     0 &  \at_2   & \ddots & \vdots      \\
\vdots &  \ddots  & \ddots & \bt_{N-1,N} \\
     0 &  \cdots  &    0   & \at_N
\end{matrix}\right)\;,
\]
where $\at_i=u_i^*a_iu_i$ and $\bt_{ij}=u_i^*b_{ij}u_j^*$.
By Proposition~\ref{prop:freeelts}, the family
\[
\{\at_i\mid1\le i\le N\},\quad(\{\bt_{ij}\})_{1\le i<j\le N}
\]
is $*$--free.
We have therefore realized a $\DT(\mu,1)$--element as a matrix having $*$--free entries of the desired form.
\end{proof}

\section{Decomposability of DT--operators}
\label{sec:decomp}

In this section, we apply the results of~\S\ref{sec:DTinLF2} to show that every DT--operator is
strongly decomposable.

\begin{lemma}
\label{lem:specAB}
If $A$ and $B$ are bounded operators on a Hilbert space and if $A$ is normal then
\[
\sigma(A+B)\subseteq\{z\in\Cpx\mid d(z,\sigma(A))\le\|B\|\}\;.
\]
\end{lemma}
\begin{proof}
If $\lambda\in\Cpx$ and $d(\lambda,\sigma(A))>\|B\|$ then $\|(A-\lambda)^{-1}\|<\|B\|^{-1}$.
Therefore $(A-\lambda)^{-1}(A+B-\lambda)=1-(A-\lambda)^{-1}B$ is invertible, and $\lambda\notin\sigma(A+B)$.
\end{proof}

\begin{thm}
\label{thm:DTspec}
Let $\mu$ be a compactly supported Borel probability measure on $\Cpx$ and let $c>0$.
If $a$ is a $\DT(\mu,c)$ operator then $\sigma(a)=\supp(\mu)$.
\end{thm}
\begin{proof}
Let $(\Mcal,\tau)$ be a W$^*$--noncommutative probability space, with $\tau$ a faithful trace.
The spectrum of an operator $a\in\Mcal$ depends only on its $*$--moment distribution with respect to $\tau$.
Let $N\in\Nats$ and let $(a_k)_{k=1}^N,\,(b_{ij})_{1\le i,j\le N}$ be a $*$--free family of elements of $\Mcal$,
where each $a_k$ is $\DT(\mu,\frac c{\sqrt N})$ and each $b_{ij}$ is circular with $\tau(|b_{ij}|^2)=\frac{c^2}N$.
Let
\[
x=\left(\begin{matrix}
   a_1 & b_{12} & \cdots & b_{1N}    \\
     0 &  a_2   & \ddots & \vdots    \\
\vdots & \ddots & \ddots & b_{N-1,N} \\
     0 & \cdots &    0   & a_N
\end{matrix}\right) \in M_N(\Mcal)\;.
\]
Then by Theorem~\ref{thm:DTfreematrix}, $x$ is $\DT(\mu,c)$ with respect to $\tau\circ\tr_N$.
Since $\Mcal$ is a finite von Neumann algebra, we have
$\sigma(x)=\bigcup_{k=1}^N\sigma(a_k)=\sigma(a_1)$.
We have shown that the spectrum of a $\DT(\mu,c)$ operator is the same as that of a $\DT(\mu,\frac c{\sqrt N})$ operator.

Let us use Theorem~\ref{thm:DTinLF2}, to realize a $\DT(\mu,c)$--operator as $D+cT$,
for $D$ a normal operator with spectrum $\supp(\mu)$ and
as noted in Remark~\ref{rem:normT}, $\|T\|\le2$.
We have $\sigma(D+cT)=\sigma(D+\frac c{\sqrt N}T)$ for all $N\in\Nats$.
Using Lemma~\ref{lem:specAB} we thus get
\[
\sigma(D+cT)\subseteq\{z\in\Cpx\mid d(z,\supp(\mu))\le\tfrac{2c}{\sqrt N}\}
\]
for all $N\in\Nats$.
Hence $\sigma(D+cT)\subseteq\supp(\mu)$.

On the other hand, suppose for contradiction that $\lambda\in\supp(\mu)\backslash\sigma(D+cT)$.
Then there is $\eps>0$ such that the ball $B_\eps(\lambda)$ of radius $\eps$ around $\lambda$
is disjoint from $\sigma(D+cT)$.
Let $N\in\Nats$ be such that $N\ge\mu(B_\eps(\lambda))^{-1}$.
Using Theorem~\ref{thm:DTfreematrix}, we may take $(\Mcal,\tau)=(M_N(\Nc),\tau_\Nc\circ\tr_N)$ and
\[
D+cT=x=\left(\begin{matrix}
   a_1 & b_{12} & \cdots & b_{1N}    \\
     0 &  a_2   & \ddots & \vdots    \\
\vdots & \ddots & \ddots & b_{N-1,N} \\
     0 & \cdots &    0   & a_N
\end{matrix}\right) \in M_N(\Nc)\;,
\]
where $\Nc$ is a von Neumann algebra with faithful normal tracial state $\tau_\Nc$ and where $a_1$
is $\DT(\mu_1,\frac c{\sqrt N})$ with $\supp(\mu_1)\subseteq B_\eps(\lambda)$.
But $\sigma(D+cT)=\bigcup_{k=1}^N\sigma(a_k)$ and $\sigma(a_1)\subseteq\supp(\mu_1)$.
Since $\sigma(a_1)$ is nonempty, this implies that $B_\eps(\lambda)$ meets the spectrum of $D+cT$,
which is a contradiction.
\end{proof}

The next result relates upper triangular decompositions to local spectral subspaces.
\begin{prop}
\label{prop:utmat}
Let $T$ be a bounded operator on a separable Hilbert space $\HEu=\HEu_1\oplus\HEu_2$
and suppose $\HEu_1\oplus0$ is $T$--invariant, so that
$T=\left(\begin{smallmatrix}T_{11}&T_{12}\\0&T_{22}\end{smallmatrix}\right)$
where $T_{ij}:\HEu_j\to\HEu_i$.
Let $A\subseteq\Cpx$.
If $\sigma(T_{11})\subseteq A$, then $\HEu_1\oplus0\subseteq\HEu_T(A)$.
On the other hand, if $\sigma(T_{22})\cap A=\emptyset$, then $\HEu_T(A)\subseteq\HEu_1\oplus0$.
\end{prop}
\begin{proof}
Write elements $\xi$ of $\HEu$ as $\xi=\left(\begin{smallmatrix}\xi_1\\ \xi_2\end{smallmatrix}\right)$ with $\xi_i\in\HEu_i$.
Suppose $\sigma(T_{11})\subseteq A$.
For $\xi_1\in\HEu_1$, taking $f(\lambda)=\left(\begin{smallmatrix}(T_{11}-\lambda)^{-1}\xi_1\\ 0\end{smallmatrix}\right)$
we have $(T-\lambda)f(\lambda)=\left(\begin{smallmatrix}\xi_1\\ 0\end{smallmatrix}\right)$, ($\lambda\in\Cpx\backslash\sigma(T_{11})$).
Thus $\sigma_T(\left(\begin{smallmatrix}\xi_1\\ 0\end{smallmatrix}\right))\subseteq\sigma(T_{11})\subseteq A$
and therefore $\left(\begin{smallmatrix}\xi_1\\ 0\end{smallmatrix}\right)\in\HEu_T(A)$.

Now suppose $\sigma(T_{22})\cap A=\emptyset$ and $\xi=\left(\begin{smallmatrix}\xi_1\\ \xi_2\end{smallmatrix}\right)\in\HEu_T(A)$.
Let $f(\lambda)=\left(\begin{smallmatrix}f_1(\lambda)\\ f_2(\lambda)\end{smallmatrix}\right)$ be analytic such that
$(T-\lambda)f(\lambda)=\xi$, ($\lambda\in\Cpx\backslash\sigma_T(\xi)$).
However, as $\sigma_T(\xi)\cap\sigma(T_{22})=\emptyset$ and as $T$ has the single--valued extension property,
(since $\HEu$ is separable), the analytic function $f_2$ extends to an entire function $\ft_2:\Cpx\to\HEu_2$ such that
$(T_{22}-\lambda)f_2(\lambda)=\xi_2$ for all $\lambda\in\Cpx$.
Since $f_2(\lambda)=(T_{22}-\lambda)^{-1}\xi_2\to0$ as $|\lambda|\to\infty$, by Liouville's
theorem $f_2$ is the zero function and $\xi_2=0$;
thus $\left(\begin{smallmatrix}\xi_1\\ \xi_2\end{smallmatrix}\right)\in\HEu_1\oplus0$.
\end{proof}

For the rest of this section, $Z$ will be a $\DT(\mu,c)$--operator in
a W$^*$--noncommutative probability space $(\Mcal,\tau)$,
with $\Mcal\subseteq\Bc(\HEu)$, $\HEu$ separable and $\tau$ faithful.

\begin{thm}
\label{thm:closures}
Let $B$ be a Borel subset of $\Cpx$
and let $p=p_Z(B)$ be the projection from $\HEu$ onto $\overline{\HEu_Z(B)}$.
Then
\renewcommand{\labelenumi}{(\roman{enumi})}
\begin{enumerate}
\item
$p\in\Mcal$, $Zp=pZp$ and $\tau(p)=\mu(B)$
\item
if $\mu(B)>0$ then $Zp$ is $\DT(\mu(B)^{-1}\mu\restrict_B,c\sqrt{\mu(B)})$
with respect to $\mu(B)^{-1}\tau\restrict_{p\Mcal p}$
\item
if $\mu(B)<1$ then $(1-p)Z$ is $\DT(\mu(B^c)^{-1}\mu\restrict_{B^c},c\sqrt{\mu(B^c)})$
with respect to $\mu(B^c)^{-1}\tau\restrict_{(1-p)\Mcal(1-p)}\;$.
\end{enumerate}
\end{thm}
\begin{proof}
$\HEu_Z(B)$ is $Z$--hyperinvariant (see~\cite[1.2.16(a)]{LN}) and therefore so is $\overline{\HEu_Z(B)}$;
this implies $p\in\Mcal$ and $Zp=pZp$.
Let $n\in\Nats$ and let $k\in\{0,1,\ldots,n-1\}$, let $K_n\subseteq\Cpx$ be compact
and let $U_n\subseteq\Cpx$ be open such that $K_n\subseteq B\subseteq U_n$ and
\begin{equation*}
\frac kn\le\mu(K_n)\le\mu(B)\le\mu(U_n)\le\frac{k+1}n\;.
\end{equation*}
Consider the probability measures
\begin{align*}
\nu_n&=\mu(K_n)^{-1}\mu\restrict_{K_n}\text{ if }\mu(K_n)>0\;, \\
\nu'_n&=\mu(U_n^c)^{-1}\mu\restrict_{U_n^c}\text{ if }\mu(U_n)<1\;, \\
\nu''_n&=n\mu\restrict_{U_n\backslash K_n}+(n\mu(K_n)-k)\nu_n+(k+1-n\mu(U_n))\nu_n'\;.
\end{align*}
Then $\mu=\frac1n(k\nu_n+\nu''_n+(n-k-1)\nu'_n)$ and hence by Theorem~\ref{thm:DTfreematrix}
there is a $\DT(\mu,c)$ operator
\[
Z_n=\left(\begin{matrix}
   a_1 & b_{12} & \cdots & b_{1n}    \\
     0 &  a_2   & \ddots & \vdots    \\
\vdots & \ddots & \ddots & b_{n-1,n} \\
     0 & \cdots &    0   & a_n
\end{matrix}\right)
\]
in $(M_n(\Nc),\tau_\Nc\circ\tr_n)$ where in $(\Nc,\tau_\Nc)$,
$a_i$ is $\DT(\nu_n,\frac c{\sqrt n})$ if $i\le k$
and is $\DT(\nu'_n,\frac c{\sqrt n})$ if $i\ge k+2$,
$a_{k+1}$ is $\DT(\nu''_n,\frac c{\sqrt n})$,
$b_{ij}$ is circular with $\tau_\Nc(|b_{ij}|^2)=\frac{c^2}n$ and the family
$(a_i)_{i=1}^n,\,(b_{ij})_{1\le i<j\le n}$ is $*$--free.

Let
\begin{align*}
q_n&=\diag(\underset k{\underbrace{1,\ldots,1}},\underset{n-k}{\underbrace{0,\ldots,0}}) \\
q'_n&=\diag(\underset{k+1}{\underbrace{1,\ldots,1}},\underset{n-k-1}{\underbrace{0,\ldots,0}})
\end{align*}
in $M_n(\Nc)$.
Applying Theorem~\ref{thm:DTfreematrix} and Theorem~\ref{thm:DTspec},
if $k>0$ then $Z_nq_n$ is $\DT(\nu_n,c\sqrt{\frac kn})$
in $(M_k(\Nc),\tau_\Nc\circ\tr_k)$ and $\sigma(Z_nq_n)\subseteq K_n$. 
Similarly, if $k<n-1$ then $(1-q'_n)Z_n$ is $\DT(\nu'_n,c\sqrt{\frac{n-k-1}n})$
in $(M_{n-k-1}(\Nc),\tau_\Nc\circ\tr_{n-k-1})$  and $\sigma((1-q'_n)Z_n)\subseteq U_n^c$.
But Proposition~\ref{prop:utmat} then gives $q_n\le p_{Z_n}(B)\le q'_n$.
Letting $n\to\infty$, this shows that $\tau(p)=\mu(B)$.
Thus~(i) is proved.

If $\mu(B)>0$ then $Z_nq_n$ converges in $*$--moments to $Zp$ as $n\to\infty$.
Indeed, for $\eps(1),\ldots,\eps(\ell)\in\{*,1\}$ we have
\[
\tau(pZ^{\eps(1)}\cdots pZ^{\eps(\ell)}p)
=\tau_\Nc\circ\tr_N(p_{Z_n}(B)Z_n^{\eps(1)}\cdots p_{Z_n}(B)Z_n^{\eps(\ell)}p_{Z_n}(B))
\]
and this quantity differs from $\tau_\Nc\circ\tr_N(q_nZ_n^{\eps(1)}\cdots q_nZ_n^{\eps(\ell)}q_n)$
by an amount no greater than $\|q_n-p_{Z_n}(B)\|_2Q_\ell(\|Z_n\|)$, where $Q_\ell$ is a polynomial
independent of $n$ and with positive coefficients.
Since $\|q_n-p_{Z_n}(B)\|_2\le n^{-1/2}\to0$ as $n\to\infty$, this shows that $Z_nq_n$ converges to $Zp$.
Clearly, $\nu_n$ converges in $*$--moments as $n\to\infty$ to $\mu(B)^{-1}\mu\restrict_B$.
Now an application of Corollary~\ref{cor:DTconv} shows that
$Zp$ is $\DT(\mu(B)^{-1}\mu\restrict_B,c\sqrt{\mu(B)})$, and~(ii) is proved.

For~(iii), one shows similarly that if $\mu(B)<1$ then $(1-q'_n)Z_n$ converges in $*$--moments
to $(1-p)Z$ as $n\to\infty$, $\nu_n'$ converges in $*$--moments to $\mu(B^c)^{-1}\mu\restrict_{B^c}$
and thus that $(1-p)Z$ is $\DT(\mu(B^c)^{-1}\mu\restrict_{B^c},c\sqrt{\mu(B^c)})$.
\end{proof}

L.G.\ Brown~\cite{B} discovered a measure supported inside the spectrum of an arbitrary element
of a finite von Neumann algebra, which in the finite dimensional case reduces to the eigenvalue distribution
weighted according to generalized multiplicity of eigenvalues.

\begin{cor}
\label{cor:Brownmeas}
The Brown measure of a $\DT(\mu,c)$ operator is $\mu$.
\end{cor}
\begin{proof}
Let $Z$ be a $\DT(\mu,c)$--operator, and let $\nu_Z$ be the Brown measure of $Z$.
From the theorem just proved and~\cite[Theorem 4.3]{B}, it follows that $\nu_Z(F)\ge\tau(p_Z(F))=\mu(F)$
for every closed subset $F$ of $\Cpx$.
Thus $\nu_Z=\mu$.
\end{proof}

The next result is a converse to Theorem~\ref{thm:closures}, which we will use in~\S\ref{sec:vN}.
\begin{thm}
\label{thm:closuresconverse}
\renewcommand{\labelenumi}{(\roman{enumi})}
\begin{enumerate}
\item Let $F\subseteq\Cpx$ be a closed subset and
suppose $q\in\Mcal$ is a projection satisfying $Zq=qZq$, $\tau(q)=\mu(F)$
and $\sigma(Z\restrict_{q\HEu})\subseteq F$.
Then $q=p_Z(F)$.
\item Let $B\subseteq\Cpx$ be a Borel subset and suppose $q\in\Mcal$
is a projection satisfying $Zq=qZq$, $\tau(q)=\mu(B)$ and that
$Zq$ is a $\DT(\mu(B)^{-1}\mu\restrict_B,c')$--operator
with respect to $\mu(B)^{-1}\tau\restrict_{q\Mcal q}$ for some $c'>0$.
Then $q=p_Z(B)$ and $c'=c\sqrt{\mu(B)}$.
\end{enumerate}
\end{thm}
\begin{proof}
In~(i), Proposition~\ref{prop:utmat} implies $q\le p_Z(F)$.
But $\tau(q)=\tau(p_Z(F))$, so $q=p_Z(F)$.

By the definition of spectral subspaces we have
\begin{equation}
\label{eq:HZB}
\HEu_Z(B)=\bigcup_F\HEu_Z(F),
\end{equation}
where the union is over all closed subsets $F$ of $\Cpx$ such that $F\subseteq B$.
Regard $Zq$ as an operator on $q\HEu$ and let $F\subseteq B$ be a closed subset of $\Cpx$.
Then by Theorem~\ref{thm:closures}, the $Zq$--invariant projection $p_{Zq}(F)$
satisfies $\tau(p_{Zq}(F))=\mu(F)$ and $\sigma(Z\restrict_{p_{Zq}(F)\HEu})=F$.
But $p_{Zq}(F)$ is also $Z$--invariant, and hence part~(i) implies
$p_Z(F)=p_{Zq}(F)\le q$.
Using~\eqref{eq:HZB} we get $p_Z(B)\le q$.
Now $\tau(p_Z(B))=\tau(q)$ implies $p_Z(B)=q$.
\end{proof}

The following result shows that $Z$ has property~(C) of Dunford; (see~\cite[1.2.18]{LN}).

\begin{lemma}
\label{lem:pC}
If $F$ is a closed subset of $\Cpx$ then $\HEu_Z(F)$ is closed.
\end{lemma}
\begin{proof}
If $\mu(F)=0$ then by Theorem~\ref{thm:closures}, $\overline{\HEu_Z(F)}=\HEu_Z(F)=\{0\}$,
so suppose $\mu(F)>0$.
We must show $\overline{\HEu_Z(F)}=\HEu_Z(F)$ or, in other words, that if $\xi\in\overline{\HEu_Z(F)}$
then there is an analytic function $f:\Cpx\backslash F\to\HEu$ such that $(Z-\lambda)f(\lambda)=\xi$,
($\lambda\in\Cpx\backslash F$).
Letting $p=p_Z(F)$ be the projection from $\HEu$ onto $\overline{\HEu_Z(F)}$, then $Zp=pZp$ and
by Theorem~\ref{thm:closures} and Theorem~\ref{thm:DTspec}, $\sigma(Zp)\subseteq F$,
where the spectrum is for $Zp$ as an element of $p\Mcal p$.
Now given $\xi\in p\HEu$ let $f(\lambda)=(pZp-\lambda p)^{-1}\xi$ ($\lambda\in\Cpx\backslash F$),
where the inverse is taken in $p\Mcal p$.
Then $f$ is analytic and $(Z-\lambda)f(\lambda)=\xi$.
\end{proof}

\begin{thm}
\label{thm:decomp}
$Z$ is strongly decomposable.
\end{thm}
\begin{proof}
We will first show that $Z$ is decomposable, and we will use the characterization
of decomposability found at~\cite[1.2.23(b)]{LN}, which is essentially due to
Jafarian and Vasilescu~\cite{JV}.
Since we already showed in Lemma~\ref{lem:pC} that $Z$ has property~(C), it remains only to show that
for every closed subset $F$ of $\Cpx$, the spectrum of the operator $Z/\HEu_Z(F)$ induced by $Z$
on the quotient space $\HEu/\HEu_Z(F)$ is contained in $\overline{\sigma(Z)\backslash F}$.
Note that $\sigma(Z)=\supp(\mu)$.
Consider first the trivial cases $\mu(F)\in\{0,1\}$.
From Theorem~\ref{thm:closures} we have $\HEu_Z(F)=\{0\}$ if and only if $\mu(F)=0$,
and in this case $\intr(F)=\emptyset$ so $\overline{\sigma(Z)\backslash F}=\sigma(Z)$,
while $\sigma(Z/\HEu_Z(F))=\sigma(Z)$.
If $\mu(F)=1$ then $\HEu_Z(F)=\HEu$ and $Z/\HEu_Z(F)$ is the operator on the Hilbert space $\{0\}$,
which has empty spectrum.

Suppose $0<\mu(F)<1$ and let $p=p_Z(F)$ be the projection from $\HEu$ onto $\HEu_Z(F)$.
Then $\HEu/\HEu_Z(F)$ is canonically isomorphic to $(1-p)\HEu$ in such a way that
$Z/\HEu_Z(F)$ corresponds to $(1-p)Z\in\Bc((1-p)\HEu)$.
By Theorems~\ref{thm:closures} and~\ref{thm:DTspec}, the spectrum of $(1-p)Z$
is $\supp(\mu\restrict_{F^c})=\overline{\supp(\mu)\backslash F}=\overline{\sigma(Z)\backslash F}$.
Thus $Z$ is decomposable.

In order to show that $Z$ is strongly decomposable we must show that for every closed subset $F$
of $\Cpx$, the restriction $Z\restrict_{\HEu_Z(F)}$ of $Z$ to $\HEu_Z(F)$ is decomposable.
However, by Theorem~\ref{thm:closures}, either $\HEu_Z(F)=\{0\}$ or $Z\restrict_{\HEu_Z(F)}$
is itself a DT--operator.
In either case, $Z\restrict_{\HEu_Z(F)}$ is decomposable.
\end{proof}

\section{Von Neumann algebras generated by DT--operators}
\label{sec:vN}

Theorem~\ref{thm:DTinLF2} allows us realize a $\DT(\mu,c)$--operator as an element
$D+cT$ in the II$_1$--factor $L(\Fb_2)$, where we have some flexibility in choosing $D$.
It is natural to ask:  what sort of II$_1$--factors do these DT--operators generate?
It is this question that is addressed in this section.
We first consider the von Neumann algebra, which we will denote $\Qc$,
generated by $T=\UT(X,\lambda)$ inside $L(\Fb_2)$ (cf Definition~\ref{def:UT}).
We show (Theorem~\ref{thm:Tfactor}) that $\Qc$ is an irreducible subfactor of $L(\Fb_2)$, but
we are presently unable to decide whether $\Qc$ is all of $L(\Fb_2)$.
This question is closely bound up with the question of whether $T$ has nontrivial
hyperinvariant subpaces.

In a tracial W$^*$--noncommutative probability space $(\Mcal,\tau)$, let $X$ be a semicircular operator
and let $\lambda:L^\infty[0,1]\to\Mcal$ be a normal, unital, injective $*$--homomorphism whose image is free from $X$.
Thus $X$ and the image of $\lambda$ together generate a copy of $L(\Fb_2)$.
\begin{thm}
\label{thm:Tfactor}
Let $T=\UT(X,\lambda)$ be the operator as constructed in Lemma~\ref{lem:Tconstr} inside
the von Neumann algebra factor $\Nc=(\{X\}\cup\lambda(L^\infty[0,1]))''\cong L(\Fb_2)$.
Then the von Neumann algebra generated by $T$ is an irreducible subfactor of $\Nc$, and is in particular a factor.
\end{thm}
\begin{proof}
We use the convention, given a state $\phi$ on a von Neumann algebra $\Mcal$, that $a\mapsto\ah$ denotes
the defining mapping of $\Mcal$ into $L^2(\Mcal,\phi)$.

The W$^*$--subalgebra of $\Nc$ generated by $X$ is canonically isomorphic to $L^\infty[-2,2]$.
Since $W^*(X)$ and the image of $\lambda$ together freely generate $\Nc$, the free product construction~\cite{V:sym}
identifies $L^2(\Nc,\tau\restrict_\Nc)$ with the Hilbert space
\begin{equation*}
\HEu=\Cpx\oplus\bigoplus_{
\substack{n\ge1\\\iota_1,\ldots,\iota_n\in\{1,2\}\\\iota_j\ne\iota_{j+1}} }
\HEu_{\iota_1}\oup\otimes\cdots\otimes\HEu_{\iota_n}\oup\;,
\end{equation*}
where
\begin{align*}
\HEu_1\oup&=L^2[0,1]\ominus\Cpx\oneh=(\lambda(L^\infty[0,1])\cap\ker\tau)\widehat{\;} \\
\HEu_2\oup&=L^2[-2,2]\ominus\Cpx\oneh=(W^*(X)\cap\ker\tau)\widehat{\;}\;.
\end{align*}
For example, in this identification, an element $(a_1b_1a_2b_2\cdots a_nb_n)\hat{\;}$
of $L^2(\Nc,\tau\restrict_\Nc)$, where $a_i\in W^*(X)\cap\ker\tau$ and $b_i\in\lambda(L^\infty[0,1])\cap\ker\tau$,
corresponds to the element
$\ah_1\otimes\bh_1\otimes\cdots\otimes\ah_n\otimes\bh_n$ of the direct summand 
$\HEu_2\oup\otimes\HEu_1\oup\otimes\cdots\otimes\HEu_2\oup\otimes\HEu_1\oup$ of $\HEu$.
The map $L^\infty[0,1]\odot L^\infty[0,1]\to\Nc$, where $\odot$ denotes the algebraic tensor product,
defined by $f\otimes g\mapsto\lambda(f)X\lambda(g)$,
gives rise to an isometry $W:L^2[0,1]\otimes L^2[0,1]\to\HEu$ at the $L^2$--level.
Specifically,
\begin{align*}
W(v_1\otimes v_2)=&\langle v_1,\oneh\rangle\Xh\langle v_2,\oneh\rangle
+(v_1-\langle v_1,\oneh\rangle\oneh)\otimes\Xh\langle v_2,\oneh\rangle \\
&+\langle v_1,\oneh\rangle\Xh\otimes(v_2-\langle v_2,\oneh\rangle\oneh)
+(v_1-\langle v_1,\oneh\rangle\oneh)\otimes\Xh\otimes(v_2-\langle v_2,\oneh\rangle\oneh).
\end{align*}
In other words, $W(v_1\otimes v_2)=v_1\otimes\Xh\otimes v_2$ if we make the convention that
$\oneh$ is absorbed by tensor products.
Then from the definition of $T$, we have $\Th=W(h)$, where $h\in L^2[0,1]\otimes L^2[0,1]$ is given,
upon identifying this Hilbert space with $L^2([0,1]^2)$, by
\begin{equation*}
h(s,t)=\begin{cases}
1&\text{if }s<t \\
0&\text{if }s\ge t\;.
\end{cases}
\end{equation*}

Suppose $a\in\Nc\cap\{T,T^*\}'$ and $\tau(a)=0$.
Since $aX=Xa$, it is well--known that we must have $a\in W^*(X)$, so $\ah\in\HEu_2\oup$.
Let $\oneh=b_0,b_1,b_2,\ldots$ be an orthonormal basis for $L^2[0,1]$.
Then
\begin{equation*}
h=\sum_{i,j=0}^\infty c_{ij}b_i\otimes b_j
\end{equation*}
for some $c_{ij}\in\Cpx$.
Clearly $c_{ij}\ne0$ for some $i,j>0$.
Thus
\begin{align*}
\Th&=c_{00}\Xh+\sum_{i>0}c_{i0}b_i\otimes\Xh+\sum_{j>0}c_{0j}\Xh\otimes b_j
+\sum_{i,j>0}c_{ij}b_i\otimes\Xh\otimes b_j \\
(aT-Ta)\widehat{\;}&=v+\sum_{i,j>0}c_{ij}
(\ah\otimes b_i\otimes\Xh\otimes b_j-b_i\otimes\Xh\otimes b_j\otimes\ah)\;,
\end{align*}
where
\begin{equation*}
v\in\Cpx\oplus\HEu_2\oup\oplus(\HEu_1\oup\otimes\HEu_2\oup)
\oplus(\HEu_2\oup\otimes\HEu_1\oup)\oplus(\HEu_2\oup\otimes\HEu_1\oup\otimes\HEu_2\oup)\;.
\end{equation*}
Since $a$ commutes with $T$, we conclude that $\ah=0$, i.e.\ $a=0$.
\end{proof}

As remarked above, we will let $\Qc$ denote the von Neumann subalgebra of $\Nc$ generated by $T=\UT(X,\lambda)$.
It is an interesting question whether the von Neumann algebra $\Qc$ generated by $T$ contains any nontrivial
invariant projections for $T$, i.e. any projections $p$ neither $0$ nor $1$ such that $Tp=pTp$.
Clearly, if $\Qc$ is all of $\Nc$, then many such projections exist.
Moreover, if there is one such projection, then by compressing one can find a coninuum of them, that together densely
span a diffuse abelian subalgebra of $\Qc$.
Finally, using this idea we see that if $\lambda(1_{[0,t]})\in\Qc$ for some $0<t<1$ then $\Qc=\Nc$.

\begin{lemma}
\label{lem:DinM}
Let $\mu$ be a compactly supported Borel measure on $\Cpx$.
Then there is $f\in L^\infty[0,1]$ whose distribution is $\mu$ and such that,
if $D=\lambda(f)$ and if $T=\UT(X,\lambda)$ is as above,
then for any $c>0$, $D$ itself lies in the W$^*$--algebra generated by $D+cT$.
\end{lemma}
\begin{proof}
Simply take $f\in L^\infty([0,1])$ whose distribution is $\mu$ and that satisfies the following:
\renewcommand{\labelenumi}{(\roman{enumi})}
\begin{enumerate}
\item for every atom $a$ of $\mu$, $f^{-1}(\{a\})$
      is a half--open interval $[c,d)$;
\item if $a_i$ ($i\in I$) are the atoms of $\mu$ then 
      the restriction of $f$ to 
      \[
      [0,1]\backslash f^{-1}(\{a_i\mid i\in I\})
      \]
      is a measure preserving isomorphism of measure spaces.
\end{enumerate}
Let $Z=D+cT$ and let $W^*(Z)$ denote the von Neumann subalgebra of $\Nc$ generated by $Z$.
If $t\in[0,1]$ is such that $f^{-1}(f([0,t)))=[0,t)$
then letting $B=f([0,t))$, and using Lemma~\ref{lem:qTq} and Theorem~\ref{thm:closuresconverse}, we get
$p_Z(B)=\lambda(1_{[0,t)})\in W^*(Z)$.
Since $p_Z(B)\in W^*(Z)$, it follows that $D\in W^*(Z)$.
\end{proof}

Using Lemma~\ref{lem:DinM} we see that if the von Neumann algebra $\Qc$ is in fact all of $\Nc$, then every DT--operator
can be embedded in $\Nc$ in such a way that it generates all of $\Nc$, which is isomorphic to $L(\Fb_2)$.

\begin{prop}
\label{prop:atomdep}
If $Z$ is a $\DT(\mu,c)$--operator, then the von Neumann algebra generated by $Z$
is a II$_1$--factor whose isomorphism class depends only on $\mu$, and in fact, only on the number and sizes of
the atoms of $\mu$.
\end{prop}
\begin{proof}
Using Lemma~\ref{lem:DinM} and
Theorem~\ref{thm:DTinLF2}, we see that the isomorphism class of the von Neumann algebra generated by a $\DT(\mu,c)$--operator
depends only on the isomorphism class of the von Neumann algebra generated by a normal element whose distribution is $\mu$.
\end{proof}

\begin{thm}
\label{thm:LF2}
If $\mu$ has no atoms and $c>0$ then a $\DT(\mu,c)$--operator generates a von Neumann algebra isomorphic to $L(\Fb_2)$.
\end{thm}
\begin{proof}
In the proof of Lemma~\ref{lem:DinM}, we see that
$D$ generates all of $\lambda(L^\infty([0,1]))$ and the constructed DT-operator generates all of $\Nc\cong L(\Fb_2)$.
\end{proof}

\section{DT--operators that are also R--diagonal}
\label{sec:DTRdiag}

In~\cite{DH} we introduced the circular free Poisson elements of parameters
$c\ge1$,
which are R--diagonal elements and 
(see Theorem~\ref{thm:cfPDT}) also DT--elements.
In this section we will show that, up to multiplication with scalars,
these are the only DT--elements that are also R--diagonal.

Let $\mu$ be a Borel probability measure compactly supported in $\Cpx$.
For every $n\in\Nats$, let $D(n)\in\MEu_n$ be a diagonal random matrix
whose entries are independent, each with distribution $\mu$,
and let $T(n)\in\UTGRM(n,\frac1n)$ be such that $D(n)$ and $T(n)$
are independent as matrix--valued random variable.
From Theorem~\ref{thm:DTdef}, (see also Remarks~\ref{rem:DTbdd} and~\ref{rem:normT}) it follows
that the pair $(D(n),T(n))$ converges in $*$--moments
as $n\to\infty$ to a pair $(D,T)$ in a tracial
W$^*$--noncommutative probability space $(\Mcal,\tau)$.
Of course, by Definition~\ref{def:DT}, $D+cT$ is a $\DT(\mu,c)$--element.

In the lemmas to follow, orthogonality and the $2$--norm, $\|\cdot\|_2$,
will be with respect to the trace $\tau$.

\begin{lemma}
\label{lem:DTorthog}
\renewcommand{\labelenumi}{(\roman{enumi})}
\begin{enumerate}
\item For all $p\in\Nats\cup\{0\}$,
      \[
      \|(DT)^pD\|_2^2=\frac{\|D\|_2^{2(p+1)}}{(p+1)!}
      \]
\item If $p,q\in\Nats\cup\{0\}$ and $p\ne q$, then $(DT)^pD$ and $(DT)^qD$
      are orthogonal.
\end{enumerate}
\end{lemma}
\begin{proof}
By the proof of Theorem~\ref{thm:DTdef},
$\tau((DT)^pDD^*(T^*D^*)^q)$ is equal to  the quantity~\eqref{eq:ENTOsum},
where
\[
\eps(1),\ldots,\eps(p+q)\quad\text{is}\quad
\underset p{\underbrace{1,\ldots,1}},
\underset q{\underbrace{*,\ldots,*}}
\]
and where $\Ec(\sigma;\cdots)$ must be properly interpreted.

If $p\ne q$ then there are no noncrossing pairings of $\{1,\ldots,p+q\}$ compatible with
$\eps(1),\ldots,\eps(p+q)$, so $\tau((DT)^pDD^*(T^*D^*)^q)=0$ and~(ii) is proved.

If $p=q$, then there exactly only one noncrossing pairing $\sigma$ of $\{1,\ldots,2p\}$
compatible with $\eps(1),\ldots,\eps(2p)$, namely $\sigma=\{\{1,2p\},\{2,2p-1\},\ldots,\{p,p+1\}\}$.
Following Algorithm~\ref{alg:NTO}, we find that the quotient graph $Q$ is the straight line having
$2p$ edges all oriented in the same direction and hence $\NTO(\sigma;\eps(1),\ldots,\eps(2p))=1$.
Furthermore, the quantity $\Ec(\sigma;\cdots)$, which is given by equation~\eqref{eq:Edlim}, is equal to
\[
\bigg(\int_\Cpx|z|^2\dif\mu(z)\bigg)^{p+1}\;.
\]
\end{proof}

\begin{lemma}
\label{lem:DTexp}
Let $M_\mu(k,\ell)=\int_\Cpx z^k\zbar^\ell\dif\mu(z)$.
Then for every $|\lambda|<(\|D\|+\|T\|)^{-1}$, we have
\begin{equation*}
\bigg\|\sum_{n=0}^\infty\lambda^n(D+T)^n\bigg\|_2^2
=\frac1{|\lambda|^2}\bigg(\exp\bigg(
\sum_{k,\ell=0}^\infty\lambda^{k+1}(\lambdab)^{\ell+1}M_\mu(k,\ell)\bigg)-1\bigg)
\end{equation*}
\end{lemma}
\begin{proof}
Because $1-\lambda D$ is invertible and
\[
\|(1-\lambda D)^{-1}\lambda T\|\le\frac{|\lambda|\|T\|}{1-|\lambda|\|D\|}<1,
\]
we have
\begin{align*}
\sum_{n=0}^\infty\lambda^n(D+T)^n&=(1-\lambda(D+T))^{-1}=((1-\lambda D)-\lambda T)^{-1} \\
&=(1-(1-\lambda D)^{-1}\lambda T)^{-1}(1-\lambda D)^{-1} \\
&=\sum_{k=0}^\infty((1-\lambda D)^{-1}\lambda T)^k(1-\lambda D)^{-1}\;.
\end{align*}
Now applying Lemma~\ref{lem:DTorthog} to $(1-\lambda D)^{-1}$ instead of $D$ gives
\begin{align*}
\bigg\|\sum_{n=0}^\infty\lambda^n(D+T)^n\bigg\|_2^2
&=\sum_{k=0}^\infty\frac{|\lambda|^{2k}}{(k+1)!}\|(1-\lambda D)^{-1}\|_2^{2(k+1)} \\
&=\frac{\exp(|\lambda|^2\|(1-\lambda D)^{-1}\|_2^2)-1}{|\lambda|^2}\;.
\end{align*}
This finishes the proof, since
\begin{align*}
|\lambda|^2\|(1-\lambda D)^{-1}\|_2^2
&=|\lambda|^2\tau((1-\lambda D)^{-1}(1-\lambdab D^*)^{-1}) \\
&=|\lambda|^2\sum_{k,\ell=0}^\infty\lambda^k(\lambdab)^\ell M_\mu(k,\ell)\;.
\end{align*}
\end{proof}

\begin{thm}
\label{thm:DTRdiag}
Let $Z$ be a $\DT(\mu,1)$--element.
Then $Z$ is R--diagonal if and only if $\mu$ is the uniform distribution on
the annulus
\begin{equation}
\label{eq:ann}
\{\lambda\in\Cpx\mid\sqrt{c-1}\le|\lambda|\le\sqrt c\}
\end{equation}
for some $c\ge1$, in which case $Z$ is a circular free Poisson element
of parameter $c$.
\end{thm}
\begin{proof}
If $\mu$ is the uniform distribution on the annulus~\eqref{eq:ann},
then by Theorem~\ref{thm:cfPDT} $\mu$ is circular free Poisson and is therefore
R--diagonal.

Suppose $Z$ is an R--diagonal element in a
$W^*$--noncommutative probability space $(\Mcal,\tau)$.
Again, $2$--norms $\|\cdot\|_2$ will be with respect to $\tau$.
We may take $Z=D+T$ with $D,T\in\Mcal$ as defined at the beginning of this section.
By~\cite{HL:br}, R--diagonality of $Z$ implies
$\|Z^n\|_2^2=\|Z\|_2^{2n}$ for every $n\in\Nats$.
Since $Z$ and $\zeta Z$ have the same $*$--moments for all $\zeta\in\Tcirc$,
we also have $\tau(Z^n(Z^*)^m)=0$ whenever $n\ne m$.
Hence for
\[
|\lambda|<(\|D\|+\|T\|)^{-1}\le\|Z\|^{-1}
\]
we have
\[
\bigg\|\sum_{n=0}^\infty\lambda^nZ^n\bigg\|_2^2=\frac1{1-|\lambda|^2\|Z\|_2^2}\;,
\]
so by Lemma~\ref{lem:DTexp},
\[
\frac1{|\lambda|^2}\bigg(\exp\bigg(\sum_{k,\ell=0}^\infty\lambda^{k+1}(\lambdab)^{\ell+1}M_\mu(k,\ell)\bigg)-1\bigg)
=\frac1{1-|\lambda|^2\|Z\|_2^2}\;,
\]
from which we get
\[
\sum_{k,\ell=0}^\infty\lambda^{k+1}(\lambdab)^{\ell+1}M_\mu(k,\ell)
=\log\left(\frac{1-|\lambda|^2(\|Z\|_2^2-1)}{1-|\lambda|^2\|Z\|_2^2}\right)\;.
\]
Since the RHS has a power series expansion in $|\lambda|^2$, we must have
$M_\mu(k,\ell)=0$ whenver $k\ne\ell$.
Thus $\mu$ is rotationally invariant.
Let $c=\|Z\|_2^2$.
Substituting $t=|\lambda|^2$ gives
\begin{align*}
\sum_{k=0}^\infty t^{k+1}M_\mu(k,k)
&=\log(1-(c-1)t)-\log(1-ct) \\
&=\sum_{j=1}^\infty\bigg(\frac{(ct)^j}j-\frac{((c-1)t)^j}j\bigg) \\
&=\sum_{k=0}^\infty\frac{c^{k+1}-(c-1)^{k+1}}{k+1}t^{k+1}\;,
\end{align*}
which implies that
\[
\tau((DD^*)^k)=M_\mu(k,k)=\frac{c^{k+1}-(c-1)^{k+1}}{k+1}t^{k+1}=\int_{c-1}^cx^k\dif x\;.
\]
Thus the distribution of $DD^*$ is the uniform distribution on the interval
$[c-1,c]$.
In particular, $c\ge1$ since $DD^*$ is a positive operator, and $\mu$ is uniform
measure on the annulus~\eqref{eq:ann}.
\end{proof}

\section{The distribution of $T^*T$}
\label{sec:TstT}

Throughout this section,
$(\Ac,\tau)$ will be a $*$--noncommutative probability space
and $T$ will denote a $\DT(\delta_0,1)$--element in $(\Ac,\tau)$.
We will investigate the
moments of the element $T^*T$, finding its distribution and R--transform.

Let $1\in B\subseteq\Ac$ be a unital $*$--subalgebra.
The following result is standard.
\begin{prop}
\label{prop:uBu}
If $u\in\Ac$ is a unitary such that $\tau(u)=0$ and $B$ and $\{u\}$ are $*$--free, then
$B$ and $u^*Bu$ are free.
\end{prop}

\begin{cor}
\label{cor:zbz}
Suppose $z\in\Ac$ is a circular element such that $B$ and $\{z\}$ are $*$--free.
Let $D$ be the $*$--subalgebra of $\Ac$ generated by $\{z^*bz\mid b\in B\}$.
Then $B$ and $D$ are free.
\end{cor}
\begin{proof}
By Voiculescu's result~\cite[Proposition 2.6]{V90}
about the polar decomposition of a circular element,
enlarging $\Ac$ if necessary we may without loss of generality assume $z=hu$
where $u$ is a Haar unitary, $h\ge0$ and $u$ and $h$ are $*$--free.
Then $D\subseteq u^*Cu$, where $C$ is the $*$--subalgebra of $\Ac$ generated
by $B\cup\{h\}$.
Since $C$ and $\{u\}$ are $*$--free, it follows from Proposition~\ref{prop:uBu}
that $D$ and $C$ are free,
and therefore that $D$ and $B$ are free.
\end{proof}

\begin{lemma}
\label{lem:RQ}
Suppose $N\in\Nats$ and $(a_{ij})_{1\le i\le j\le N}$ is a $*$--free family
of circular elements in $(\Ac,\tau)$ such that $\tau(a_{ij}^*a_{ij})$
is the same for all $i$ and $j$.
In the $*$--noncommutative probability space $(M_N(\Ac),\tau\circ\tr_N)$,
consider the elements
\[
R=\left(\begin{matrix}
     0 &      0 &      0 & \cdots & 0      \\
     0 & a_{22} & a_{23} & \cdots & a_{2N} \\
     0 &      0 & a_{33} & \cdots & a_{3N} \\
\vdots & \vdots & \ddots & \ddots & \vdots \\
     0 &      0 & \cdots &      0 & a_{NN}
\end{matrix}\right)\;,
\qquad
Q=\left(\begin{matrix}
 a_{11} & a_{12} & \cdots & a_{1N} \\
      0 &      0 & \cdots &      0 \\
 \vdots & \vdots & \cdots & \vdots \\
      0 &      0 & \cdots &      0      
\end{matrix}\right)\;.
\]
Then $R$ and $Q^*Q$ are $*$--free with respect to $\tau\circ\tr_N$.
\end{lemma}
\begin{proof}
Expanding $\Ac$ if necessary, let $b_{k\ell},c_{rs}\in\Ac$ be such that
\[
(a_{ij})_{1\le i\le j\le N,}\quad(b_{k\ell})_{1\le k,\ell\le N,}\quad
(c_{rs})_{1\le r,s\le N}
\]
is a $*$--free family of circular elements in $(\Ac,\tau)$
such that $\tau(a_{ij}^*a_{ij})=\tau(b_{k\ell}^*b_{k\ell})=\tau(c_{rs}^*c_{rs})$
for all $i$, $j$, $k$, $\ell$, $r$ and $s$.
Let
\[
Y=\left(\begin{matrix}
 b_{11}& b_{12} & b_{13} &    \cdots & b_{1N} \\
 b_{21}& a_{22} & a_{23} &    \cdots & a_{2N} \\
 b_{31}& b_{32} & a_{33} &    \cdots & a_{3N} \\
\vdots & \vdots & \ddots &    \ddots & \vdots \\
 b_{N1}& b_{N2} & \cdots & b_{N-1,N} & a_{NN}
\end{matrix}\right)\;,
\qquad
Z=\left(\begin{matrix}
 a_{11} & a_{12} & \cdots & a_{1N} \\
 c_{21} & c_{22} & \cdots & c_{2N} \\
 \vdots & \vdots & \vdots & \vdots \\
 c_{N1} & c_{N2} & \cdots & c_{NN}      
\end{matrix}\right)\;.
\]
Moreover, for $1\le j\le N$ let $p_j\in M_N(\Ac)$ be the diagonal $N\times N$
matrix with $1$ in the $j$th diagonal entry and zeros elsewhere.
Then by results of Voiculescu~\cite{V90}, with respect to $\tau\circ\tr_N$,
$Y$ and $Z$ are circular and the family
\[
\{Y\},\quad\{Z\},\quad\{p_j\mid 1\le j\le N\}
\]
of sets of random variables is $*$--free.
Since $R$ belongs to the subalgebra of $\Ac$ generated by $\{Y\}\cup\{p_j\mid 1\le j\le N\}$
and since $Q=p_1Z$, $*$--freeness of $R$ and $Q^*Q$ follows from Corollary~\ref{cor:zbz}.
\end{proof}

\begin{lemma}
\label{lem:ST}
In a W$^*$--noncommutative probability space $(\Ac,\tau)$
let $(Y_{ij})_{1\le i\le j\le N}$ be a $*$--free family of circular elements
with $\tau(Y_{ij}^*Y_{ij})=1$.
For $N\in\Nats$, consider the elements
\begin{gather*}
S_N=\frac1{\sqrt N}\left(\begin{matrix}
0      & Y_{12} & Y_{13} & \cdots & Y_{1N}    \\
0      &      0 & Y_{23} & \cdots & Y_{2N}    \\
0      & \ddots & \ddots & \ddots & \vdots    \\
\vdots & \ddots & 0      &      0 & Y_{N-1,N} \\
0      & \cdots & 0      &      0 & 0
\end{matrix}\right) \\[2ex]
T_N=\frac1{\sqrt N}\left(\begin{matrix}
Y_{11} & Y_{12} & Y_{13} & \cdots & Y_{1N}    \\
0      & Y_{22} & Y_{23} & \cdots & Y_{2N}    \\
0      & \ddots & \ddots & \ddots & \vdots    \\
\vdots & \ddots & 0      & \ddots & Y_{N-1,N} \\
0      & \cdots & 0      &      0 & Y_{N,N}
\end{matrix}\right)
\end{gather*}
in the W$^*$--noncommutative probability space $(\Ac\otimes M_N(\Cpx),\tau_N)$,
where $\tau_N=\tau\otimes\tr_N$ with $\tr_N$ the normailzed trace on $M_N(\Cpx)$.
Then
\renewcommand{\labelenumi}{(\alph{enumi})}
\begin{enumerate}

\item
for all $p,N\in\Nats$ with $N\ge2$,
\[
\tau_N((S_N^*S_N)^p=\bigg(\frac{N-1}N\bigg)^{p+1}\tau_{N-1}((T_{N-1}^*T_{N-1})^p),
\]

\item
if $\nu_N$ and $\mu_N$ denote the distributions of $T_N^*T_N$ and $S_N^*S_N$,
respectively, computed with respect to $\tau_N$,
then their R--transforms satisfy
\[
R_{\nu_N}(z)=R_{\mu_N}(z)+\frac1{N(1-z)}
\]
for all $z$ in a neighborhood of $0$ in $\Cpx$.

\end{enumerate}
\end{lemma}
\begin{proof}
It is clear that $S_N$ has the same $*$--moments as
\[
S_N'=\frac1{\sqrt N}\left(\begin{matrix}
0      & Y_{11} & Y_{12} & \cdots & Y_{1,N-1}   \\
0      &      0 & Y_{22} & \cdots & Y_{2,N-1}   \\
0      & \ddots & \ddots & \ddots & \vdots      \\
\vdots & \ddots & 0      &      0 & Y_{N-1,N-1} \\
0      & \cdots & 0      &      0 & 0
\end{matrix}\right).
\]
Moreover,
\[
(S_N')^*S_N'=\frac{N-1}N\left(\begin{matrix}
0 & 0 \quad \cdots \quad 0 \\
0 \\
\vdots & \quad T_{N-1}^*T_{N-1} \\
0
\end{matrix}\right)\;.
\]
Thus
\[
\tau_N((S_N^*S_N)^p)=\tau_N(((S_N')^*S_N')^p)
=\bigg(\frac{N-1}N\bigg)^{p+1}\tau_{N-1}((T_{N-1}^*T_{N-1})^p),
\]
which proves~(a).

For~(b), write $T_N=P_N+Q_N$, where
\[
P_N=\frac1{\sqrt N}\left(\begin{matrix}
0      & 0      & \cdots & 0      \\
0      & Y_{22} & \cdots & Y_{2N} \\
0      & \ddots & \vdots          \\
0      & \cdots &      0 & Y_{N,N}
\end{matrix}\right),
\qquad
Q_N=\frac1{\sqrt N}\left(\begin{matrix}
Y_{11} & Y_{12} & \cdots & Y_{1N} \\
     0 &      0 & \cdots &      0 \\
\vdots & \vdots & \vdots & \vdots \\
     0 &      0 & \cdots &      0
\end{matrix}\right).
\]
Since $P_N^*Q_N=0$, we have
\[
T_N^*T_N=P_N^*P_N+Q_N^*Q_N.
\]
By Lemma~\ref{lem:RQ}, $P_N^*P_N$ and $Q_N^*Q_N$ are free with respect to $\tau_N$.
Moreover, $P_N^*P_N=(S_N'')^*S_N''$, where
\[
S_N''=\frac1{\sqrt N}\left(\begin{matrix}
0      & Y_{22} & Y_{23} & \cdots & Y_{2N} \\
0      &      0 & Y_{33} & \cdots & Y_{3N} \\
0      & \ddots & \ddots & \ddots & \vdots \\
\vdots & \ddots & 0      &      0 & Y_{NN} \\
0      & \cdots & 0      &      0 & 0
\end{matrix}\right).
\]
Note that $S_N''$ has the same $*$--moments as $S_N$.
Hence $P_N^*P_N$ has the same distribution as $S_N^*S_N$.
By~\cite[Theorem 1.6]{V96a}, the distribution of $Q_N^*Q_N$ with respect to $\tau_N$
is free Poisson of parameter $1/N$.
In particular, (cf.~\cite[p.\ 34]{VDN} or~\cite[Example 3.3.5]{HP}), if $\sigma$ is the distribution
of $Q_N^*Q_N$ computed with respect to $\tau_N$, then
\[
R_\sigma(z)=\frac1{N(1-z)}.
\]
This proves part~(b).
\end{proof}

\begin{lemma}
\label{lem:alphabeta}
Let $N\in\Nats$, let $\alpha_N(0)=\beta_N(0)=1$ and for $p\in\Nats$ let
\begin{align*}
\alpha_N(p)&=\frac{(p-\frac1N)(p-\frac2N)\cdots(p-\frac pN)}{(p+1)!} \\
\beta_N(p)&=\frac{(p+\frac1N)(p+\frac2N)\cdots(p+\frac pN)}{(p+1)!}\;.
\end{align*}
Then
\renewcommand{\labelenumi}{(\alph{enumi})}
\begin{enumerate}

\item
the function
\[
K_N(t)=\frac t{1-\sum_{p=0}^\infty\alpha_N(p)t^{p+1}}
\]
is defined, analytic and invertible with respect to composition in a neighborhood of $0$
in $\Cpx$, and its inverse is given by
\[
K_N^{\langle-1\rangle}(z)=z(1+\tfrac zN)^{-N}
\]
in some neighborhood of $0$ in $\Cpx$;

\item
the function
\[
L_N(t)=\frac t{1-\sum_{p=0}^\infty\beta_N(p)t^{p+1}}
\]
is defined, analytic and invertible with respect to composition in a neighborhood of $0$
in $\Cpx$, and its inverse is given by
\[
L_N^{\langle-1\rangle}(z)=z(1-\tfrac zN)^N
\]
in some neighborhood of $0$ in $\Cpx$.

\end{enumerate}
\end{lemma}
\begin{proof}
For~(a), let $g(z)=z(1+\frac zN)^{-N}$, $z\in\Cpx\backslash\{-N\}$.
Since $g(0)=0$ and $g'(0)=1$, $g$ is invertible with respect to composition
in a neighborhood of $0$.
We must show that
\[
g^{\langle-1\rangle}(t)=\frac t{1-\sum_{p=0}^\infty\alpha_N(p)t^{p+1}}
\]
in some neighborhood of $0$ in $\Cpx$, which is equivalent to showing that
\begin{equation}
\label{eq:tgt}
\frac t{g^{\langle-1\rangle}(t)}=1-\sum_{p=0}^\infty\alpha_N(p)t^{p+1}
\end{equation}
in some neighborhood of $0$.
Since $t/g^{\langle-1\rangle}(t)$ is analytic in a neighborhood of $0$,
it is of the form
\[
\frac t{g^{\langle-1\rangle}(t)}=\sum_{p=0}^\infty a_pt^p,\qquad0<|t|<r
\]
for some $r>0$.
Since $g(z)=z-z^2+O(|z|^3)$ as $z\to0$, one has
\begin{align*}
g^{\langle-1\rangle}(t)&=t+t^2+O(|t|^3) \\
\frac t{g^{\langle-1\rangle}(t)}&=1-t+O(|t|^2)
\end{align*}
as $t\to0$.
Therefore, $a_0=1$ and $a_1=-1=-\alpha_N(0)$.

For $p\ge2$, we use the formula
\begin{equation}
\label{eq:ap}
a_p=\frac1{2\pi i}\int_C\frac t{g^{\langle-1\rangle}(t)t^{p+1}}\dif t,
\end{equation}
where $C$ is the positively oriented path around a circle centered at $0$
with radius less than $r$.
If the radius of $C$ is chosen small enough, then the image, $C'$, of $C$ under
$g^{\langle-1\rangle}$ will be a smooth, simple, closed curve winding once around $0$ in 
the positive direction.
The substitution $t=g(z)$ in the integral~\eqref{eq:ap} yields
\[
a_p=\frac1{2\pi i}\int_{C'}\frac{g'(z)}{zg(z)^p}\dif z.
\]
Since
\[
\dif(g(z)^{1-p})=(1-p)\frac{g'(z)}{g(z)^p}\dif z,
\]
integrating by parts yields
\begin{alignat*}{2}
(1-p)a_p&\quad=\quad\frac1{2\pi i}\int_{C'}\frac1z\dif(g(z)^{1-p})
& &\quad=\quad-\frac1{2\pi i}\int_{C'}g(z)^{1-p}\dif\big(\tfrac1z\big) \\
&\quad=\quad\frac1{2\pi i}\int_{C'}\frac{g(z)^{1-p}}{z^2}\dif z
& &\quad=\quad\frac1{2\pi i}\int_{C'}\frac{(1+\frac zN)^{N(p-1)}}{z^{p+1}}\dif z \\
&\quad=\quad\Res\bigg(\frac{(1+\frac zN)^{N(p-1)}}{z^{p+1}}\;;\;z=0\bigg)
& &\quad=\quad\binom{N(p-1)}p\frac1{N^p},
\end{alignat*}
where the last quantity is the coefficient of $z^p$ in the polynomial expansion of
$(1+\frac zN)^{N(p-1)}$.
Thus for $p\ge1$,
\[
a_{p+1}=-\frac1p\binom{Np}{p+1}\frac1{N^{p+1}}=
-\frac{(Np-1)(Np-2)\dots(Np-p)}{(p+1)!\,N^p}=-\alpha_N(p).
\]
We have shown $a_0=1$ and $a_{p+1}=-\alpha_N(p)$ for all $p\ge0$,
which proves equation~\eqref{eq:tgt} and finishes the proof of part~(a).

Part~(b) follows by minor modifications of the above proof.
Put
\[
h(z)=z(1-\tfrac zN)^N
\]
and let
\[
\frac t{h^{\langle-1\rangle}(t)}=\sum_{p=0}^\infty b_pt^p
\]
be the power series expansion, which is
valid for $0<|t|<s$ for some $s>0$.
As before, one finds
$b_0=1$, $b_1=-1=-\beta_N(0)$ and, for $p\ge2$,
\[
(1-p)b_p\quad=\quad\Res\bigg(\frac{(1-\frac zN)^{-N(p-1)}}{z^{p+1}}\;;\;z=0\bigg)
\quad=\quad\binom{-N(p-1)}p\frac{(-1)^p}{N^p}\;,
\]
where the latter quantity is the coefficient of $z^p$ in the power series expansion of
$(1-\frac zN)^{-N(p-1)}$.
Thus, for $p\ge1$,
\[
b_{p+1}=\frac{(-1)^p}{pN^{p+1}}\binom{-Np}{p+1}
=-\frac{(Np+1)(Np+2)\dots(Np+p)}{(p+1)!\,N^p}=-\beta_N(p),
\]
which shows that
\[
\frac t{h^{\langle-1\rangle}(t)}=1-\sum_{p=0}^\infty\beta_N(p)t^{p+1}
\]
for $0<|t|<s$.
This proves part~(b).
\end{proof}

\begin{lemma}
\label{lem:STmom}
Let $S_N$, $T_N$, $\tau_N$, $\mu_N$ and $\nu_N$ be as in Lemma~\ref{lem:ST}, ($N\in\Nats$).
Then
\renewcommand{\labelenumi}{(\roman{enumi})}
\begin{enumerate}

\item
for all $p\in\Nats$,
\[
\tau_N((S_N^*S_N)^p)=\frac{(p-\frac1N)(p-\frac2N)\cdots(p-\frac pN)}{(p+1)!}\;,
\]

\item
for all $p\in\Nats$,
\[
\tau_N((T_N^*T_N)^p)=\frac{(p+\frac1N)(p+\frac2N)\cdots(p+\frac pN)}{(p+1)!},
\]

\item for all $z$ in a neighborhood of $0$ in $\Cpx$,
\[
R_{\mu_N}(z)=\frac1{N(1-z)((1-z)^{-1/N}-1)}-\frac1z,
\]

\item for all $z$ in a neighborhood of $0$ in $\Cpx$,
\[
R_{\nu_N}(z)=\frac1{N(1-z)(1-(1-z)^{1/N})}-\frac1z.
\]

\end{enumerate}
\end{lemma}
\begin{proof}
We proceed by induction on $N$.
We have $S_1^*S_1=0$, so~(i) and~(iii) hold when $N=1$.
Moreover, $(T_1^*T_1)^{1/2}=(Y_{11}^*Y_{11})^{1/2}$ is a quarter--circular
element, hence
\[
\tau_N((T_1^*T_1)^p)=\frac1\pi\int_0^2x^{2p}\sqrt{4-x^2}\dif x
=\frac1{p+1}\binom{2p}p,
\]
(the $p$th Catalan number), which proves~(ii) when $N=1$.
But $T_1^*T_1$ has the free Poisson distribution with parameter $1$, which
implies that $R_{\nu_N}=1/(1-z)$, proving~(iv) when $N=1$.

Now suppose $N\ge2$.
By the induction hypothesis, for every $p\in\Nats$ we have
\[
\tau_{N-1}((T_{N-1}^*T_{N-1})^p)=\frac{(p+\frac1{N-1})(p+\frac2{N-1})\dots(p+\frac p{N-1})}{(p+1)!}.
\]
Then by part~(a) of Lemma~\ref{lem:ST},
\begin{align*}
\tau_N((S_N^*&S_N)^p)
=\bigg(\frac{N-1}N\bigg)^{p+1}
 \frac{(p+\frac1{N-1})(p+\frac2{N-1})\dots(p+\frac p{N-1})}{(p+1)!} \\[1ex]
&=\bigg(\frac{N-1}{N^{p+1}}\bigg)
 \frac{\big((N-1)p+1\big)\big((N-1)p+2\big)\dots\big((N-1)p+(p-1)\big)Np}{(p+1)!} \\[1ex]
&=\frac{(pN-p)(pN-(p-1))\dots(pN-1)}{N^p(p+1)!} \\[1ex]
&=\frac{(p-\frac pN)(p-\frac{p-1}N)\dots(p-\frac1N)}{(p+1)!}.
\end{align*}
This proves~(i) for this particular $N$.

We now show (i)$\implies$(iii).
By~(i), the Cauchy transform of $\mu_N$ is given by
\[
G_{\mu_N}(\lambda)=\sum_{p=0}^\infty\alpha_N(p)\lambda^{-p-1}
\]
for $\lambda\in\Cpx$ with $|\lambda|$ large, where $\alpha_N(p)$ is as defined
in Lemma~\ref{lem:alphabeta}.
Thus
\[
K_N(t)=\frac t{1-G_{\mu_N}(1/t)}
\]
for $t$ in some neighborhood of $0$, where $K_N$ is as defined
in Lemma~\ref{lem:alphabeta}.
Since by Lemma~\ref{lem:alphabeta},
\[
K^{\langle-1\rangle}_N(w)=w(1+\frac wN)^{-N}
\]
in a neighborhood of $0$, we have for $|w|$ small
\[
w=K_N(K^{\langle-1\rangle}_N(w))=\frac{w(1+\frac wN)^{-N}}{1-G_{\mu_N}\big(\frac1w(1+\frac wN)^N\big)}.
\]
Hence for $|w|$ small,
\begin{equation}
\label{eq:GfrK}
G_{\nu_N}\bigg(\frac1w\bigg(1+\frac wN\bigg)^N\bigg)=1-\bigg(1+\frac wN\bigg)^{-N}.
\end{equation}
Letting $z=1-(1+\frac wN)^{-N}$ for small values of $|w|$, we have
$w=N((1-z)^{-1/N}-1)$ and $(1+\frac wN)^N=1/(1-z)$.
Thus for small values of $|z|$, we have
\[
G_{\nu_N}\bigg(\frac1{N((1-z)^{-1/N}-1)(1-z)}\bigg)=z.
\]
Therefore
\[
R_{\mu_N}(z)=G_{\mu_N}^{\langle-1\rangle}(z)-\frac1z
=\frac1{N(1-z)((1-z)^{-1/N}-1)}-\frac1z,
\]
which proves~(iii) for this particular $N$.

Part~(iii), together with part~(b) of Lemma~\ref{lem:ST}, proves~(iv)
for this particular $N$.

Finally, we prove~(ii) (still for this $N$).
Since there is a one--to--one correspondence between moment series and R--series,
it will suffice to prove that if~(ii) holds, then the corresoponding R--transform
$R_{\nu_N}$
is given by~(iv).
This in turn follows by minormodifications of the above proof of (i)$\implies$(iii).
Indeed, if~(ii) holds then the Cauchy transform of $\nu_N$ is
\[
G_{\nu_N}(\lambda)=\sum_{p=0}^\infty\beta_N(p)\lambda^{-p-1}.
\]
Thus
\[
L_N(t)=\frac t{1-G_{\nu_N}(1/t)}
\]
for $t$ in a neighborhood of $0$, where $L_N$ is the function defined in Lemma~\ref{lem:alphabeta}.
Thus
\[
L_N^{\langle-1\rangle}(w)=w\bigg(1-\frac wN\bigg)^N
\]
for $|w|$ small, and arguing as in the proof of (i)$\implies$(iii) shows
\[
G_{\nu_N}\Bigg(\frac1w\bigg(1-\frac wN\bigg)^{-N}\Bigg)=1-\bigg(1-\frac wN\bigg)^N
\]
for $|w|$ small.
Setting $z=1-(1-\frac wN)^N$ yields $w=N(1-(1-z)^{1/N})$ and
\[
G_{\nu_N}\bigg(\frac1{N(1-(1-z)^{1/N})(1-z)}\bigg)=z,
\]
proving~(iv), as desired.
This concludes the proof of the induction step, and of the lemma.
\end{proof}

Recall that $T$ is a $\DT(\delta_0,1)$--operator in a W$^*$--noncommutative probability space
$(\Ac,\tau)$.

\begin{thm}
\label{thm:TstT}
\renewcommand{\labelenumi}{(\alph{enumi})}
\begin{enumerate}

\item
For every $p\in\Nats$,
\[
\tau((T^*T)^p)=\frac{p^p}{(p+1)!}.
\]

\item
If $\nu$ is the distribution of $T^*T$, then its R--transform is given by
\[
R_\nu(z)=\frac{-1}{(1-z)\log(1-z)}-\frac1z
\]
for $|z|$ small.

\end{enumerate}
\end{thm}
\begin{proof}
Without loss of generality assume $(\Ac,\tau)$ is a W$^*$--noncommutative probability space with $\tau$ faithful.
By Theorem~\ref{thm:DTfreematrix}, enlarging $(\Ac,\tau)$ if necessary,
$S_N$ differs from a $DT(\delta_0,1)$--operator by a diagonal matrix
\[
A=\diag(A_1,\dots,A_N)\in\Ac\otimes M_N(\Cpx),
\]
where each $A_j$ is $\DT(\delta_0,1/\sqrt N)$.
Thus $\|A\|\le2/\sqrt N$.
Also, $T_N$ differs from $S_N$ by a diagonal matrix of norm $2/\sqrt N$,
and we conclude that $T_N$ converges in $*$--moments to $T$ as $N\to\infty$.
Taking the limit as $N\to\infty$ of the value of $\tau_N((T_N^*T_N)^p)$ given in part~(ii)
of Lemma~\ref{lem:STmom} yields~(a) above.

The convergence in $*$--moments of $T_N$ to $T$ ensures that
the coefficients of the R--series of the distribution $\nu_N$ of $T_N^*T_N$
converge to the coefficients of the R--series of the distribution $\nu$ of $T^*T$.
Writing $(1-z)^{1/N}=\exp(\log(1-z)/N)$ and
taking the limit as $N\to\infty$ of the formula in
part~(iv) of Lemma~\ref{lem:STmom} proves part~(b).
\end{proof}

\begin{lemma}
\label{lem:L}
Let $\gamma(0)=1$ and for $p\in\Nats$ let
\[
\gamma(p)=\frac{p^p}{(p+1)!}.
\]
Then the function
\[
L(t)=\frac t{1-\sum_{p=0}^\infty\gamma(p)t^{p+1}}
\]
is defined, analytic and invertible with respect to composition in a neighborhood of $0$ in $\Cpx$,
and its inverse is given by
\begin{equation}
\label{eq:L-1}
L^{\langle-1\rangle}(z)=ze^{-z}
\end{equation}
for $|z|$ small.
\end{lemma}
\begin{proof}
It is clear that $L$ is defined and analytic in a neighborhood of $0$.
Since $L(0)=0$ and $L'(0)=1$, $L$ is invertible in a neighborhood of $0$.
Turning now to Lemma~\ref{lem:alphabeta}, since $\lim_{N\to\infty}\beta_N(p)=\gamma(p)$
for all $p\in\Nats\cup\{0\}$,
it follows that the coefficients in the power series expansion for $L_N(t)$ around $0$
converge as $N\to\infty$
termwise to the ceoffiecients in the power series expansion of $L(t)$.
Hence the coefficients in the power series for $L^{\langle-1\rangle}_N(z)$
converge to those in the power series expansion $L^{\langle-1\rangle}(z)$.
Since $L_N^{\langle-1\rangle}(z)=z(1-\frac zN)^N$,
we get the formula~\eqref{eq:L-1}.
\end{proof}

\begin{thm}
\label{thm:phi}
The distribution $\nu$ of $T^*T$ is absolutely continuous with respect to Lebesgue measure
and has support equal to the interval $[0,e]$.
Its density function $\phi$ is defined on the interval $(0,e)$ by
\begin{equation}
\label{eq:phidef}
\phi\bigg(\frac{\sin v}v\exp(v\cot v)\bigg)
=\frac1\pi\sin v\exp(-v\cot v),\qquad0<v<\pi.
\end{equation}
\end{thm}
\begin{proof}
By part~(a) of Theorem~\ref{thm:TstT}, the Cauchy transform of $\nu$ is
\[
G_\nu(\lambda)=\sum_{p=0}^\infty\gamma(p)\lambda^{-p-1}.
\]
Using Lemma~\ref{lem:L} and the technique used in the proof of equation~\eqref{eq:GfrK}
in Lemma~\ref{lem:STmom}, we find
\[
G_\nu(w^{-1}e^w)=1-e^{-w}
\]
for $|w|$ small.
Define $\rho:(0,\pi)\to\Reals$ by
\begin{equation}
\label{eq:rhodef}
\rho(v)=\frac{\sin v}v\exp(v\cot v).
\end{equation}
Then $\rho(v)>0$ and $\rho$ is strictly decreasing, because for $0<v<\pi$,
\begin{align*}
\frac{\dif}{\dif v}\log\rho(v)&=2\cot v-\frac1v-\frac v{\sin^2 v} \\
&=-\bigg(\frac{(v-\sin v)^2+2v\sin v(1-\cos v)}{v\sin^2 v}\bigg)\quad<0.
\end{align*}
Moreover,
\[
\lim_{v\to0^+}\rho(v)=e,\qquad\lim_{v\to\pi^-}\rho(v)=0,
\]
so $\rho$ is an order--reversing bijection from $(0,\pi)$ onto $(0,e)$.
Hence we may define
\[
\phi:(0,e)\to(0,\infty)
\]
by equation~\eqref{eq:phidef}.
Using substitution and integration by parts, we have
\begin{align*}
\int_0^e\phi(x)\dif x
&=-\int_0^\pi\phi(\rho(v))\rho^{\,\prime}(v)\dif v \\
&=-\int_0^\pi\frac1\pi\sin v\exp(-v\cot v)\dif\bigg(\frac{\sin v}v\exp(v\cot v)\bigg) \\
&=\frac{-\sin^2 v}{\pi v}\bigg|^\pi_0
 +\frac1\pi\int_0^\pi\frac{\sin v}v\exp(v\cot v)\dif\big(\sin v\exp(-v\cot v)\big) \\
&=\frac1\pi\int_0^\pi\bigg(\frac{\sin v}v\bigg)\bigg(\frac v{\sin v}\bigg)\dif v=1.
\end{align*}
So $\phi$ is the density of a probability measure on $(0,e)$.

It remains to prove that for all $p\in\Nats$,
\[
\int_0^ex^p\phi(x)\dif x=\gamma(p)=\frac{p^p}{(p+1)!}.
\]
Let $g(w)=we^{-w}$.
By Lemma~\ref{lem:L}, the inverse of $g$ in a neighborhood of $0$ is
\[
g^{\langle-1\rangle}(t)=\frac t{1-\sum_{p=0}^\infty\gamma(p)t^{p+1}}.
\]
As in the proof of Lemma~\ref{lem:alphabeta}, $\gamma(p)$ can be recovered for $p\ge1$ as
$\gamma(p)=-c_{p+1}$ where
\[
c_{p+1}=\frac1{2\pi i}\int_{C'}\frac{g'(w)}{wg(w)^{p+1}}\dif w,
\]
for some smooth simple closed path $C'$ winding once around $0$
with positive orientation.
In our case, $g(w)=we^{-w}$ and $g'(w)=(1-w)e^{-w}$, so
\begin{equation}
\label{eq:gamp}
\gamma(p)=-c_{p+1}=\frac1{2\pi i}\int_{C'}\frac{w-1}{w^{p+2}}e^{pw}\dif w.
\end{equation}
Since the integrand is analytic in $\Cpx\backslash\{0\}$, we may choose $C'$ to be
any piecewise smooth closed path winding once around $0$ with positive orientation.
Let $\alpha\in(\frac\pi2,\pi)$ and put $C'=C_1^{(\alpha)}\cup C_2^{(\alpha)}$,
where $C_1^{(\alpha)}$ and $C_2^{(\alpha)}$ are parameterized by
\begin{alignat*}{3}
C_1^{(\alpha)}&: &\quad v &\mapsto v\cot v+iv,          &\qquad -\alpha\le&v\le\alpha \\
C_2^{(\alpha)}&: &\quad v &\mapsto \alpha\cot\alpha-iv, &\qquad -\alpha\le&v\le\alpha,
\end{alignat*}
setting $v\cot v=1$ when $v=0$.
Since $\alpha\cot\alpha\to-\infty$ as $\alpha\to\pi^-$, it is clear that
\[
\bigg|\int_{C_2^{(\alpha)}}\frac{w-1}{w^{p+2}}e^{pw}\dif w\bigg|\to0
\]
as $\alpha\to\pi^-$, for all $p\ge1$.
Hence
\[
\gamma(p)=\frac1{2\pi i}\int_{C_1}\frac{w-1}{w^{p+2}}e^{pw}\dif w,\qquad(p\in\Nats),
\]
where $C_1$ is the path parametrized by
\[
v\mapsto w=v\cot v+iv=\frac v{\sin v}e^{iv},\qquad-\pi<v<\pi,
\]
with the convention $v\cot v=1=v/\sin v$ when $v=0$.
We have
\begin{align*}
\bigg(\frac1we^w\bigg)^p&=\bigg(\frac{\sin v}v\exp(v\cot v)\bigg)^p \\[1ex]
\frac{w-1}{w^2}e^w\dif w&=\dif\bigg(\frac1we^w\bigg)=\dif\bigg(\frac{\sin v}v\exp(v\cot v)\bigg) \\[1ex]
e^{-w}&=e^{-v\cot v-iv},
\end{align*}
which give
\[
\gamma(p)=\frac1{2\pi i}\int_{-\pi}^\pi e^{-v\cot v-iv}\rho(v)^p\rho^{\,\prime}(v)\dif v,
\]
where $\rho$ is as defined in equation~\eqref{eq:rhodef}.
Since $\rho$ is an even function, $\rho^{\,\prime}$ is odd and using equation~\eqref{eq:phidef} we get
\begin{align*}
\gamma(p)&=\frac{-1}{2\pi}\int_{-\pi}^\pi e^{-v\cot v}(\sin v)\rho(v)^p\rho^{\,\prime}(v)\dif v \\[1ex]
&=\frac{-1}\pi\int_0^\pi(\sin v)e^{-v\cot v}\rho(v)^p\rho^{\,\prime}(v)\dif v
=\int_0^e\phi(x)x^p\dif x,
\end{align*}
as required.
\end{proof}

\begin{remark}\rm
\label{rem:phi0}
From Theorem~\ref{thm:phi} one gets
\begin{alignat*}{2}
\phi(x)&\sim\frac1{x(\log x)^2}&\quad\text{as }&x\to0^+ \\
\phi(x)&\sim\frac{\sqrt2}{\pi e^{3/2}}(e-x)^{1/2}&\quad\text{as }&x\to e^-.
\end{alignat*}
Note that the asymptotic behevior as $x\to0^+$ implies
\[
\int_0^e\phi(x)\log x\dif x=-\infty.
\]
This equality also follows from the fact that the Fuglede--Kadison determinant
of any quasi--nilpotent operator in a II$_1$--factor is equal to $0$.
\end{remark}

\begin{cor}
\label{cor:nmT}
$\|T\|=\sqrt e$.
\end{cor}

\begin{cor}
\label{cor:kerT}
$\ker(T)=\ker(T^*T)=\{0\}$.
\end{cor}
\begin{proof}
This follows from $\nu(\{0\})=0$, where $\nu$ is the distribution of $T^*T$.
\end{proof}

\section{$*$--moments of the operator $T$}
\label{sec:Tmom}
\setcounter{table}{0}

In this section, $T$ will denote a $\DT(\delta_0,1)$--element.
Algorithm~\ref{alg:NTO} allows computation of an arbitrary $*$--moment of $T$,
though for large moments this becomes somewhat arduous.
We will prove a
recursion formula for general $*$--moments of $T$.
The proof will use the construction $T=\UT(X,\lambda)$, (see Definition~\ref{def:UT}),
and will involve asymptotic behavior of certain products of free semicircular
variables.

We begin with a few preparatory lemmas.

\begin{lemma}
\label{lem:XtYt}
In a W$^*$--noncommutative probability space $(\Ac,\tau)$, let $B\subseteq\Ac$
be a unital, commutative W$^*$--subalgebra and let $\Xt$ and $\Yt$ be semicircular
elements having the same second moments and such that the family
$B,\{\Xt\},\{\Yt\}$ is free.
Let $p\in B$ be a projection and let
\[
X=p\Xt p+p\Yt(1-p)+(1-p)\Yt p+(1-p)\Yt(1-p).
\]
Then $X$ is semicircular, has the same second moment as $\Xt$ and $\Yt$,
and the pair $\{X\},B$ is free.
\end{lemma}
\begin{proof}
We may without loss of generality assume that $\tau$ is a faithful trace,
$B\cong L^\infty[0,1]$ and $\tau(\Xt^2)=1=\tau(\Yt^2)$.
Let $k\in\Nats$ and let $B_k\subseteq B$ be the linear span of a set of $k$
orthogonal projections in $B$, each having trace $1/k$.
We will show that if $\pt\in B_k$ is a projection, $0<\tau(\pt)<1$,
and if
\[
Z=\pt\Xt\pt+\pt\Yt(1-\pt)+(1-\pt)\Yt\pt+(1-\pt)\Yt(1-\pt),
\]
then $Z$ is semicircular with $\tau(Z^2)=1$ and the pair $\{Z\},B$ is free.
By approximation, this will suffice to prove the lemma.

Let $(D,\tau_D)$ be a W$^*$--noncommutative probability space having a $*$--free
family $(x_{ij})_{1\le j\le k},\,(y_{ij})_{1\le j\le k}$ where for all $i$,
$x_{ii}$ and $y_{ii}$ are semicircular elements with $\tau_D(x_{ii}^2)=\tau_D(y_{ii}^2)=1/k$
and for all $i\ne j$, $x_{ij}$ and $y_{ij}$ are circular elements with
$\tau_D(x_{ij}^*x_{ij})=\tau_D(y_{ij}^*y_{ij})=1/k$.
By results~\cite{V:RM} of Voiculescu arising from his matrix model,
we may take, in the W$^*$--noncommutative probability space $(M_k(D),\tau_D\circ\tr_k)$,
$\Xt=(x_{ij})_{1\le i,j\le n}$, $\Yt=(y_{ij})_{1\le i,j\le n}$, $B_k$ the set of diagonal
$k\times k$ matrices with scalar entries and
$\pt=\diag(1,\ldots,1,0,\ldots,0)$ with $\tr_k(\pt)=\ell/k$.
Then
\[
Z=\left(\begin{matrix}
x_{11} & \cdots & x_{1\ell} & y_{1,\ell+1} & \cdots & y_{1k} \\
\vdots & \vdots & \vdots    & \vdots       & \vdots & \vdots \\
x_{\ell1} & \cdots & x_{\ell\ell} & y_{\ell,\ell+1} & \cdots & y_{\ell k} \\
y_{\ell+1,1} & \cdots & y_{\ell+1,\ell} & y_{\ell+1,\ell+1} & \cdots & y_{\ell+1,k} \\
\vdots & \vdots & \vdots    & \vdots       & \vdots & \vdots \\
y_{k1} & \cdots & y_{k\ell} & y_{k,\ell+1} & \cdots & y_{kk}
\end{matrix}\right).
\]
By these same results in~\cite{V:RM}, $Z$ is semicircular and the pair $B_k,\{Z\}$ is free.
\end{proof}

\begin{lemma}
\label{lem:1pXp}
Let $(\Ac,\tau)$ be a W$^*$--noncommutative probability space, let $B\subseteq\Ac$
be a unital, commutative W$^*$--subalgebra and let $X\in \Ac$ be a semicircular
element such that the pair $B,\{X\}$ is free.
Let $p\in B$ be a projection, $0<\tau(p)<1$.
Then the family
\begin{equation}
\label{eq:pXBfree}
pB,\quad\{pXp\},\quad pW^*(\{X(1-p),(1-p)X\}\cup(1-p)B)p
\end{equation}
is free with respect to $\tau(p)^{-1}\tau\restrict_{p\Ac p}$.
\end{lemma}
\begin{proof}
Without loss of generality we may take $X=p\Xt p+p\Yt(1-p)+(1-p)\Yt p+(1-p)\Yt(1-p)$
with $\Xt$, $\Yt$ and $B$ as in Lemma~\ref{lem:XtYt} and with $\tau$ a faithful trace.
Let $D=\Cpx p+\Cpx(1-p)$ and let $E_D:\Ac\to D$ be the $\tau$--preserving conditional expectation.
Then the family
\begin{equation}
\label{eq:wD}
B,\quad W^*(\{\Xt\}\cup D),\quad W^*(\{\Yt\}\cup D)
\end{equation}
is free over $D$ with respect to $E_D$.
Let us show that the family
\begin{equation}
\label{eq:wDfree}
pB,\quad pW^*(\{\Xt\}\cup D)p,\quad pW^*(\{\Yt\}\cup(1-p)B)p
\end{equation}
is free with respect to $\tau(p)^{-1}\tau\restrict_{p\Ac p}$.
Let
\begin{align*}
\YEu\oup&=W^*(\{\Yt\}\cup D)\cap\ker E_D, \\
\XEu\oup&=W^*(\{\Xt\}\cup D)\cap\ker E_D, \\
B\oup&=B\cap\ker E_D.
\end{align*}
Let $C=\Cpx p+(1-p)B$ and $C\oup=C\cap B\oup$.
By freeness of~\eqref{eq:wD} over $D$,
$W^*(\{\Yt\}\cup(1-p)B)\cap\ker E_D$ is the weak closure of the
linear span of $\Lambdao(\YEu\oup,C\oup)$.
(The definition of this notation can be found near equation~\eqref{eq:Lambdao}.)
From this and the fact that $p$ is a minimal projection in $C$ and $D$,
we see that $pW^*(\{\Yt\}\cup(1-p)B)p\cap\ker\tau$ is the weak
closure of the linear span of $p\Theta p$, where $\Theta$ is the set of
all words belonging to $\Lambdao(\YEu\oup,C\oup)$ whose first and last letters
come from $\YEu\oup$.
Therefore, to prove freeness of~\eqref{eq:wDfree}, it will suffice to show
\[
\Lambdao\big(p\Theta p,\,pW^*(\{\Xt\}\cup D)p\cap\ker\tau,\,pB\cap\ker\tau\big)\subseteq\ker\tau.
\]
Since $p$ is a minimal projection in $D$,
\[
pW^*(\{\Xt\}\cup D)p\cap\ker\tau\subseteq\XEu\oup.
\]
Using this and the freeness of~\eqref{eq:wD} over $D$, we have
\begin{align*}
\Lambdao\big(p\Theta p,\,pW^*(\{\Xt\}\cup D)&p\cap\ker\tau,\,pB\cap\ker\tau\big)\subseteq \\[0.5ex]
&\subseteq p\Lambdao\big(\Theta,\,pW^*(\{\Xt\}\cup D)p\cap\ker\tau,\,pB\cap\ker\tau\big)p \\[0.8ex]
&\subseteq p\Lambdao(\YEu\oup,\,\XEu\oup,\,B\oup)p\subseteq p(\ker E_D)p\subseteq\ker\tau,
\end{align*}
so the family~\eqref{eq:wDfree} is free.
However,
\begin{align*}
pXp&\in pW^*(\{\Xt\}\cup D)p \\
(1-p)X,\,(1-p)X&\in W^*(\{\Yt\}\cup(1-p)B)
\end{align*}
and the freeness of the family~\eqref{eq:pXBfree} follows from that
of the family~\eqref{eq:wDfree}.
\end{proof}

\begin{lemma}
\label{lem:Xpasympt}
In a noncommutative probability space $(\Qc,\tau)$, let $X$ be a semicircular
element satisfying $\tau(X^2)=1$ and for $n\in\Nats$ define the function
$h_n:(0,1)\to[0,\infty]$ by
\[
h_n(t)=\tau((Xp_tX(1-p_t))^n),
\]
where $p_t\in\Ac$ is an idempotent such that $\tau(p_t)=t$ and the pair $X,\,p_t$
is free.
Then the asymptotic behavior of $h_n(t)$ as $t\to0$ is
\begin{equation}
\label{eq:hnas}
h_n(t)=t+O(t^2).
\end{equation}
\end{lemma}
\begin{proof}
We will use Speicher's moment--cummulant formula~\cite[Theorem 2.17]{Sp:Gr},
(see alternatively equation~(3.14) of~\cite{NS}),
which gives
\begin{equation}
\label{eq:momcum}
\tau\big((Xp_tX(1-p_t))^n\big)=\sum_{\pi\in\NC(2n)}
k_\pi[X,X,\ldots,X]\phi_{K(\pi)}[p_t,1-p_t,\ldots,p_t,1-p_t].
\end{equation}
Here the sum is over all non--crossing partitions $\pi$ of $\{1,\ldots,2n\}$,
$k_\pi[X,\ldots,X]$ is the corresponding cummulant and $K(\pi)\in\NC(2n)$
is the Kreweras complement of $\pi$.
Since
\[
k_m[X,\ldots,X]=\begin{cases}1&\text{if }m=2 \\0&\text{otherwise},\end{cases}
\]
the sum~\eqref{eq:momcum} becomes
\[
\tau\big((Xp_tX(1-p_t))^n\big)=\sum_{\pi\in\NC_2(2n)}\phi_{K(\pi)}[p_t,1-p_t,\ldots,p_t,1-p_t],
\]
where $\NC_2(2n)$ is the set of non--crossing pairings
of $\{1,\ldots,2n\}$.
For any $\pi\in\NC_2(2n)$, $K(\pi)$ will never have a block containing
both even and odd numbers, and
\[
\phi_{K(\pi)}[p_t,1-p_t,\ldots,p_t,1-p_t]=t^a(1-t)^b,
\]
where $a$ is the number of blocks of $K(\pi)$ containing odd numbers
and $b$ is the number of blocks containing even numbers.
Thus, if $a\ge2$, then $\phi_{K(\pi)}[p_t,1-p_t,\ldots,p_t,1-p_t]=O(t^2)$ as $t\to0$.
There is exactly one non--crossing pairing $\pi$ of $\{1,\ldots,2n\}$
such that $K(\pi)$ groups all odd numbers into one block, namely this one:
\begin{center}
\begin{picture}(225,50)(0,0)

\put(0,40){$Xp_tX(1-p_t)Xp_tX(1-p_t)X\cdots p_tX(1-p_t)$.}

\drawline(4,30)(4,0)(182,0)(182,30)
\drawline(24,30)(24,15)(74,15)(74,30)
\drawline(94,30)(94,15)(144,15)(144,30)

\put(155,20){$\cdots$}

\end{picture}
\end{center}
This yields $t$ in the asymptotic expansion~\eqref{eq:hnas} of $h_n(t)$.
\end{proof}

Recall we let $T$ be a $\DT(\delta_0,1)$--element in a $*$--noncommutative probability space
$(\Ac,\tau)$.
For $n\in\Nats$ and $k_1,\ldots,k_n,\ell_1,\ldots,\ell_n\in\Nats\cup\{0\}$,
let
\begin{equation}
\label{eq:Mkl}
M(k_1,\ell_1,\ldots,k_n,\ell_n)=\tau((T^*)^{k_1}T^{\ell_1}\ldots(T^*)^{k_n}T^{\ell_n}).
\end{equation}
The following properties of these moments are easily seen
from Lemma~\ref{lem:utmom} and Algorithm~\ref{alg:NTO}, and the fact (Proposition~\ref{prop:constadj})
that also $T^*$ is a $\DT(\delta_0,1)$--element.
\begin{prop}
\label{prop:Mprop}
Let $n\in\Nats$ and $k_1,\ldots,k_n,\ell_1,\ldots,\ell_n\in\Nats\cup\{0\}$.
\renewcommand{\labelenumi}{(\roman{enumi})}
\begin{enumerate}

\item If $k_1+\cdots k_n\ne\ell_1+\cdots+\ell_n$ then $M(k_1,\ell_1,\ldots,k_n,\ell_n)=0$.

\item $M(k_1,\ell_1,\ldots,k_n,\ell_n)=M(\ell_1,k_2,\ell_2,\ldots,k_n,\ell_n,k_1)$.

\item $M(k_1,\ell_1,\ldots,k_n,\ell_n)=M(\ell_n,k_n,\ell_{n-1},k_{n-1},\ldots,\ell_1,k_1)$.

\item If $k_1=0$ and $n\ge2$ then
\[
M(k_1,\ell_1,\ldots,k_n,\ell_n)=M(k_2,\ell_2,\ldots,k_{n-1},\ell_{n-1},k_n,\ell_n+\ell_1).
\]

\end{enumerate}
\end{prop}

We are now ready to state and prove the recursion formula.
\begin{thm}
\label{thm:Tmom}
Let  $n\in\Nats$ and $k_1,\ldots,k_n,\ell_1,\ldots,\ell_n\in\Nats$, (all nonzero),
be such that $m:=k_1+\cdots+k_n=\ell_1+\cdots+\ell_n$.
Then
\begin{align}
M(k_1,&\ell_1,\ldots,k_n,\ell_n)=\frac1{m+1}\sum_{r=1}^n
\sum_{1\le j(1)<\cdots<j(r)\le n} \label{eq:recur} \\
&\begin{aligned}[t]
\bigg(
M&(k_1,\ell_1,\ldots,k_{j(1)-1},\ell_{j(1)-1},
k_{j(1)}-1,\ell_{j(r)}-1,k_{j(r)+1},\ell_{j(r)+1},\ldots,k_n,\ell_n) \\
&\prod_{i=1}^{r-1}M(\ell_{j(i)}-1,k_{j(i)+1},\ell_{j(i)+1},
\ldots,k_{j(i+1)-1},\ell_{j(i+1)-1},k_{j(i+1)}-1)\bigg).
\end{aligned} \notag
\end{align}
\end{thm}
\begin{proof}
We assume $(\Ac,\tau)=(L(\Fb_2),\tau)$ and let
$T=\UT(X,\lambda)$ with $X$ and $\lambda$ as in Lemma~\ref{lem:Tconstr}.
We will let $B$ denote the image of $\lambda$.
For $0<t<1$ let $p_t=\lambda(1_{[0,t]})\in B$.
Then $\tau(p_t)=t$.
Let $S_t=(1-p_t)T$ and $Q_t=p_tT$.
From Lemmas~\ref{lem:Tbasics} and~\ref{lem:qTq}, we have that
\renewcommand{\labelenumi}{(\alph{enumi})}
\begin{enumerate}

\item $(1-t)^{-1/2}S_t$ is an element of $W^*((1-p_t)X(1-p_t)\cup(1-p_t)B)$
and is a $\DT(\delta_0,1)$--operator in
$((1-p_t)\Ac(1-p_t),(1-t)^{-1}\tau\restrict_{(1-p_t)\Ac(1-p_t)})$

\item $t^{-1/2}Q_tp_t=t^{-1/2}Tp_t$ is an element of $W^*(p_tXp_t\cup p_tB)$
and is a $\DT(\delta_0,1)$--operator in
$(p_t\Ac p_t,t^{-1}\tau\restrict_{p_t\Ac p_t})$.

\end{enumerate}
Substituting $T=S_t+Q_t$ into~\eqref{eq:Mkl} and distributing,
$M(k_1,\ell_1,\ldots,k_n,\ell_n)$ is written as a sum of the $2^{2m}$ terms
obtained by substituting $S_t$ or $Q_t$ variously for each $T$ in~\eqref{eq:Mkl}.
One of these terms is
\begin{equation}
\label{eq:Sterm}
\tau((S_t^*)^{k_1}S_t^{\ell_1}\cdots(S_t^*)^{k_n}S_t^{\ell_n}).
\end{equation}
Using the observation~(a) about $S_t$, we see that the term~\eqref{eq:Sterm}
equals
\[
(1-t)^{m+1}M(k_1,\ell_1,\ldots,k_n,\ell_n).
\]
Thus
\[
\big(1-(1-t)^{m+1}\big)M(k_1,\ell_1,\ldots,k_n,\ell_n)
\]
is equal to the sum of the remaining $2^{2m}-1$ terms, in each of which $Q_t$ or its adjoint appears at least once.
As $t\to0$, we have the asymptotic behavior $1-(1-t)^{m+1}=(m+1)t+O(t^2)$.
We will show that the asympototic
behavior as $t\to0$ of each of the remaining $2^{2m}-1$ terms is $ct+O(t^{3/2})$ for some
constant $c$.
The value of $M(k_1,\ell_1,\ldots,k_n,\ell_n)$ will then be the sum of all these constants,
divided by $m+1$.

Since $p_tS_t=0=S_tp_t$, any of these $2^{2m}-1$ terms is zero unless it is of the form
\begin{equation}
\label{eq:SQform}
\tau((S_t^*)^{k_1-a_1}(Q_t^*)^{a_1}Q_t^{b_1}S_t^{\ell_1-b_1}
\cdots(S_t^*)^{k_n-a_n}(Q_t^*)^{a_n}Q_t^{b_n}S_t^{\ell_n-b_n})
\end{equation}
for some $a_j\in\{0,1,\ldots,k_j\}$ and $b_j\in\{0,1,\ldots,\ell_j\}$,
and where for all $j$, $a_j\ne0\Longleftrightarrow b_j\ne0$.
Writing $Q_t=p_tTp_t+p_tX(1-p_t)$ we see
\[
\tau(Q_t^*Q_t)=\tau(p_tT^*p_tTp_t)+\tau((1-p_t)Xp_tX(1-p_t)).
\]
From the obervation~(b) about $p_tTp_t$, we find $\tau(p_tT^*p_tTp_t)=t^2/2$.
Moreover, using freeness of $X$ and $p_t$, we find
$\tau((1-p_t)Xp_tX(1-p_t))=t-t^2$, so $\tau(Q_t^*Q_t)=t-(t^2/2)$ and $\|Q_t\|_2<t^{1/2}$.
On the other hand,
\[
\|Q_t^2\|_2=\|p_tTp_tT\|_2\le\|p_tTp_t\|\,\|Q_t\|_2
\le2t^{1/2}\|Q_t\|_2<2t.
\]
Consequently,
\[
\|Q_t^*Q_t^2\|_1\le\|Q_t^*\|_2\|Q_t^2\|_2<2t^{3/2}.
\]
Since $\|S_t\|,\,\|Q_t\|\le\|T\|\le2$ for all $t$,
if in~\eqref{eq:SQform} $a_j\ge1$ and $b_j\ge2$ for some $j$, then
\begin{align*}
|\tau((S_t^*)^{k_1-a_1}(Q_t^*)^{a_1}&Q_t^{b_1}S_t^{\ell_1-b_1}
\cdots(S_t^*)^{k_n-a_n}(Q_t^*)^{a_n}Q_t^{b_n}S_t^{\ell_n-b_n})|\le \\
&\le\|(S_t^*)^{k_1-a_1}(Q_t^*)^{a_1}Q_t^{b_1}S_t^{\ell_1-b_1}
\cdots(S_t^*)^{k_n-a_n}(Q_t^*)^{a_n}Q_t^{b_n}S_t^{\ell_n-b_n}\|_1 \\
&\le2^{2m-3}\|Q_t^*Q_t^2\|_1\le2^{2m-3}t^{3/2},
\end{align*}
and similarly if $a_j\ge2$ and $b_j\ge1$.
Therefore
\[
\tau((S_t^*)^{k_1-a_1}(Q_t^*)^{a_1}Q_t^{b_1}S_t^{\ell_1-b_1}
\cdots(S_t^*)^{k_n-a_n}(Q_t^*)^{a_n}Q_t^{b_n}S_t^{\ell_n-b_n})=O(t^{3/2})
\]
as $t\to0$, except possibly if for all $j$ either $a_j=0=b_j$ or $a_j=1=b_j$;
these are the terms we shall examine in more detail.
Each of them can be written in the form
\begin{equation}
\label{eq:FQform}
\tau(F_1(Q_t^*Q_t)^{d_1}F_2(Q_t^*Q_t)^{d_2}\cdots F_s(Q_t^*Q_t)^{d_s})
\end{equation}
for some $s\in\Nats$ and $d_1,\ldots,d_s\in\Nats$, where each $F_j$ is a monomial
in $S_t$ and $S_t^*$ with $F_j\ne1$,
or, when ``$s=0$,'' in the form $\tau((Q_t^*Q_t)^n)$.
Note that we always have $F_j=(1-p_t)F_j(1-p_t)$.
Let us show that for all $s\in\Nats$ and $d_1,\ldots,d_s\in\Nats$,
the asymptotic behavior
\begin{equation}
\label{eq:QQasympt}
\tau\big((1-p_t)(Q_t^*Q_t)^{d_1}(1-p_t)(Q_t^*Q_t)^{d_2}\cdots(1-p_t)(Q_t^*Q_t)^{d_s})
=t+O(t^2)
\end{equation}
holds as $t\to0$.
We have
\[
Q_t^*Q_t=(1-p_t)Xp_tX(1-p_t)+(1-p_t)Xp_tTp_t+p_tT^*p_tX(1-p_t)+p_tT^*p_tTp_t,
\]
so $(1-p_t)(Q_t^*Q_t)^{d_j}(1-p_t)$ is equal to $((1-p_t)Xp_tX(1-p_t))^{d_j}$
plus some other monomials in $p_t$, $(1-p_t)$, $X$, $T^*$ and $T$, each of
which contains at least one of the submonomials $p_tT^*p_tTp_t$ or $p_tTp_tT^*p_t$.
Since
\[
\|p_tT^*p_tTp_t\|_1=\|p_tTp_tT^*p_t\|_1=\tau(p_tT^*p_tTp_t)=t^2/2,
\]
we obtain the asymptotic behavior
\begin{align*}
\tau\big((1-p_t)(Q_t^*Q_t)^{d_1}(1-p_t)&(Q_t^*Q_t)^{d_2}\cdots(1-p_t)(Q_t^*Q_t)^{d_s}) \\
&-\tau\big(((1-p_t)Xp_tX(1-p_t))^{d_1+\cdots+d_s}\big)=O(t^2)
\end{align*}
as $t\to0$.
Now Lemma~\ref{lem:Xpasympt} shows that the asymptotic behavior~\eqref{eq:QQasympt} holds.
Similar reasoning shows that for every $n\in\Nats$, $\tau((Q_t^*Q_t)^n)=t+O(t^2)$ as $t\to0$.

Using observation~(b) above we have $Q_t\in W^*(p_tX\cup p_tB)$.
Using observation~(a) and Lemma~\ref{lem:1pXp}, we see that the pair
\begin{equation}
\label{eq:SQQfree}
\{S_t\},\quad \{(1-p_t)(Q_t^*Q_t)^k(1-p_t)\mid k\in\Nats\}
\end{equation}
is $*$--free with respect to $(1-t)^{-1}\tau\restrict_{(1-p_t)\Ac(1-p_t)}$.
Therefore, Speicher's moment--cummulant formula~\cite[Theorem 2.17]{Sp:Gr}
gives
\begin{multline*}
\tau(F_1(Q_t^*Q_t)^{d_1}F_2(Q_t^*Q_t)^{d_2}\cdots F_s(Q_t^*Q_t)^{d_s})=\sum_{\pi\in\NC(s)} \\
k_\pi[F_1,\ldots,F_s]
\tau_{K(\pi)}[(1-p_t)(Q_t^*Q_t)^{d_1}(1-p_t),\ldots,(1-p_t)(Q_t^*Q_t)^{d_s}(1-p_t)].
\end{multline*}
However, the asymptotic behavior~\eqref{eq:QQasympt} implies
\begin{multline*}
\tau_{K(\pi)}[(1-p_t)(Q_t^*Q_t)^{d_1}(1-p_t),\ldots,(1-p_t)(Q_t^*Q_t)^{d_s}(1-p_t)]= \\
=\begin{cases}
t+O(t^2)&\text{if }K(\pi)=\oneb_s, \\
O(t^2)&\text{otherwise,}
\end{cases}
\end{multline*}
where $\oneb_s$ is the (trivial) partition of $\{1,\ldots,s\}$ into only one block.
If $K(\pi)=\oneb_s$ then $\pi=\zerob_s$, which is the partition of $\{1,\ldots,s\}$ into $s$ blocks, and
$k_{\zerob_s}[F_1,\ldots,F_s]=\prod_{j=1}^s\tau(F_j)$.
Therefore,
\[
\tau(F_1(Q_t^*Q_t)^{d_1}F_2(Q_t^*Q_t)^{d_2}\cdots F_s(Q_t^*Q_t)^{d_s})=
\bigg(\prod_{j=1}^s\tau(F_j)\bigg)t+O(t^2)
\]
as $t\to0$.

Now let us put this together to obtain the formula~\eqref{eq:recur}.
Of the $2^{2m}-1$ terms in which $Q_t$ or its adjoint appears, after dividing by $t$, all contribute zero
as $t\to0$ except those of the form
\begin{equation}
\label{eq:Skcterm}
\tau\big((S_t^*)^{k_1-c_1}(Q_t^*Q_t)^{c_1}S_t^{\ell_1-c_1}(S_t^*)^{k_2-c_2}(Q_t^*Q_t)^{c_2}S_t^{\ell_2-c_2}
\cdots(S_t^*)^{k_n-c_n}(Q_t^*Q_t)^{c_n}S_t^{\ell_n-c_n}\big),
\end{equation}
where $c_1,\ldots,c_n\in\{0,1\}$.
We perform the summation $\sum_{r=1}^n\sum_{1\le j(1)<\cdots<j(r)\le n}$,
let $c_{j(1)}=c_{j(2)}=\cdots=c_{j(r)}=1$ and let all other $c_i=0$.
Then the term~\eqref{eq:Skcterm} becomes
\begin{align*}
\tau\big(&(S_t^*)^{k_1}S_t^{\ell_1}\cdots(S_t^*)^{k_{j(1)-1}}S_t^{\ell_{j(1)-1}}(S_t^*)^{(k_{j(1)})-1}(Q_t^*Q_t) \\
&S_t^{(\ell_{j(1)})-1}(S_t^*)^{k_{j(1)+1}}S_t^{\ell_{j(1)+1}}\cdots
(S_t^*)^{k_{j(2)-1}}S_t^{\ell_{j(2)-1}}(S_t^*)^{(k_{j(2)})-1}(Q_t^*Q_t) \\
&\qquad\qquad\vdots \\
&S_t^{(\ell_{j(r-1)})-1}(S_t^*)^{k_{j(r-1)+1}}S_t^{\ell_{j(r-1)+1}}\cdots
(S_t^*)^{k_{j(r)-1}}S_t^{\ell_{j(r)-1}}(S_t^*)^{(k_{j(r)})-1}(Q_t^*Q_t) \\
&S_t^{(\ell_{j(r)})-1}(S_t^*)^{k_{j(r)+1}}S_t^{\ell_{j(r)+1}}\cdots
(S_t^*)^{k_n}S_t^{\ell_n}\big)= \displaybreak[2] \\
=\tau\big(&
\begin{aligned}[t]
S_t^{(\ell_{j(r)})-1}(S_t^*)^{k_{j(r)+1}}S_t^{\ell_{j(r)+1}}\cdots
(S_t^*)^{k_n}S_t^{\ell_n}
&(S_t^*)^{k_1}S_t^{\ell_1}\cdots \\
&\cdots(S_t^*)^{k_{j(1)-1}}S_t^{\ell_{j(1)-1}}(S_t^*)^{(k_{j(1)})-1}(Q_t^*Q_t)
\end{aligned} \\
&S_t^{(\ell_{j(1)})-1}(S_t^*)^{k_{j(1)+1}}S_t^{\ell_{j(1)+1}}\cdots
(S_t^*)^{k_{j(2)-1}}S_t^{\ell_{j(2)-1}}(S_t^*)^{(k_{j(2)})-1}(Q_t^*Q_t) \\
&\qquad\qquad\vdots \\
&S_t^{(\ell_{j(r-1)})-1}(S_t^*)^{k_{j(r-1)+1}}S_t^{\ell_{j(r-1)+1}}\cdots
(S_t^*)^{k_{j(r)-1}}S_t^{\ell_{j(r)-1}}(S_t^*)^{(k_{j(r)})-1}(Q_t^*Q_t)\big).
\end{align*}
Divided by $(m+1)t$, this quantity tends in the limit as $t\to0$ to
\begin{align*}
\frac1{m+1}&M(\ell_{j(r)}-1,k_{j(r)+1},\ell_{j(r)+1},\ldots,k_n,\ell_n,
k_1,\ell_1,\ldots,k_{j(1)-1},\ell_{j(1)-1},k_{j(1)}-1)\cdot \\
\cdot&\prod_{i=1}^{r-1}M(\ell_{j(i)}-1,k_{j(i)+1},\ell_{j(i)+1},\ldots,k_{j(i+1)-1}\ell_{j(i+1)-1},k_{j(i+1)}-1) \\
=\frac1{m+1}&M(k_1,\ell_1,\ldots,k_{j(1)-1},\ell_{j(1)-1},k_{j(1)}-1,
\ell_{j(r)}-1,k_{j(r)+1},\ell_{j(r)+1},\ldots,k_n,\ell_n)\cdot \\
\cdot&\prod_{i=1}^{r-1}M(\ell_{j(i)}-1,k_{j(i)+1},\ell_{j(i)+1},\ldots,k_{j(i+1)-1}\ell_{j(i+1)-1},k_{j(i+1)}-1)\;.
\end{align*}
Thus we obtain~\eqref{eq:recur}.
\end{proof}

We finish this section and the paper with a conjecture.

\begin{conj}
\label{conj}
For all $k,n\in\Nats$,
\[
\tau\big(((T^*)^kT^k)^n\big)=\frac{n^{nk}}{(nk+1)!}.
\]
\end{conj}
The conjecture was proved in the case $k=1$ in Theorem~\ref{thm:TstT}.
It is easy to prove in the case $n=1$, either from the recursion formula~\eqref{eq:recur}
or using Lemma~\ref{lem:utmom},
and either of these techniques with a little more work
can be used to prove it
in the case $n=2$, (or see the proof by contour integration of
a generating function in~\cite{DY}).
We have, without making an effort to go very far,
checked some additional instances of the conjecture on a computer
using the recursion formula~\eqref{eq:recur}.
The following table lists the cases of the conjecture that
we have either proved or checked by computation:
\bigskip
\begin{center}
\begin{tabular}{||l||l||l||l||l||l||l||l||}
\hline
  $n=1$ &   $n=2$ &    $n=3$ &    $n=4$ &    $n=5$ &     $n=6$ &   $n=7$ & all $n$  \\ \hline
all $k$ & all $k$ & $k\le60$ & $k\le25$ & $k\le15$ &  $k\le10$ & $k\le5$ &   $k=1$. \\ \hline
\end{tabular}
\end{center}
\bigskip


\bibliographystyle{plain}

\end{document}